\documentclass{article}

\usepackage{amsmath, amssymb, anysize, datetime, chemarr, latexsym, stmaryrd, xspace}
\usepackage{stackrel, color}
\usepackage[usenames,dvipsnames]{xcolor}
\usepackage{amsfonts,epsfig,colordvi,graphics,textcomp}
\usepackage{amsthm}
\usepackage{bm}
\usepackage{longtable}
\usepackage{array}
\usepackage{verbatim}
\usepackage{afterpage}
\usepackage{pdfpages}
\usepackage[font=small,justification=justified]{caption}
\usepackage{bbm}
\usepackage{listings}
\usepackage{hyperref}
\usepackage[nopostdot,acronym,nomain,toc,nonumberlist]{glossaries}
\usepackage{cases}
\usepackage{titlesec}
\usepackage{chemfig}
\usepackage[font=small,justification=centering]{subcaption}
\usepackage{mathrsfs,mathtools}
\usepackage{enumerate}
\usepackage{makeidx}
\usepackage[ruled]{algorithm2e}
\usepackage[document]{ragged2e}
\usepackage[version=4]{mhchem}
\usepackage{physics}
\usepackage{booktabs}
\usepackage{array}
\usepackage{cite}

\newtheorem{thm}{Theorem}[section]
\newtheorem{lem}{Lemma}[section]

\newtheorem{cor}{Corollary}[section]
\newtheorem{rem}{Remark}[]

\newtheorem{assume}[subsection]{Assumption}
\usepackage{ragged2e}

\usepackage[letterpaper,top=2cm,bottom=2cm,left=3cm,right=3cm,marginparwidth=1.75cm]{geometry}

\title{Mathematical Modeling of a pH Swing Precipitation Process and its Optimal Design}
\author{Sandesh Athni Hiremath$^1$, Chinmay Hegde$^2$ and Andreas Voigt$^2$ \\ 
{Correspondence: \it sandesh.hiremath@rptu.de} \\
{$^1$\it Department of Mathematics},\\
{\it Rhineland-Palatinate Technical University of Kaiserslautern-Landau,}\\ 
{\it Kaiserslautern, Germany.} \\[.5ex]
{$^2$ \it Institute of Process Engineering},\\
{\it Otto von Guericke University, Magdeburg, Germany.}}
\newcommand{\psic}{\Psi^C}
\newcommand{\psiq}{\Psi^Q}
\newcommand{\psih}{\Psi^H}
\newcommand{\psif}{\Psi^F}
\newcommand{\R}{\mathbb{R}}
\newcommand{\E}{\mathbb{E}}
\newcommand{\N}{\mathbb{N}}
\newcommand{\cnt}{\mathcal{C}}
\newcommand{\dom}{\mathfrak{D}}

\newcommand{\F}{\mathcal{F}}
\newcommand{\W}{\mathcal{W}}
\newcommand{\Ham}{\mathcal{H}}
\newcommand{\x}{\mathbf{x}}
\newcommand{\Xalp}{X^{\alpha}}
\newcommand{\lam}{\lambda}
\newcommand{\blam}{{\bm{\lambda}}}
\newcommand{\const}{{\mbox{\large{\rm k}}}}

\newcommand{\vsig}{\scalebox{1.5}{$\varsigma$}}
\newcommand{\bvsig}{\scalebox{1.5}{$\bm{\varsigma}$}}

\def \bs {\boldsymbol}
\newcommand{\diag}{\rm{diag}}

\newcommand{\Om}{\Omega}
\newcommand{\om}{\omega}
\renewcommand{\P}{\mathbb{P}}

\newcommand{\indc}{\mathbbm{1}}

\newcommand{\mat}[1]{ 
	\begin{bmatrix} #1
        \end{bmatrix}
}
\newcommand{\smat}[1]{
	\left[\begin{smallmatrix} #1
	\end{smallmatrix}\right]
    }

\def \bs {\boldsymbol}

\newcommand{\caion}{{\rm Ca}^{2+}}
\newcommand{\coion}{{\rm CO}^{2-}_3}

\newcommand{\caco}{{\rm CaCO}_3}

\newcommand{\csat}{C_{{\rm sat}}}
\newcommand{\Ksp}{K_{\rm sp}}
\newcommand{\dt}{{\rm dt}}
\newcommand{\dWt}[1]{{{\rm d}W^{#1}_t}}
\newcommand{\dx}{{\rm dx}}
\newcommand{\dtau}{{\rm d}\tau}
\newcommand{\ds}{{\rm ds}}

\newcommand{\dr}{{\rm dr}}

\date{}

\begin{document}
\maketitle

\justifying

\begin{abstract}
In this work we consider the semi-batch process of precipitation of calcium carbonate solids from a solution containing calcium ions by adjusting the pH of the solution. The change in pH is induced either by the addition of alkaline solution such as sodium hydroxide (NaOH) or by the addition of a carbon dioxide gas (CO$_2$) to the given ionic solution. Under this setup we propose a system of degenerate stochastic partial differential equations that is able to explain the dynamical behavior of the key components of precipitation process. In particular, we propose a semi-linear advection equation for the dynamics of particle size distribution (PSD) of the precipitated particles. This is in turn coupled with a system of stochastic differential equations (SDEs) that is able to explain the chemical kinetics between calcium ions ($\caion$), calcium carbonate ($\caco$) in aqueous state and pH of the solution. The resulting coupled system is first mathematically studied, in particular conditions for the  existence of a mild-solution is established and also the long time behavior of the system is established. Following this we consider the validation of the model for which we consider experimentally obtained lab-scale data. Using this data we propose three methods to fit the model with the data which then also validates the suitability of the proposed model. The three methods include manual intuitive tuning, classical forward backward SDE (FBSDE) method and finally a DNN based method. The FBSDE method is based on the stochastic optimal control formulation for which we provide the necessary and sufficient condition for the existence of an optimal solution. Lastly, we compare the three methods and show that DNN method is the best in terms of the lowest error and as the most economical in terms of computational resources necessary during online use. 
\end{abstract}

\section{Introduction \label{sec:intro}}

Precipitation is a ubiquitous and highly significant unit operation across various industries, ranging from the production of high-value pharmaceuticals and specialty chemicals to large-scale commodity manufacturing and environmental engineering applications \cite{Feng2007Effect, Ma2002Optimal, Woodall2024}. The significance of precipitation can be primarily attributed to its ability to form solid particles with specific desired characteristics, such as size, morphology, and purity, from homogeneous solutions. Apart from typical manufacturing applications, precipitation plays a pivotal role in emerging sustainable technologies, notably in carbon capture and utilization (CCU) for mineral carbonation, where CO$_2$ is sequestered by reacting with alkaline materials to produce stable carbonates like calcium carbonate (CaCO$_3$) \cite{Seifritz1990CO2, Goff1998Environmental, Bobicki2012Carbon, Olajire2013Carbon, Sanna2014Carbonation, Hegde2024Towards}. Furthermore, the increasing demand for advanced materials and the need for sustainable industrial processes in turn necessitates an in depth understanding and a more precise modeling and control of precipitation phenomena. The quality of the precipitated product, be it a drug substance, a catalyst support, or a CO$_2$ sequestration material, is critically dependent on the intricate interplay of physical and chemical processes occurring during precipitation. Some examples of these include nucleation (the birth of new particles), crystal growth (enlargement of existing particles), aggregation (clumping of particles), and breakage (fragmentation of particles). The delicate balance and strong coupling between these mechanisms make precipitation processes inherently complex, challenging to model, and difficult to control to meet stringent product specifications.

The realization of high-quality precipitates with consistent characteristics is often hindered by several fundamental challenges inherent to the process. Precipitation involves a series of concurrently occurring and interacting phenomena such as nucleation, particle growth, aggregation, breakage, and polymorphism. Nucleation, the process of formation of new solid phases, can be primary (from solution) or secondary (due to attrition or breakage of existing crystals). Nucleation rates are highly sensitive to supersaturation and often occur rapidly, making them difficult to control \cite{Koutsoukos1984Precipitation}. Particle growth is the process by which existing solid particles or nuclei increase in size through the deposition of additional material from the supersaturated solution. Growth rates can be anisotropic, leading to diverse particle morphologies, and are affected by surface reactions and mass transfer limitations \cite{Eisenschmidt2015Face, Borchert2012Efficient}. Aggregation and breakage describe processes where particles can agglomerate to form larger clusters, or large particles can break into smaller ones due to shear forces, both of which significantly impact the final particle size distribution (PSD) and are often difficult to characterize and model accurately. Lastly, polymorphism occurs when many substances can crystallize into multiple solid forms with different particle structures. This mechanism is particularly crucial in pharmaceuticals, as different polymorphs can have varying solubility, bio-availability, and stability, necessitating precise modeling and control of the process \cite{Simone2017Systematic}. Beyond these intrinsic complexities, precipitation is extremely sensitive to operating parameters such as temperature, pH, reactant concentrations, stirring rate, and residence time. Small variations in these conditions can lead to significant changes in particle size, morphology, and purity \cite{Feng2007Effect, Rauscher2005Influence}. Furthermore, the presence of impurities, even at trace levels, can drastically alter nucleation and growth kinetics, influencing crystal habit and even promoting undesired polymorphic forms \cite{Falini2009Calcium, Park2008Effects, Davis2000Calcium}.

In this challenging context, the controlled selective precipitation of calcium carbonate (CaCO$_3$) is particularly vital when dealing with industrial waste streams and mine tailing solutions containing a variety of dissolved ions, such as calcium (\ce{Ca^2+}), magnesium (\ce{Mg^2+}), potassium (\ce{K^+}), sodium (\ce{Na^+}), and iron (\ce{Fe^2+}) \cite{Bobicki2012, Khosa2019, Gomes2025}. In such multicomponent systems, achieving selective precipitation of CaCO$_3$ becomes a critical challenge and a valuable capability, ensuring that calcium is efficiently removed as solid carbonates, leaving behind other dissolved ions that could interfere with downstream processes or product purity. The process control for the selective precipitation of CaCO$_3$ from solutions with different dissolved ions is influenced by several factors, including solubility, ion concentration, and reaction kinetics. The relatively low solubility of CaCO$_3$ in water and alkaline pH ranges enables its preferential precipitation. Controlling the pH in the alkaline range promotes the conversion of dissolved CO$_2$ into carbonate (CO$_3^{2-}$) and bicarbonate (HCO$_3^{-}$) ions, which facilitate CaCO$_3$ precipitation. However, very high pH can induce co-precipitation of other species like MgCO$_3$ or iron hydroxides. This necessitates optimal regulation of the pH as it directly influences the solubility of CaCO$_3$, its nucleation and crystal growth kinetics, thereby enhancing its selectivity and avoiding co-precipitation of other substances. Semi-batch operation enhances this control by allowing gradual addition of reactants, such as CO$_2$ gas or carbonates and pH regulating solutions, ensuring uniform mixing and avoiding localized supersaturation that can lead to undesired by-products. CO$_2$ dissolution in water plays a dual role, supplying carbonate ions while moderating the pH through carbonic acid formation. The concentration of ions in the solution also governs the precipitation process. High calcium levels increase the supersaturation of CaCO$_3$, promoting nucleation, while managing competing ion concentrations, such as \ce{Mg^2+} or \ce{Fe^2+}, prevents their interference in the precipitation pathway. This precise control of parameters ensures the selective and efficient precipitation of CaCO$_3$, essential for processes like water treatment, mineral recovery, and carbon sequestration.

Addressing these complexities and sensitivities necessitates robust predictive frameworks, and math based dynamical models serve as indispensable tools. Population Balance Equations (PBEs) are a primary technique for mechanistic modeling of particulate systems, describing the evolution of the number density function of particles with respect to their intrinsic properties (e.g., size, shape, composition) and time. A general form of the PBE for particle size ($x$) in a well-mixed system is: 
$${\partial_t n(t,x)} + {\partial_x (G ~ n(t,x))} = B(t,x) - D(t,x) + \text{terms for aggregation/breakage},$$
where $n(t,x)$ is the number density function, $G$ is the linear growth rate, $B(t,x)$ represents birth due to nucleation or breakage, and $D(t,x)$ represents death due to removal or aggregation \cite{Ramkrishna2014Population}. PBEs require constitutive models for nucleation rates, growth rates, and aggregation/breakage kernels, which are often empirically derived or based on thermodynamic/kinetic theories \cite{Koutsoukos1984Precipitation, Liendo2022Nucleation}. Solving PBEs is computationally intensive due to their integro-partial differential nature. Common numerical methods include method of moments, finite difference methods, and fixed pivot methods \cite{John2009Numerical, Durr2020Approximate}. Moment methods reduce the PBE to a set of ordinary differential equations (ODEs) for statistical moments (e.g., total number, mean size), but suffer from closure problems, requiring approximations for higher-order moments. To capture the influence of fluid dynamics and mixing on local supersaturation and particle interactions, PBEs are increasingly coupled with Computational Fluid Dynamics (CFD). These multiscale models provide a more realistic representation of industrial precipitators by resolving spatial variations in concentration, velocity, and turbulence, which in turn affect nucleation, growth, and aggregation rates \cite{Bartsch2019Multiscale}. Furthermore, models for particle or crystal shape evolution can be coupled, allowing for the prediction of two-dimensional or three-dimensional particle morphology distributions \cite{Borchert2012Efficient, Eisenschmidt2015Face}. Complementary to PBEs, detailed thermodynamic and kinetic models are essential. Thermodynamic models describe chemical equilibria, solubility limits, and speciation in complex electrolyte solutions, which are critical for determining the driving force for supersaturation and thus precipitation \cite{Chen1986Generalized, Qian2012Calculation}. Kinetic Models focus on the rates of elementary steps, such as dissolution, surface reaction, and transport, particularly relevant for mineral carbonation processes that involve solid-liquid-gas reactions with complex dissolution kinetics of raw materials \cite{Teir2007Dissolution, Teir2009Thermodynamic, Liendo2022Nucleation}. The increasing availability of process-specific data and advancements in machine learning have opened new avenues for modeling precipitation. Reduced-Order Models can be used to simplify complex high-fidelity models, creating computationally cheaper surrogates for real-time applications. Machine learning algorithms, particularly those based on sparse regression, can identify governing partial differential equations (PDEs) directly from data, offering insights into complex processes where mechanistic understanding is incomplete \cite{Rudy2017DataDriven}. In particular deep neural networks have demonstrated their ability to solve high-dimensional PDEs and stochastic PDEs (SPDEs), which are often encountered in population balance modeling. Approaches like Physics-Informed Neural Networks (PINNs) \cite{Raissi2018Deep} and methods leveraging Neural Ordinary Differential Equations (NeuralODEs) \cite{E2021Deep, Han2018Solving} are also effective for forward and inverse modeling of precipitation dynamics, circumventing the need for traditional numerical solvers for complex PBEs. These ideas were recently incorporated in \cite{HiremathECC2024, HiremathEscape2024} where DNNs were used both as predictors and controllers to autonomously control the residence time for single and multi-compound precipitation process. 

In this work we shall be primarily focused on modeling and optimal design of calcium carbonate precipitation process controlled by varying pH of the solution. Subsequently, the design of an optimal controller shall be studied in a separate independent work. For the sake of operational ease in laboratory setting, we shall restrict to the case of a semi-batch process where the ionic solution and process operational regime are primarily kept fixed, and only the control inputs, specifically pH-varying reactants, are gradually added, followed by uniform mixing. Under this pretext, we model and analyze the dynamical evolution of the particle size density (PSD) of carbonate particles under the influence of varying pH. Since the effective precipitation requires changing pH direction, we call this technique of precipitation as the pH swing (precipitation) process. To this end, we shall use the population balance modeling technique \cite{Ramkrishna2014Population} to model the PSD. The novel aspect of the model lies in the way the population balance equation (PBE) is coupled with the kinetic equations of the involved chemical reactions. Accordingly, the main chemical components of interest are: the pH of the solution (or proton concentration [H$^+$]), concentration of calcium ions ([Ca$^{2+}$]), and the concentration of calcium carbonate (CaCO$_3$) in aqueous state. Based on the availability of CO$_2$ in the solution and the pH, bicarbonate [HCO$_3^-$] and carbonate [CO$_3^{2-}$] ions are formed to react with the available calcium ions [Ca$^{2+}$] in the solution to form CaCO$_3$ in the aqueous phase, which will then precipitate out as a solid based on the pH-dependent saturation level of the solution. Thus, these chemical kinetics inherently influence the nucleation of new particles \cite{Koutsoukos1984Precipitation, Liendo2022Nucleation} as well as the growth of formed particles \cite{Eisenschmidt2015Face, Borchert2012Efficient}. Consequently, the growth rates and nucleation rates are intimately coupled with the reaction kinetic components. Due to the presence of impurities, due to the uncertain effects or the lack of a detailed (involving more than the above-mentioned reactants) chemical mass balance equation, and the presence of inherent noisy thermodynamical effects, the dynamics of the key chemical components are modeled using stochastic differential equations (SDEs). Consequently, coupling the PBE with the SDEs, even though PBE itself is deterministic, results in a degenerate stochastic partial differential equation (SPDE). Here degeneracy has to be understood in the sense that the noise coefficient for one of the state variable (i.e. PSD variable) is zero. Though, mathematically, it would be very interesting to also consider a noise term for the PBE, it would then go beyond the assumed stable experimental setup of a semi-batch process, thus not relevant in the current scope of investigation. Nevertheless, the way of coupling noisy chemical kinetics with the PBE via the growth and nucleation rates is quite meaningful, since the PBE is already a macroscale description whose dynamics are affected by the noisy microscopic dynamics, and the effects of the latter are only seen as probabilistic variations of the macro-variable, namely the particle size distribution (PSD).

The semi-batch setup also motivates the use of stochastic optimal control problem (SOCP) formulation which then enables for algorithmic parameter estimation and model validation. The algorithm for model validation is obtained by converting the SOCP to a system of forward-backward stochastic differential equations (FBSDEs) via the stochastic maximum principle. Based on this the successive forward backward sweeping method is used to iteratively reduce the discrepancy between model predictions and experimental data by suitably updating the model parameters. The setting of a semi-batch process, where pH is varied in a continuous manner, facilitates the effective application of FBSDEs since it provides the continuous input-output time series required for solving the inverse (control/design) problem. Following this approach but going beyond the traditional optimization techniques and resorting the use of deep neural network (DNN) based methods, \cite{Raissi2018Deep, E2021Deep, Han2018Solving}, we can further obtain generalized validation and estimation method that shall then pave way of autonomous control. 

Based on this rest of paper is organized as follows. In Section \ref{sec:proc_desc} we shall describe the semi-batch process and the involved process phenomenons we consider for the modeling. Subsequently, in Section \ref{sec:wellposedness} we shall provide a detailed analysis of the precipitation model where in we provide conditions for the existence of a unique solution for the proposed system of equations along with the sufficient conditions for stationarity. In Section \ref{sec:ocf} we shall provide sufficient conditions for the existence of optimal design/control functions that enables development of semi-automated and automated methods for model validation. Following this in Section \ref{sec:numerics} we provide numerical investigations where in simulation results are provided along with different methods namely- manual, semi-automated and automated ways to validate the model based on experimental data. In Section \ref{sec:sim_results} we shall provide the simulation and validation results which are eventually discussed in Section \ref{sec:discuss} followed by concluding remarks and outlook in Section \ref{sec:conlcusion}.

\section{Process description \label{sec:proc_desc}}

The process of CaCO$_3$ precipitation using atmospheric CO$_2$ into an ionic solution rich in calcium ion (Ca$^{2+}$) involves several physical and chemical steps, where the atmospheric CO$_2$ reacts with Ca$^{2+}$ ions in the solution to form CaCO$_3$. The detailed description of the process is as follows:

\subsection{Precipitation of CaCO$_3$}
Calcium carbonate (CaCO$_3$) precipitation involves the reaction between Ca$^{2+}$ ions and CO$_3^{2-}$ ions in a solution. These carbonate ions dissociate from dissolved CO$_2$ gas in the atmosphere:
\begin{equation}
\label{chem:caco}
\text{Ca}^{2+}(\text{aq}) + \text{CO}_3^{2-}(\text{aq}) \rightleftharpoons \text{CaCO}_3(\text{aq})
\end{equation}
As the concentration of CaCO$_3$ in the solution increases, saturation is reached, and solid carbonate precipitates. The amount of precipitate formed depends on the solubility product constant ($K_{\rm sp}$), which represents the equilibrium between the dissolved ions and the solid phase of CaCO$_3$ \cite{Qian2012Calculation}:
\begin{equation}
K_{\rm sp} = \left[\text{Ca}^{2+}\right]_{\text{sat}} \left[\text{CO}_3^{2-}\right]_{\text{sat}}
\end{equation}
Here, $[\text{Ca}^{2+}]_{\rm sat}$ and $[\text{CO}_3^{2-}]_{\rm sat}$ represent the saturation concentrations of calcium and carbonate ions in a solution at a given pH and temperature. The saturation state of a calcium carbonate ($\caco$) solution is determined by comparing the ion activity product with the solubility product constant, $K_{\rm sp}$. At 25$^{\circ}$C, the value of $K_{\rm sp}$ for water is 4.8$\times$10$^{-9}$ \cite{Larson1973}. Supersaturation occurs when the activities of calcium ($\caion$) and carbonate ($\coion$) ions exceed their equilibrium values, leading to potential precipitation. This deviation from equilibrium can be quantified using the solubility ratio, $s$, defined as:
\begin{equation} s = \sqrt{\frac{a(\caion) a(\coion)}{K_{\rm sp}}}\end{equation}
where $a(\caion)$ and $a(\coion)$ are the activities of calcium and carbonate ions, respectively. The supersaturation ratio, $\sigma$, is then expressed as $\sigma = s - 1$.  Changes in pH, temperature, or ion concentrations can alter the saturation state, inducing either supersaturation $(\sigma > 0)$ or undersaturation $(\sigma < 0)$. Our supersaturation definition, $\sigma = s-1$, provides a signed relative measure compared to the ratio $s$ defined in \cite{Liendo2022}.

From the above equation, it is evident that the amount of precipitate formed, for a given concentration of $[\text{Ca}^{2+}]$, depends on the concentration of carbonate ions ($[\text{CO}_3^{2-}]$). The role of pH in determining the concentration of $[\text{CO}_3^{2-}]$ ions, and thus the solubility and precipitation of CaCO$_3$, is crucial for understanding and thus modeling the process of precipitation. The equilibrium of carbonate ion in the calcium ion-rich solution is governed by the dissociation of carbonic acid (H$_2$CO$_3$) which in turn depends on the availability of CO$_2$ gas. The CO$_2$ gas, dissolved in the ion-rich solution, forms H$_2$CO$_3$, which dissociates further into bicarbonate (HCO$_3^-$) and carbonate (CO$_3^{2-}$) ions. This process is highly pH-dependent:
\begin{equation}
\text{CO}_2 + \text{H}_2\text{O} \overset{K_H}{\rightleftharpoons} \text{H}_2\text{CO}_3 \overset{K_{a_1}}{\rightleftharpoons} \text{H}^+ + \text{HCO}_3^- \overset{K_{a_2}}{\rightleftharpoons} 2\text{H}^+ + \text{CO}_3^{2-}
\end{equation}
where $K_H$ is the Henry constant for CO$_2$ dissolution in water, while $K_{a_1}$ and $K_{a_2}$ are the first and the second dissolution constants for H$_2$CO$_3$, respectively. 
From the above equation, the concentration of [$\text{CO}_3^{2-}$] (in equilibrium) can be expressed as:
\begin{equation}
\left[\text{CO}_3^{2-}\right] = \frac{C_T K_{a_1} K_{a_2}}{\left[\text{H}^+\right]^2 + K_{a_1} \left[\text{H}^+\right] + K_{a_1} K_{a_2}}
\end{equation}
where $C_T$ is the total dissolved carbon (H$_2$CO$_3$ + HCO$_3^-$ + CO$_3^{2-}$).
The pH dependency of CaCO$_3$ precipitation is evident from the carbonate system equilibria. As the concentration of $[\text{H}^+]$ ions increases (lower pH), the equilibrium shifts towards H$_2$CO$_3$ and HCO$_3^-$, reducing $[\text{CO}_3^{2-}]$ and dissolving CaCO$_3$. Conversely, as the pH increases (lower $[\text{H}^+]$ concentration), $[\text{CO}_3^{2-}]$ increases, allowing for greater CaCO$_3$ precipitation due to higher supersaturation ($\sigma$).

\subsection{Population Balance Equation}
In a semi-batch system, the general Population Balance Equation (PBE) can be expressed as \cite{Liendo2022}:
\begin{equation}
{\partial_t F} + \nabla \cdot \left(\mathbf{V}F\right) = \hat{N} - \hat{D},
\end{equation}
where $F$ is the number density of the CaCO$_3$ crystals (\#/m$^3$ of the solution), and $\hat{N}$ and $\hat{D}$ are the birth and death rates (\#/m$^3$s) of the crystals due to nucleation and dissolution of the CaCO$_3$ precipitates. In a well-mixed reactor, the PBE simplifies to:\
\begin{equation}
\label{eq:psdpde_pre}
{\partial_t F(t,x)} + a(t) {\partial_x F(t,x)} = \hat{N}(t,x) - \hat{D}(t,x),
\end{equation}
where $a$ is the growth rate (m/s), $N(t,x)$ is the growth term due to nucleation process while $D(t,x)$ is the death term due to disassociation process. Both the terms are modeled in a separable form with $N(t,x) = N(t) F(t,x)$ and $D(t,x) = D(t) F(t,x)$, where $N(t)$ and $D(t)$ are size independent nucleation and disassociation terms modeling respective processes.
\begin{equation}
\label{eq:psdpde}
{\partial_t F(t,x)} + a(t) {\partial_x F(t,x)} = N(t)F(t,x) - D(t)F(t,x),
\end{equation}
The pH control on the PBE is dictated by its dependency on supersaturation. The magnitude of $\sigma$ directly influences the precipitation dynamics of CaCO$_3$. The size-independent growth rate is generally assumed in the crystal growth kinetics study for the precipitation reaction. The pH dependency of the growth rate in the supersaturated solution of CaCO$_3$ is governed by:
\begin{align}
\label{eq:grate}
\begin{aligned}
a(t) &= k_g \tanh\left((\sigma(t))^p\right).
\end{aligned}
\end{align}
Considering the stoichiometric ratio of $2:1$ between [$\caion$] and [$\caco$] in the reaction \eqref{chem:caco}, we can express $s$ in terms of $C(t)$ and $C_{\rm sat}$ as 
\begin{align}
\begin{aligned}
s &= {C(t) \over C_{\rm sat}(t)}, \quad C_{\rm sat}(t) = \sqrt{[\coion](t) \over K_{\rm sp}},\quad a(t) = k_g \tanh\left(\left( {C(t) \over C_{\rm sat}}  - 1\right)^p\right),
\end{aligned}
\end{align}
where $C(t)$ is the actual concentration of CaCO$_3$ at a given time, $C_{\rm sat}(t)$ is the pH dependent (via $\coion$) saturation level of $\caco$ solution, $k_g$ is the growth constant, and $p$ is the exponential constant depending on the growth mechanism. While the growth exponential $p$ typically ranges from 1 to 2, we set $p = 2$ to ensure that the growth coefficient $a$ is always a positive real number.

Nucleation, the process of forming new solid particles, is driven by the energy barrier associated with forming a critical nucleus. Letting $\delta$ denote the constant representing the energy required to create the interface between the newly formed solid nucleus and the surrounding solution, the nucleation rate $N(t)$ is given by:
\begin{align}
\label{eq:nucrate}
\begin{aligned}
N(t) &= k_N \: \exp \left(-\delta \Big{/} \left({\sigma(t)}\right)^n \right) 
= k_N \: \exp \left(-\delta \Big{/}{\left(\left({C(t) \over C_{\rm sat}} - 1\right)\right)^n}\right).
\end{aligned}
\end{align}
In order to ensure that the nucleation rate goes to zero as $C$ approaches saturation level $\csat$ and is also well-behaved near this value, it is sufficient to keep $\sigma(t)$ as a positive real number. Consequently, we set the nucleation coefficient as $n=2$.

\subsection{Semi-Batch Process Modelling}
For the semi-batch experimental setup, the concentration modelling was implemented. The change in concentration of calcium ions and carbonate ions over time reflects the precipitation process. Due to the alternative addition of CO$_2$ and NaOH, the dilution effect needs to be considered in the process modelling.

The general mass balance equation inside the reactor is defined by:
\begin{equation}
\text{Accumulation} = \text{In} - \text{Out} + \text{Generation} - \text{Consumption},
\end{equation}
\begin{equation}
\frac{dn_A(t)}{dt} = \sum \dot{n}_{A,\ \text{in}} - \sum \dot{n}_{A,\ \text{out}} + \sum \dot{n}_{A,\ \text{gen}} - \sum \dot{n}_{A,\ \text{cons}},
\end{equation}
where the number of moles $n = C R$, $C$ is the concentration (mol/L), and $R$ is the volume (L). The mass balance equation for a species $A$ in terms of concentration is:
\begin{equation}
R(t)\frac{dC_A(t)}{dt} + C_A(t)\dot{R}_{\rm in} = \sum \const_{A,\text{in}}\dot{R}_{\rm in} - \sum \const_{A,\text{out}}\dot{R}_{\text{out}} \pm r_A(t)R(t),
\end{equation}
where $C_A(t)\dot{R}$ accounts for the dilution effect, $r_A(t)$ is the reaction rate for formation/consumption of species $A$, and $\dot{R}_{\rm in}$ is the inflow rate of NaOH.

Once CO$_2$ is dissolved and NaOH is added, the [Ca$^{2+}$] concentration decreases due to CaCO$_3$ precipitation. Since only NaOH is added, $[\text{Ca}^{2+}]_{\text{in}} = 0$ and $[\text{Ca}^{2+}]_{\text{out}} = 0$. The mass balance for the change in concentration of [Ca$^{2+}$] is:
\begin{equation}
R(t)\frac{d[\text{Ca}^{2+}]}{dt} + [\text{Ca}^{2+}]\dot{R}_{\rm in} = - {r}(t)R(t),
\end{equation}
Thus, the time dependent change in [Ca$^{2+}$] is:
\begin{equation}
\label{eq:caODE}
\frac{d[\text{Ca}^{2+}]}{dt} = -{r}(t) - [\text{Ca}^{2+}]\frac{\dot{R}_{\rm in}}{R(t)},
\end{equation}
where $r (t)$ is the pH dependent reaction rate of CaCO$_3$ formation, defined as:
\begin{align}
\begin{aligned}
r(t) &= \hat{r}(t) U_r(t), \quad \hat{r}(t) = k_c[\text{Ca}^{2+}][\text{CO}_3^{2-}].
\end{aligned}
\end{align}
The pH dependence of the reaction is governed by the function $U_r(t)$, which in turn is dependent on the rate of change, $U_H$, of pH (see \eqref{eq:pHODE}). 
Due to the addition of CO$_2$, the carbon content of the system increases making the solution more acidic due the formation of bicarbonate ions. As a consequence, due to lesser availability of carbonate ions, the solubility ratio of $\caco$ decreases which results in the emergence of free $\caion$ in the solution. 
\begin{equation}
\label{eq:ca2_dominant}
\begin{aligned}
&\ce{CO_2 (\text{g}) \leftrightarrow CO_2 (\text{aq} )} \\
&\ce{CO_2 (\text{aq}) + H_2O(\text{i}) \leftrightarrow HCO_3^-(\text{aq}) + H^+(\text{aq})} \\
&\ce{CaCO_3(\text{aq})} \rightarrow \ce{\ce{Ca^{2+}(\text{aq})} + CO_3^{2-}(\text{aq})} \\
&\ce{ CO_3^{2-}(\text{aq}) + H^+(\text{aq}) \rightarrow HCO_{3}^{-}(\text{aq})}
\end{aligned}
\end{equation}

\noindent Similarly, during the NaOH addition there are more hydroxyl ions [OH$^-$], so free protons are used up to form water molecules, thus the reaction towards bicarbonate ions is less dominant, thus resulting in more free carbonate ions [CO$_3^{2-}$] which can react with free calcium ions [Ca$^{+2}$] to from calcium carbonate in aqueous state.
\begin{equation}
\label{eq:caco3_dominant}
\begin{aligned}
&\ce{NaOH (\text{aq}) \leftrightarrow Na^+ (\text{aq}) + OH^-(\text{aq})} \\
&\ce{OH^-(\text{aq}) + H^+(\text{aq}) \leftrightarrow H_2O(\text{i})}\\
&\ce{HCO_{3}^{-}(\text{aq})\rightarrow CO_3^{2-}(\text{aq}) + H^+(\text{aq})}\\
&\ce{\ce{Ca^{2+}(\text{aq})} + CO_3^{2-}(\text{aq})\rightarrow \ce{CaCO_3(\text{aq})}}
\end{aligned}
\end{equation}
Thus, altogether we have that during NaOH addition free $\caion$ ions are consumed to form $\caco$ in aqueous state (i.e. dissolved in solution) while the opposite is more dominant during the addition of CO$_2$. This swing in the reaction direction is captured by the changing signs, from positive to negative, respectively, of the rate function $r$ depending on the pH, which is modeled via $U_r$. This phenomenon is depicted in the reaction below:
\begin{equation}
\label{chem:caco}
\text{Ca}^{2+}(\text{aq}) \overset{r(t)}{\longleftrightarrow} \text{CaCO}_3(\text{aq}) \overset{a(t) + N(t)}{\longrightarrow} \text{CaCO}_3(\text{solid})
\end{equation}


\noindent While \eqref{eq:caODE} covers the left half of the reaction \eqref{chem:caco}, right of the reaction indicates the conversion of $\caco$ from aqueous state to solid state i.e. the precipitation of $\caco$ via nucleation and growth. 
Consequently, the dynamical evolution of the concentration of CaCO$_3$, denoted by $C(t)$, is described according to the mass balance:
\begin{equation}
\label{eq:cacoODE}
\frac{d{C}(t)}{dt} =  r(t) - {C}(t)\frac{\dot{R}_{\rm in}}{R(t)} - \Big(\rho S(t)a(t)\Big)\frac{1}{R(t)} - \Big(\rho v_{\text{nuc}}N(t)\Big)\frac{1}{R(t)},
\end{equation}
where $r(t)$ is the rate of $\caco$ increase depending on the rate of change of pH, due to CO$_2$ or NaOH addition, $\rho$ is the molar density of CaCO$_3$ (mol/m$^3$), $S(t)$ is the second moment of the particle size distribution for the precipitated CaCO$_3$. The coefficients $a(t)$ and $N(t)$ are the growth and nucleation rates provided in \eqref{eq:grate} and \eqref{eq:nucrate} respectively.   

Lastly, the change in pH due to addition of a basic or acidic solution such as NaOH or [CO$_2$]$_{\rm aq}$ respectively is given by a simple ODE
\begin{equation}
\label{eq:pHODE}
    \dv{t} H(t) = k_H U_H
\end{equation}
where $H(t)$ denotes the pH of the solution, $k_H$ denotes the proportionality rate constant for change of pH, $U_H$ is the amount of increase or decrease of acidity (-pH) and serves as a mechanism to control the precipitation process.

\begin{figure}[!h]
\begin{center}
    \includegraphics[scale=.25]{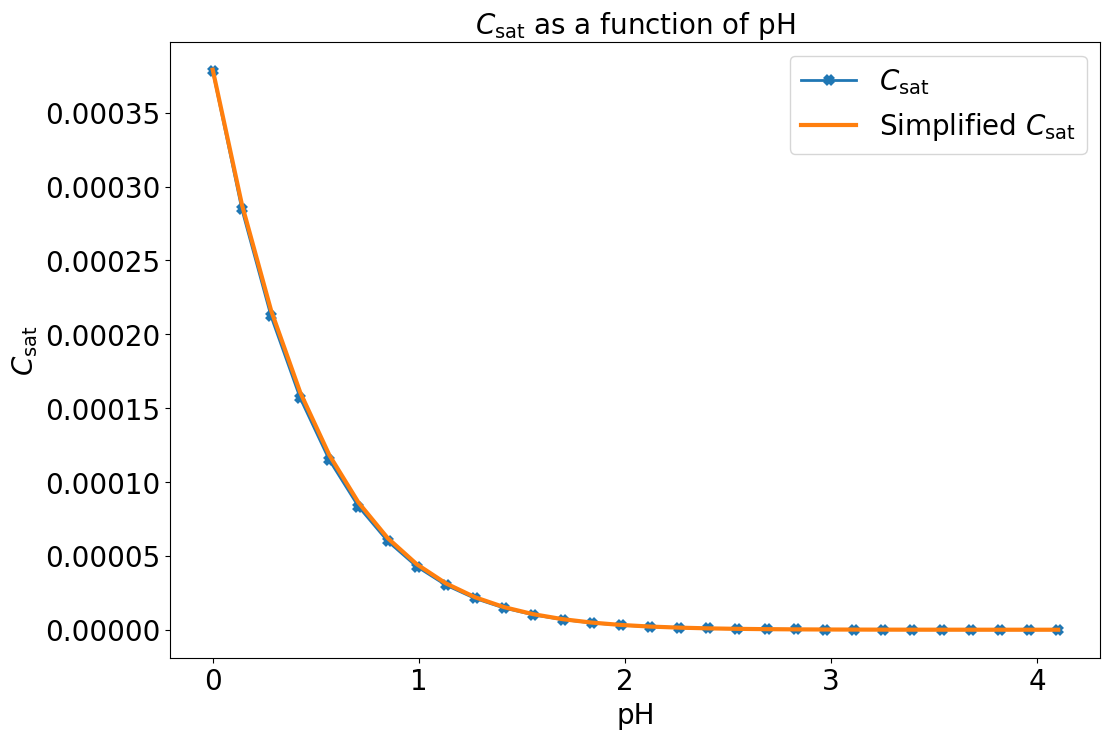} \qquad
    \includegraphics[scale=.25]{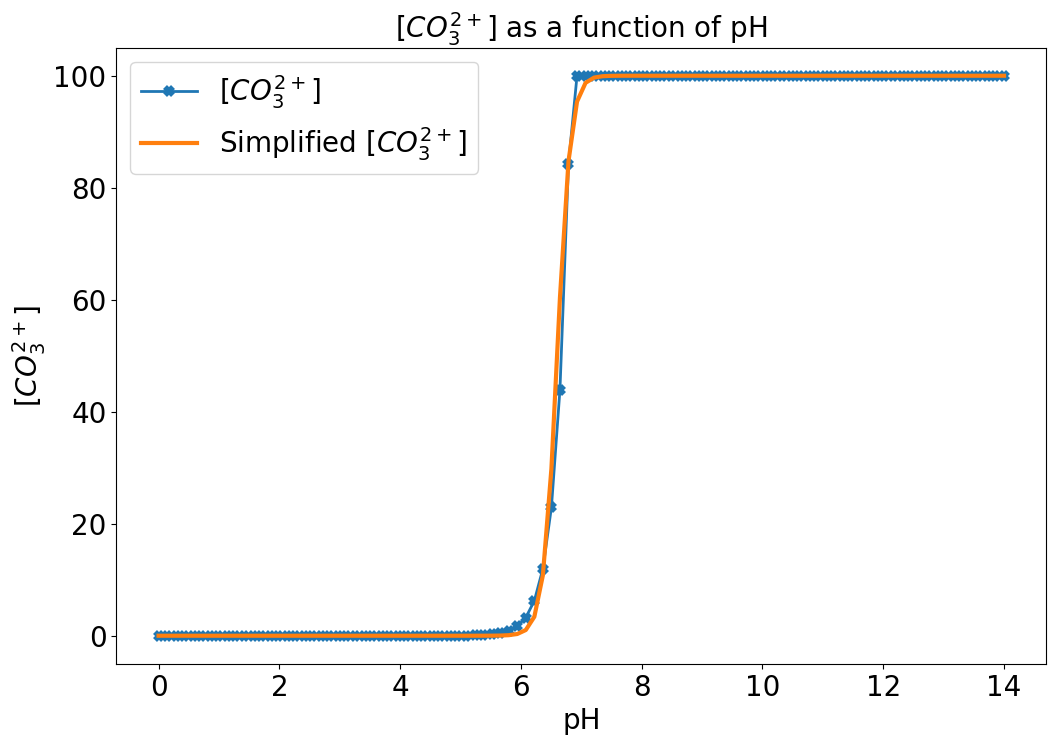} \\
\caption{\label{fig:simp_csat_co_func} Comparison between original $\csat$, $\text{CO}_3^{2-}$ functions and their simplified versions. Here blue stared curve shows the original functions as described in \eqref{eq:proc_coefs}, while the orange curve shows the respective simplified versions as mentioned in \eqref{eq:simp_csat} and \eqref{eq:simp_co3} respectively.}
\end{center}
\end{figure}

Altogether, combining the kinetic equations \eqref{eq:caODE}, \eqref{eq:cacoODE} and the particle density equation \eqref{eq:psdpde} we get the following system representing the precipitation process
\begin{align}
\label{eq:psd}
{\partial_t}  F(t,x) &+ a(t){\partial_x}  F(t,x) = N(t)F(t,x) - D(t) F(t,x), \quad  t, x > 0 \notag\\
\dv{t} C(t) &=  r(t) - C(t)\frac{\dot{R}_{\rm in}}{R(t)} - \Big(\rho S(t)a(t)\Big)\frac{1}{R(t)} - \Big(\rho v_{\text{nuc}}N(t)\Big)\frac{1}{R(t)}, \quad t > 0 \notag\\
\dv{t}[\text{Ca}^{2+}](t) &= -r(t) - [\text{Ca}^{2+}]\frac{\dot{R}_{\rm in}}{R(t)}, \quad t > 0 \\
\dv{t}H(t) &= k_H U_H,\quad \dot{R_{\rm in}} = k_v(t), \quad  t > 0 \notag\\
F(0,x) &= F_0(x), \quad C(0) = C_0, \quad R(0) = R_0 \notag\\
H(0) &= H_0, \quad [\text{Ca}^{2+}](0)=[\text{Ca}^{2+}]_0, \quad F(t,0) = 0. \notag \\[-5ex] \notag 
\end{align}
where, \\[-3ex]
\begin{equation}
\label{eq:proc_coefs}
\begin{aligned}
    a(C,H) &= k_g \tanh\left(\left(\frac{C(t)}{C_{\rm sat}} - 1\right)^{2}\right),  \quad D(t) = D(C(t)) := {k_d \left(\frac{C(t)}{C_{\rm sat}(t)}-1\right)^m \over 1 + \left(\frac{C}{C_{\rm sat}}\right)^m},\\
     S(F) &= \int_0^{\infty} \hspace*{-.25cm} x^2F(t,x)dx, \:\: \hat{r}(t) = K_{\rm sp} [\text{Ca}^{2+}] [\text{CO}_3^{2-}], \:\: r(t) = \hat{r}(t) U_r(t), \\
    N(t) &= N(C(t),H(t)) = k_N \: \exp \left(-\delta \Big{/}{\left(1 - \left({C(t) \over C_{\rm sat}(t)}\right)\right)^2}\right),\: \delta := \frac{16\pi\gamma^3_{in}v^2}{3K_B^3T^3}, \\
    [\text{CO}_3^{2-}](H(t)) &= \frac{C_T K_{a_1} K_{a_2}}{([\text{H}]^+(t))^2 + K_{a_1} ([\text{H}]^+(t)) + K_{a_1}K_{a_2}}, \quad C_{\rm sat} = \sqrt{K_{\rm sp}\over [\text{CO}_3^{2-}]},\\
    C_T(t) = C_T(H(t)) &= K_{\ce{[CO_2]}} \left( 1 + {K_{a_1} \over ([\text{H}]^+(t))} + {K_{a_2} K_{a_2} \over ([\text{H}]^+(t))^2} \right), \:\: [\text{H}]^+(t) := 10^{-H(t)}.
\end{aligned}
\end{equation}
Next we make the following simplifications. 
\begin{enumerate}
    \item In the current setting, we do not consider disassociation phenomenon of the formed precipitates thus we set $D(t,x) = 0$, by setting $k_d = 0$
    \item Since the mapping $H \mapsto [\text{CO}_3^{2-}](H)$ qualitatively resembles a sigmoid function, we shall replace it with a simpler function
\begin{align} 
\label{eq:simp_co3}
[\text{CO}_3^{2-}](H) = {\tilde K_{\ce{[CO_2]}} \over 1 + e^{-\tilde K_{a_1}\left((H-\tilde K_{a_2})(H + \tilde K_{a_3})\right)}}
\end{align}
The constants $\tilde K_{\ce{[CO_2]}} = 100, \tilde K_{a_1} = 0.45, \tilde K_{a_2} = 6.6, \tilde K_{a_3} = 13.5$ are chosen to fit the original function.
    \item Replace the unbounded function $H \mapsto C_{\rm sat}(H)$ with the following bounded function that qualitatively and quantitatively matches the original function away from infinity.
\begin{align} 
\label{eq:simp_csat}
C_{\rm sat}(H) = \frac{K_{1, \rm sat}}{\left(1+ K_{2, \rm sat} \sqrt{[\text{CO}_3^{2-}]/K_{\rm sp}}\right)}, \quad K_{1, \rm sat}:= .01, \quad K_{2, \rm sat}:= 1.2 \cdot 10^6. 
\end{align}
In Figure \ref{fig:simp_csat_co_func} below we see the comparison plot between the original and the simplified versions of the function.
\end{enumerate}

\noindent Let 
\begin{align}
\label{eq:redefs}
\begin{array}{lll}
    Q := [\text{Ca}^{2+}],& P := [\text{CO}_3^{2-}],& [H]^+ := 10^{-H},\\
    \tilde{k}_v(t) = {{k}_v(t) \over R(t)},&  \tilde \rho(t) = {\rho \over R(t)},& \tilde C = \max((C / C_{\rm sat} - 1),0)
\end{array}
\end{align} 
then we have the following equivalent reformulation of the above system:
\begin{align}
\label{eq:psd_trans}
\begin{aligned}
{\partial_t}  F(t,x) &+ a(t){\partial_x}  F(t,x) = N(t) F(t,x) , \quad  t, x > 0\\
\dv{t}C(t) = r(t) &- C(t) \tilde{k}_v(t)  - {\tilde \rho(t) a(t)}{S(t)} - {\tilde \rho(t) v_{\rm nuc}}N(t), \quad t > 0\\
\dv{t}Q(t) &= -r(t) -\tilde{k}_v(t) Q(t) ,  \quad  t > 0 \\
\dv{t}H(t) &= k_H U_H,\quad \dot{R} = k_v(t), \quad  t > 0 \\
F(0,x) &= F_0(x), \quad C(0) = C_0, \quad R(0) = R_0 \\
H(0) = &H_0, \quad Q(0)=Q_0, \quad F(t,0) = 0 \\
\end{aligned}
\end{align}
where, 
\begin{align}
\label{eq:precip_coefs}
\begin{aligned}
    \hspace*{-1cm} C_{\rm sat}(H) &=  \frac{K_{1, \rm sat}}{\left(1+ K_{2, \rm sat} \sqrt{P(H)/K_{\rm sp}}\right)}, \quad \tilde {C} = \max((C / \tilde C_{\rm sat} - 1),0),\\
    \hspace*{-.5cm} P(H) &= {\tilde K_{\ce{[CO_2]}} \over 1 + e^{-\tilde K_{a_1}\left((H-\tilde K_{a_2})(H + \tilde K_{a_3})\right)}}, \quad a(C,H) = k_g \tanh(\tilde{C}^2), \\
    N(C,H) &= k_N \exp \left(-\delta \Big{/}{\tilde{C}^2}\right),\quad {r}(t) = {\hat r}(t) U_r, \quad \hat{r}(t) = K_{\rm sp} Q P.
\end{aligned}
\end{align}
Due to the presence of solution impurities the kinetic reactions undergo uncertain dynamics. Such effects are modeled via independent white noise terms for the kinetic equations. Consequently, the pH controlled precipitation process is modeled via the following degenerate stochastic partial differential equation:
\begin{align}
\label{eq:noisy_psd}
{\partial_t}  F(t,x) &+ a(t){\partial_x}  F(t,x) = N(t) F(t,x) , \quad  t, x > 0 \notag\\[1ex]
\dd C(t) = \Big( r(t) &- C(t) \tilde{k}_v(t)  - {\tilde \rho(t) a(t)}{S(t)} - {\tilde \rho(t) v_{\rm nuc}} N(t) \Big)~\dt + \sigma_C C(t) ~ \dWt{1}, \quad t > 0 \notag\\[1ex]
\dd Q(t) &= \Big(-r(t) -\tilde{k}_v(t) Q(t) \Big) ~ \dt + \sigma_Q Q(t) ~ \dWt{2}, \quad  t > 0 \\[1ex]
\dd H(t) &= k_H U_H ~\dt + \sigma_H H(t) ~ \dWt{3}, \quad \dot{R} = k_v(t), \quad  t > 0 \notag\\[1ex]
F(0,x) &= F_0(x), \quad C(0) = C_0, \quad R(0) = R_0 \notag\\[1ex]
H(0) = &H_0, \quad Q(0)=Q_0, \quad F(t,0) = 0. \notag 
\end{align}
The involved coefficient terms are defined as in \eqref{eq:psd_trans} and $(\dWt{})_{t\ge0}$ with $\dWt{} = (\dWt{1}, \dWt{2}, \dWt{3}, \dWt{4})^{\top}$ is the white noise process corresponding to the Wiener process $(W_t)_{t\ge0}$ with $W_t = (W_t^1, W_t^2, W_t^3, W_t^4)^{\top}$.

In the following we shall first study the mathematical properties of the stochastic model followed by some numerical investigations. In the former we are mainly interested establishing the wellposedness and stability properties of the model while in the latter as we focus on the establishing an optimal design problem to quantitatively fit the model predictions with the observed data.

\section{Wellposedness Analysis \label{sec:wellposedness}}

We introduce the main function spaces and notations used throughout this work. Let $\dom \subset \mathbb{R}^n$ be open, $k \in \mathbb{N}_0$, and $1 \leq p \leq \infty$. The Sobolev space $\W^{k,p}(\dom)$ consists of functions $w \in L^p(\dom)$ whose weak derivatives up to order $k$ also belong to $L^p(\dom)$. For multi-indices $\alpha$ with $|\alpha| \leq k$, the norm is given by
\[
\|w\|_{\W^{k,p}(\dom)} = \left( \sum_{|\alpha| \leq k} \|D^\alpha w\|_{L^p(\dom)}^p \right)^{1/p}, \quad 1 \leq p < \infty,
\]
and for $p = \infty$,
\[
\|w\|_{\W^{k,\infty}(\dom)} = \max_{|\alpha| \leq k} \|D^\alpha w\|_{L^\infty(\dom)}.
\]
\noindent For $p = 2$, $\W^{k,2}(\dom)$ is a Hilbert space, denoted $H^k(\dom)$. The notation $\|w\|_{L^p(\dom)}$ may be abbreviated as $\|w\|_p$. The dual space $\W^{-k,q}(\dom)$ is defined by $1/p + 1/q = 1$.

\noindent For time-dependent spaces, we use $L^p_T(\dom) := L^p([0,T]; L^p(\dom))$. The space $\cnt^{k}_T(\R)$ consists of functions $w: [0,T] \to \R$ that are $k$-times continuously differentiable, and $\cnt^{k,\gamma}_T(\R)$ consists of functions whose $k$th derivative is $\gamma$-H\"older continuous. Similarly, $\cnt^{k,\gamma}_T(Z)$ denotes functions $w: [0,T] \to Z$ with $k$th time derivative $\gamma$-H\"older continuous in $Z$. The notation $\|w\|_{\cnt_T(Z)}$ and $\|w\|_{\cnt^{\gamma}_T(Z)}$ is used for the respective norms.

\noindent Let $(\Omega, \mathcal{F}, \P)$ be a probability space. $L^p(\Om)$ denotes the space of $p$-integrable random variables. Let $W = (W^1, W^2, W^3)$ be a vector-valued Wiener process with independent components, adapted to the filtration $(\mathcal{F}_t)_{t \ge 0}$. For all $0 \leq s \leq t$,
\[
\E(W_t | \mathcal{F}_0) = 0, \quad \E \big[ (W_t - W_s) (W_t - W_s)^{\top} | \mathcal{F}_s \big] = (t - s) I.
\]

\noindent Let $Z$ be a Banach space. For $T > 0$, $L_{\F}^p(Z)$ is the space of $Z$-valued, $\mathcal{F}_t$-adapted processes $(w_t)_{t \ge 0}$ with
\[
\|w\|_{L_{\F}^p} := \left( \int_0^T \E \big[ \|w(t, \cdot)\|_Z^p \big] \dt \right)^{1/p} < \infty.
\]
$L^p_{\F}(\cnt(Z))$ denotes the class of $Z$-valued, $\mathcal{F}_t$-adapted processes with continuous sample paths. The state space is $\mathbf{X} := \mathbf{X}_1 \times \mathbf{X}_1 \times \mathbf{X}_1 \times \mathbf{X}_2$, with $\mathbf{X}_1 := L^p_{\F}(\Om; \cnt([0,T];\R))$ and $\mathbf{X}_2 := L^p_{\F}(\Om; \cnt([0,T];W^{2,p}(\dom)))$. The Banach space $\mathcal{X}_p := \R^3 \times \W^{2,p}(\dom,\R)$ represents the state vector $X_t = (H_t, Q_t, C_t, F_t)^{\top}$. Its dual is $\mathcal{X}_q^{'} := \R^3 \times \W^{-2,q}(\dom,\R)$, with duality pairing $\langle \cdot, \cdot \rangle_{\mathcal{X}_p}$. The adjoint space is $\mathbf{Z} := L^{\infty}_{\F}(\Om; \cnt([0,T];\R^4))$. The admissible control set is $\mathcal{U} := \{ u \in L^{\infty}_{\F}(\Om;L^{\infty}([0,T];\R)): u \in [u_{\min}, u_{\max}] \}$.

\begin{assume}
\label{asm:coef_cond}
The coefficients defined in \eqref{eq:proc_coefs} satisfy the following regularity and boundedness properties. Specifically, the growth rate $a$ and nucleation term $N$ are bounded and weakly differentiable functions, the carbonate ion term $P$ and the saturation term $C_{\rm sat}$ are bounded and continuously differentiable, and the rate modulation functions $U_r$ and $U_H$ are uniformly bounded. These properties are summarized as follows:
\begin{subequations}
\label{eq:coef_props}
\begin{align}
    &a: \R^+ \times \R^+ \to [0,1],~ a \in \W^{1,\infty}\big(\R^+ \times \R^+;[0,1]\big),~~ \|a\|_{\W^{1,\infty}} < \const_a \label{eq:coef_growth_rate_bd}\\
    &N: \R^+ \times \R^+ \to [0,1],~ N \in \W^{1,\infty}\big(\R^+ \times \R^+; [0,1]\big),~~\|N\|_{\W^{1,\infty}} < \const_N \label{eq:coef_nuc_bd}\\
    &P: \R^+ \to [0, \const_P],~ P \in \cnt^1\big(\R^+;[0, \const_P]\big),~~\|P\|_{\cnt^1} < \const_P,~ \const_P \gg 0 \label{eq:coef_carb_bd}\\
    &C_{\rm sat}:\R^+ \to \R^+,~ C_{\rm sat} \in \cnt^1\big(\R^+;\R^+\big),~~\|C_{\rm sat}\|_{\cnt^1} < \const_{C_{\rm sat}} \label{eq:coef_cotwo_bd}\\
    &U_r: [0,T] \to [u_{\rm min},u_{\rm max}],~ U_r \in \cnt^1 \Big([0,T];[u_{\rm min},u_{\rm max}]\Big),~~ \|U_r\|_{\cnt^1} < \const_{U_r} \label{eq:coef_cont_mod_bd}\\
    &U_H: [0,T] \to [u_{h_{\rm min}},u_{h_{\rm max}}],~ U_H \in L^{\infty}\Big([0,T];[u_{h_{\rm min}},u_{h_{\rm max}}]\Big),~~ \|U_H\|_{L^{\infty}} < \const_{U_H} \label{eq:coef_cont_mod_bd2}
\end{align}
\end{subequations}
\end{assume}

\noindent Based on the above assumptions, we have the following Lipschitz continuity of $a$ and $N$.

\begin{lem}\label{lem:growth_lip} Let $a$ and $N$ be defined as in \eqref{eq:precip_coefs}, then it holds that
\begin{align*}
    |a(C_1, H_1) - a(C_2, H_2)|^{2p} \le L_{a_C} |C_1 - C_2|^p + L_{a_H}|H_1 - H_2|^p \\
    |N(C_1, H_1) - N(C_2, H_2)|^{2p} \le L_{N_C} |C_1 - C_2|^p + L_{N_H} |H_1 - H_2|^p
\end{align*}
where $L_{a_C} , L_{N_C}, L_{a_H}, L_{N_H} \in (0, \infty)$ are fixed finite constants.
\end{lem}
\begin{proof}
Since and $a, N \in L^{\infty}(\R)$, we observe that
\begin{align*}
    |a(C_1, H_1) - a(C_2, H_2)|^{2p} \le 2\|a\|^p_{\infty} |a(C_1, H_1) - a(C_2, H_2)|^p \\
    |N(C_1, H_1) - N(C_2, H_2)|^{2p} \le 2\|N\|^p_{\infty} |N(C_1, H_1) - N(C_2, H_2)|^p
\end{align*}
Next, let $\hat C = (C/\csat -1)$, $\hat C^+ = \max(\hat C, 0),$  and from \eqref{eq:precip_coefs} we have 
$$P(H) = {\tilde K_{\ce{[CO_2]}} \over 1 + e^{\tilde K_{a_1}(H-\tilde K_{a_2})}}, \quad C_{\rm sat}(H) = \frac{K_{\rm sat}}{\left(1+ K_{\rm sat} \sqrt{P(H)/K_{\rm sp}}\right)}$$
From this we have that
\begin{align*}
    |a(C_1, H_1) - a(C_2, H_2)| &= k_g \left|\tanh((\hat C_1^+)) - \tanh((\hat C_2^+))\right| \le k_g | \tanh^{\prime}| ~ |\hat C_1^+ - \hat C_2^+ |\\
    &\le  \left|\underbrace{k_g  \sech^2(\hat C^+)}_{\text{Term 1}} \right| ~ \left|\hat C_1 - \hat C_2 \right| \\
    |N(C_1, H_1) - N(C_2, H_2)| &= k_N\left|\exp(-\delta/(\hat C_1^+)) - \exp(-\delta/(\hat C_2^+))\right| \\
    &\le  k_N \left|\delta \exp(-\delta/(\hat C^+))/ (\hat C^+)^3\right| ~ |\hat C_1^+ - \hat C_2^+ |\\
    &\le \left|\underbrace{k_N \delta \exp(-\delta/(\hat C^+))/ (\hat C^+)^3}_{\text{Term 2}}\right| ~  \left|\hat C_1 - \hat C_2 \right|
\end{align*}
and
\begin{align*}
    \left|\hat C_1 - \hat C_2 \right| &\le |C_1 - C_2| +  |\partial_H \csat(H)| |H_1 - H_2|
\end{align*}
where, 
\begin{align*}
\frac{\partial \csat}{\partial H} &= \underbrace{-\frac{K_{sat}^2}{2 \sqrt{K_{sp} P(H)}} \left(1 + K_{sat} \sqrt{P(H)/K_{ip}}\right)^{-2}}_{\text{Term3}} \frac{\partial P}{\partial H}, \quad 
\frac{\partial P}{\partial H} = -\frac{\tilde{K}_{(CO_2)} \tilde{K}_1 e^{\tilde{K}_1 (H - \tilde{K}_2)}}{\left(1 + e^{\tilde{K}_1 (H - \tilde{K}_2)}\right)^2}.
\end{align*}
Since 
$$ \|\text{Term1}\|_{\infty} < \infty, \:\: \|\text{Term2}\|_{\infty} < \infty,  \:\: \|\text{Term3}\|_{\infty} < \infty, \:\: \text{ and } \left \|\frac{\partial P}{\partial H} \right\|_{\infty}$$
we have that 
\begin{align}
\label{eq:growth_nuc_lip}
\begin{aligned}
    |a(C_1, H_1) - a(C_2, H_2)|^p \le \ell^p_{a_C} |C_1 - C_2|^p + \ell^p_{a_H} |H_1 - H_2|^p \\
    |N(C_1, H_1) - N(C_2, H_2)|^p \le \ell^p_{N_C} |C_1 - C_2|^p + \ell^p_{N_H} |H_1 - H_2|^p
\end{aligned}
\end{align}
where,
$$ \ell_{a_C} := \|\text{Term1}\|_{\infty}, \quad \ell_{a_H} :=  \|\text{Term1}\|_{\infty} \left\|\frac{\partial \csat}{\partial H} \right\|_{\infty} $$
$$ \ell_{a_C} := \|\text{Term2}\|_{\infty},\quad \ell_{a_H} :=  \|\text{Term2}\|_{\infty} \left\|\frac{\partial \csat}{\partial H}\right\|_{\infty} $$
Thus, altogether we get that
\begin{align*}
    \Big|a(C_1, H_1) - a(C_2, H_2)\Big|^{2p} &\le 2\|a\|^p_{\infty} \Big|a(C_1, H_1) - a(C_2, H_2)\Big|^p \\
    &\le 2\|a\|^p_{\infty} \Big(\ell^p_{a_C} |C_1 - C_2|^p + \ell^p_{a_H} |H_1 - H_2|^p \Big)\\
    &= L_{a_C} |C_1 - C_2|^p + L_{a_H} |H_1 - H_2|^p \\
    \Big|N(C_1, H_1) - N(C_2, H_2)\Big|^{2p} &\le 2\|N\|^p_{\infty} \Big|N(C_1, H_1) - N(C_2, H_2)\Big|^p \\
    &\le 2\|N\|^p_{\infty} \Big(\ell^p_{N_C} |C_1 - C_2|^p + \ell^p_{N_H} |H_1 - H_2|^p\Big) \\
    &= L_{N_C} |C_1 - C_2|^p + L_{N_H}|H_1 - H_2|^p
\end{align*}
where, $L_{a_C} := 2\|a\|^p_{\infty} \ell^p_{a_C},\ L_{a_H} :=  2\|a\|^p_{\infty} \ell^p_{a_H},\ L_{N_C} := 2\|N\|^p_{\infty} \ell^p_{N_C},\ L_{N_H} :=  2\|N\|^p_{\infty} \ell^p_{N_H}$.
\end{proof}

\subsection{A priori analysis}
\begin{lem}
\label{lem:F_mom_rec}
Let $F(t,x)$ be as per $\eqref{eq:psd}$, let $\Upsilon(t) = \int_{\dom} F(t,x) \dx$ with $\Upsilon_0 := \Upsilon(0) = \int_{\dom} F_0(x) \dx$, then 
\begin{align}
    \Upsilon(t) \le e^{\const_{N} t} \Upsilon_0.
\end{align}
\noindent Similarly, let $\Upsilon^n(t) := \int_{\dom} x^n F(t,x) \dx$ for $n \in \N_0$, with $\Upsilon^0(t) := \Upsilon(t)$, then 
\begin{align}
    \Upsilon^n(t) &= e^{\int_0^t N(s) \ds} \Upsilon^n_0 + n \int_0^t e^{\int_s^t N(r) \dr} a(s) \Upsilon_{n-1}(s) \ds.
\end{align}
Furthermore, defining $ A_{k}(t) = \int_0^t a(s) A_{k-1}(s)~\ds$, $A_1 (t) = \int_0^t a(s)~\ds$ and $A_0 = 1$ with the convention $\Upsilon^0_0 := \Upsilon^0(0) := \Upsilon_0$, it holds that
\begin{align}
\label{eq:n_moment_est}
    \Upsilon^n(t)&= e^{\int_0^t N(s) \ds} \sum_{k=0}^n {n! \over k!} \Upsilon^k_0 A_{n-k}(t).
\end{align}
\end{lem}
\begin{proof}
For $\Upsilon^0(t) = \int_{\dom} F(t,x) \dx = \int_{\R} \mathbbm{1}_{\dom} F(t,x) \dx = \int_{\R} \tilde F(t,x) \dx,$ we have
\begin{equation}
\label{eq:mean_psd}
\begin{aligned}
    \dv{t}\Upsilon^0(t) &= \int_{\R} \left[-a(t)\dv{x} \tilde F(t,x) + N(t)\tilde F(t,x) \right] \dx \\
    &= \left. -a(t) \tilde F(t,x) \right|_{\R}  + N(t) \Upsilon(t) \\
    &= N(t) \Upsilon^0(t) \\ 
\implies \Upsilon^0(t) &= e^{\int_0^t N(s) \ds} \Upsilon^0_0.
\end{aligned}
\end{equation}
Thus, for $N \in L^{\infty}([0,\infty))$, with $\|N\|_{\infty} \le \const_{N}$, we obtain that 
$$\Upsilon(t) \le e^{\const_{N} t} \Upsilon_0.$$
For $n\in \N$ and $\Upsilon^n(t) = \int_{\dom} x^n F(t,x) \dx$, we have
\begin{equation}
\label{eq:moment_psd}
\begin{aligned}
\dv{t}\Upsilon^n(t) &= \int_{\R} x^n \dv{t} \tilde F(t,x) \dx \\
&= \int_{\R} x^n \left[-a(t)\dv{x} \tilde F(t,x) + N(t) \tilde F(t,x) \right] \dx \\
&= -a(t) \int_{\R} x^n \dv{x} \tilde F(t,x) \dx + N(t) \int_{\dom} x^n F(t,x) \dx  \\
&= -a(t) \left[ \left. x^n \tilde F(t,x) \right|_{\R} - n \int_{\dom} x^{n-1} F(t,x) \dx \right] + N(t) \Upsilon^n(t) \\
&= -a(t) \left[ 0 - 0 - n \int_{\dom} x^{n-1} F(t,x) \dx \right] + N(t) \Upsilon^n(t) \\
&= n~a(t)~\Upsilon_{n-1}(t) + N(t)~\Upsilon^n(t) \\
\implies \Upsilon^n(t) &= e^{\int_0^t N(s) \ds} \Upsilon^n_0 + n \int_0^t e^{\int_s^t N(r) \dr} a(s) \Upsilon_{n-1}(s) \ds.
\end{aligned}
\end{equation}
Now, let $\Upsilon^n_0 \le \const_{\Upsilon}$ for all $n \in \N_0$, then using \eqref{eq:mean_psd} in \eqref{eq:moment_psd} we get
\begin{align*}
    \Upsilon_1(t) &= e^{\int_0^t N(s) \ds} \Upsilon^1_0 + n \int_0^t e^{\int_s^t N(r) \dr} a(s) \Upsilon(s) \ds \\
    &= e^{\int_0^t N(s) \ds} \Upsilon^1_0 + \int_0^t e^{\int_s^t N(r) \dr} a(s) e^{\int_0^s N(r) \dr} \Upsilon_0 \ds \\
    &= e^{\int_0^t N(s) \ds} \Upsilon^1_0 + e^{\int_0^t N(r) \dr} \Upsilon_0 \int_0^t a(s) \ds \\
    &= e^{\int_0^t N(s) \ds} \Upsilon^1_0 + e^{\int_0^t N(r) \dr} \Upsilon_0 A_1(t).
\end{align*}
Arguing by induction, let
\begin{align*}
    \Upsilon^n(t)&= e^{\int_0^t N(s) \ds} \sum_{k=0}^n {n! \over k!} \Upsilon^k_0 A_{n-k}(t)
\end{align*}
Then,
\begin{align*}
    \Upsilon^{n+1}(t)&= e^{\int_0^t N(s) \ds} \Upsilon^{n+1}_0 + (n+1) \int_0^t e^{\int_s^t N(r) \dr} a(s) \Upsilon^n(s) \ds \\
    &= e^{\int_0^t N(s) \ds} \Upsilon^{n+1}_0 + (n+1) \int_0^t e^{\int_s^t N(r) \dr} a(s)  e^{\int_0^s N(r) \dr} \sum_{k=0}^n {n! \over k!} \Upsilon^k_0 A_{n-k}(s) \ds \\
    &= e^{\int_0^t N(s) \ds} \Upsilon^{n+1}_0 + (n+1) e^{\int_0^t N(s) \ds} \int_0^t  a(s)  \sum_{k=0}^n {n! \over k!} \Upsilon^k_0 A_{n-k}(s) \ds \\
    &= e^{\int_0^t N(s) \ds} \Upsilon^{n+1}_0 + (n+1) e^{\int_0^t N(s) \ds} \sum_{k=0}^n {n! \over k!} \Upsilon^k_0 \int_0^t  a(s) A_{n-k}(s) \ds \\
    &= e^{\int_0^t N(s) \ds} \Upsilon^{n+1}_0 + e^{\int_0^t N(s) \ds} \sum_{k=0}^n {(n+1)! \over k!} \Upsilon^k_0 A_{n+1-k} \\
    &= e^{\int_0^t N(s) \ds} \sum_{k=0}^{n+1} {(n+1)! \over k!} \Upsilon^k_0 A_{n+1-k}(t).
\end{align*}
 Thus, we obtain the following closed form formula
\begin{align*}
    \Upsilon^n(t)&= e^{\int_0^t N(s) \ds} \sum_{k=0}^n {n! \over k!} \Upsilon^k_0 A_{n-k}(t), \forall n \in \N_0
\end{align*}
where, $ A_{k}(t) = \int_0^t a(s)A_{k-1}(s)~\ds$, $A_1 (t) = \int_0^t a(s)~\ds$, $A_0 = 1$ with the convention $\Upsilon^0_{0} := \Upsilon_0$.
\end{proof}
\begin{lem}
\label{lem:F_mom_bd}
Let $\Upsilon^n_0 \le \const_{\Upsilon^0_0} \const^{n}_{\Upsilon}$, then we have that  
$$\Upsilon^n(t) \le \const_{\Upsilon^0_0} e^{\const_{N} t} (\const_{a} \const_{\Upsilon} + t)^n.$$ 
\end{lem}
\begin{proof} By convention, we have that $A_0 = 1$ and $A_k(t) = \int_0^t a(s) A_{k-1}(s) \ds.$ Thus, we obtain that 
    \begin{align*}
        A_1 &= \int_0^t a(s) \ds \le \const_{a} t, \quad 
        A_2 = \int_0^t a(s) A_1(s) \ds \le \const^2_{a} {t^2 \over 2} 
    \end{align*}
and by mathematical induction
    \begin{align*}
        A_k &= \int_0^t A_{k-1}(s) \ds \le \const^k_{a} {t^k \over k!}. 
    \end{align*}
Substituting in \eqref{eq:n_moment_est} we get
\begin{align}
    \Upsilon^n(t) &= e^{\int_0^t N(s) \ds} \sum_{k=0}^n {n! \over k!} \Upsilon^k_0 A_{n-k}(t)  \le  e^{\int_0^t N(s) \ds} \const_{\Upsilon^0_0} \sum_{k=0}^n {n! \over k!} \const^k_{a} \const^k_{\Upsilon} {t^{n-k} \over (n-k)!} \notag\\
    &= e^{\int_0^t N(s) \ds} \const_{\Upsilon^0_0} \sum_{k=0}^n {n! \over (n-k)! k!} \const^k_{a} \const^k_{\Upsilon} t^{n-k}  
    = e^{\int_0^t N(s) \ds} \const_{\Upsilon^0_0} (\const_{a} \const_{\Upsilon} + t)^n \notag \\
    &\le \const_{\Upsilon^0_0} e^{\const_{N} t} (\const_{a} \const_{\Upsilon} + t)^n.
\end{align}
\end{proof}
\begin{rem} For $0 \le F_0 \in L^1(\R)$ with compact support $\dom \subset \R$, we see that
    \begin{align*} 
        \Upsilon^n_0 &:= \int_{\dom} x^n F_0(x) dx \le |\dom|^n \int_{\dom} F_0(x) dx \le \const_{F_0} |\dom|^n.
    \end{align*}
    Thus, the assumptions of Lemma \eqref{lem:F_mom_bd} are satisfied with $\const_{\Upsilon} := |\dom|$ and $\const_{\Upsilon^0_0} := \const_{F_0}.$
\end{rem}

\begin{lem}
\label{lem:mom_lin_op}
For $F \in L^{\infty}\left(\Om; \cnt([0,T];\W^{1,p}(\dom;\R))\right)$ let the mapping $F \mapsto S^F$ be defined as
\begin{align}
\label{eq:F_mom}
S^F := S[F] := \int_{\dom} x^2 F(\cdot,x) \dx. 
\end{align}
Then the following holds:
\begin{enumerate}
    \item For each $t \in [0,T]$, $S^F_t := S^{F(t,\cdot)} := S[F(t,\cdot)] \in \W^{-1}(\dom;\R)$ for almost all $\om \in \Om$.
    \item The mapping $t \mapsto S^F_t$, $S^F:[0,T] \to \W^{-1}(\dom;\R)$ is an element of $\cnt([0,T]; \W^{-1}(\dom;\R))$ for almost all $\om \in \Om$.
    \item $S \in L^{\infty}\left(\Om;\cnt([0,T];\W^{-1}(\dom;\R))\right)$
\end{enumerate}
\end{lem}
\begin{proof}
    For $F \in  L^{\infty}\left(\Om; \cnt([0,T];\W^{1,p}(\dom;\R))\right)$, since $F(\om,t) \in \W(\dom;\R)$ for each $t \in [0,T]$ and almost all $\om \in \Om$ and since 
    $$ \left | \int_{\dom} x^2 F(t,x) \dx \right | \le \|F(t,\cdot)\|_{\cnt} |\dom|^3 \le \|F(t,\cdot)\|_{\W^{1,p}} |\dom|^3$$
    we have that $\|S_t\|_{\mathcal{L}(\W^{1,p})} \le |\dom|^3$ for each fixed $t \in [0,T]$. 
    Now from linearity of $F \mapsto S^F$ we get that $S_t \in \W^{-1}(\dom;\R)$ for every  $t \in [0,T]$.
    Accordingly, since 
    $$ \sup_{t \in [0,T]} \|F(t,\cdot)\|_{\W^{1,p}} \le \const_F < \infty $$ 
    it implies that 
    $$ \sup_{t \in [0,T]} \|S_t\|_{\mathcal{L}(\W^{1,p})} \le |\dom|^3.$$
    Thus, $S \in \cnt([0,T];\W^{-1}(\dom;\R))$ for almost all $\om \in \Om$. Consequently, we also get that $$S \in L^{\infty}\left(\Om;\cnt([0,T];\W^{-1}(\dom;\R))\right).$$
\end{proof}
\subsection{Analysis of PSD dynamics}
Next we shall show the existence of a solution for the dynamics of particle size distribution. 

\begin{lem}
\label{lem:iso_psd_sol}
Let $a$, $N$ and $\csat$ be given functions of time as defined in \eqref{eq:precip_coefs}. The family of mappings $(\Psi^a_{s,t})_{0 \le s\le t \le T}$ and $(\Psi^{N}_{s,t})_{0 \le s\le t \le T}$ defined as:
\begin{align}
\begin{aligned}
\Psi^a_{s,t} g = \indc_{\dom} g \left(\cdot - \int_s^t a(\tau) \dtau \right), \quad \forall g \in L^p(\dom),\: 0 \le s \le t \le T, \\
\Psi^{N}_{s,t} g = e^{\int_s^t N(\tau) \dtau} g, \quad \forall g \in L^p(\dom),\: 0 \le s \le t \le T \\
\end{aligned}
\end{align}
forms a two-parameter evolution semigroup on $L^p(\dom)$. Based on this, the solution to the equation
\begin{align*}
\label{eq:psd_iso}
{\partial_t}  F(t,x) + a(t){\partial_x}  F(t,x) &= N(t) F(t,x) \\
F(0,x) &= F_0(x)
\end{align*}
for every $F_0 \in \W^{m,p}(\dom)$, $m, p \ge 1$, is given as
\begin{align}
\begin{aligned}
    F(t,x) = \Psi^{N}_{0,t} \Psi^a_{0,t} F_0(x) = e^{\int_0^t N(s)} F_0 \Big(x - \int_0^t a \Big) \\
    \|F\|_{\cnt_{\W^{m,p}}} \le e^{\const_{N}T}\|F_0\|_{\W^{m,p}(\dom)}.
\end{aligned}
\end{align}
\end{lem}
\begin{proof}
Firstly, based on the assumption \eqref{eq:coef_nuc_bd}, we have that $\|N\|_{\cnt} \le \|N\|_{\W^{1,p}} \le \const_{N},$ where $\const_{N} := \const_N + \const_D$. Based on this we have that 
\begin{align*}
\|\Psi^{N}_{s,t}\|_{\cnt} = \sup_{s,t \in [0,T]} e^{\int_s^t N(\tau) \dtau} \le e^{\const_{N}T} < \infty.
\end{align*}
Thus, for any $g \in L^{p}(\dom)$ we get
\begin{align*}
\|\Psi^{N}_{s,t} g\|_{L^p(\dom)} \le \|\Psi^{N}_{s,t}\|_{\cnt} \|g\|_{L^p(\dom)} \le e^{\const_{N}T} \|g\|_{L^p(\dom)}.
\end{align*}
Thus, we $\Psi^{N}_{s,t} \in \mathcal{L}(L^p;L^p)$ for all $ s,t \in [0,T]$. The corresponding semigroup property follows immediately from the additive property of the integral and semigroup property of the exponential function. 
Next, from \eqref{eq:coef_growth_rate_bd}, for $\Psi^{a}$ we first observe that, we have that \begin{align*}
    a &\le \const_a \text{ and } \int_s^t a(\tau) \le  \const_a T  \quad \forall s,t \in [0,T], \text{ thus } \left\|\int_s^t a(\tau) \right\|_{\cnt} \le \const_a T.
\end{align*}
Thus, for any $g \in L^p(\dom)$, consider an extension $\tilde g = \indc_{x \in \dom} g$ so that $\tilde g \in L^p(\R)$. Since $\Psi^{a}$ forms a two parameter evolution semigroup on $L^p(\R)$, by defining 
\begin{align*}
    \Psi^{a}_{s,t} g = \tilde g\left(\cdot - \int_s^t a\right)
\end{align*}
we get that $\Psi^{a}$ is a two parameter evolution semigroup on $L^p(\dom).$
Thus, combining the two, we can verify that 
\begin{align*}
\Psi^F_{0,t} F_0 &:= \Psi^{N}_{0,t} \Psi^a_{0,t} F_0(x) = e^{\int_0^t N(s)} \tilde F_0 \Big(x - \int_0^t a \Big), \quad \text{ with } \\
 \|\Psi^F_{t,s} F_0\|_{\cnt} &\le e^{\const_N(t-s)}\|F_0\|_{\cnt(\dom)} \le e^{\const_N(t-s)}\|F_0\|_{\W^{2,p}}
\end{align*}
satisfies the transport equation \eqref{eq:psd_iso} for any $F_0 \in \cnt^1(\bar \dom)$. 
Now taking $(F^n_0)_{n\in\N} \in \cnt^1(\bar \dom)$ such that $F^n_0 \to F_0$ in $\W^{1,p}(\dom)$, from the strong-continuity of $\Psi^F_{0,t} := \Psi^{N}_{0,t} \Psi^a_{0,t}$ we get that $$\Psi^F_{0,t} F^n_0 \to \Psi^F_{0,t} F_0 \text{ as } n \to \infty \text{ in } \W^{1,p}(\dom), \text{ uniformly } \forall t \in [0,T].$$ 
Thus, $t \mapsto F(t,\cdot) := \Psi^F_{0,t} F_0$ is an element of $\cnt([0,T];\W^{m,p}(\dom))$.
%

\end{proof}

\begin{cor} 
\label{cor:iso_psd_sol}
Let $C,H \in \cnt(T;\R)$, then the mapping $(C,H) \mapsto \Phi^F_{0,t}$ defined as 
\begin{align}
\label{eq:iso_psd_sg}
\begin{aligned}
\Phi^F_{0,t} &: \cnt(T;\R)^{\times 2} \to \mathcal{L}\Big(\W^{m,p}(\dom);\cnt\big(T;\W^{m,p}(\dom)\big)\Big), \: m \ge 1\\
\Phi^F_{0,t}(C,H) &:= \Psi^F_{0,t} := \Psi^{N(C,H)}_{0,t} \Psi^{a(C,H)}_{0,t}, ~\forall t \in [0,T]
\end{aligned}
\end{align}
serves as a solution operator for the transport equation \eqref{eq:psd_iso}. Thus, for any given $F_0 \in \W^{m,p}(\dom)$, 
\begin{align*}
    F^{C,H}(t,\cdot) = \Phi^F_{0,t} F_0 = \Psi^F_{0,t} F_0 = e^{\int_0^t N(C,H)} F_0 \Big(\cdot - \int_0^t a(C,H) \Big)
\end{align*}
solves the equation \eqref{eq:psd_iso}.
\end{cor}
\begin{proof}
Firstly, from Lemma \ref{lem:iso_psd_sol}, we observe that for more regular initial data, i.e. for $F_0 \in \W^{m,p}(\dom)$, we get that $F(t,\cdot) := \Psi^F_{0,t} F_0$ is an element of $C([0,T];W^{m,p}(\dom))$, $m \ge 1$. In view of this, letting $a(t) = a(C(t),H(t))$ and $N(t) = N(C(t),H(t))$ we directly obtain the above claim.
\end{proof}

\begin{cor}
\label{cor:iso_psd_sol_lp}
Let $(C,H) \in L^p\Big(\Om;\cnt(T;\R)^{\times 2}\Big)$ be a $(\F_t)_{t\ge0}$ adapted stochastic process, then the mapping $(C,H) \mapsto \Phi^F_{0,t}$ defined as 
\begin{align}
\label{eq:sol_psd_sg}
\begin{aligned}
\Phi^F_{0,t} &: L^p\left(\Om;\cnt(T;\R)^{\times 2}\right) \to L^{\infty}\Bigg(\Om;\mathcal{L}\Big(\W^{m,p}(\dom);\cnt\big([0,T];\W^{m,p}(\dom)\big)\Big) \Bigg), m \ge 1, \\
\Phi^F_{0,t}(C,H) &:= \Psi^F_{0,t} := \Psi^{N(C,H)}_{0,t} \Psi^{a(C,H)}_{0,t}, ~\forall t \in [0,T]
\end{aligned}
\end{align}
serves as a solution operator for the transport equation \eqref{eq:psd_iso} such that, for any given $F_0 \in \W^{m,p}(\dom)$, $\big(\Psi^F_{0,t}(C,H)\big)_{t\ge0}$ is a $(\F_t)_{t\ge0}$ adapted process that solves \eqref{eq:psd_iso} and has the explicit form
\begin{align*}
    F^{C,H}(t,\cdot) = \Phi^F_{0,t} F_0 =  \Psi^F_{0,t} F_0 &= e^{\int_0^t N(C,H)} F_0 \Big(\cdot - \int_0^t a(C,H) \Big) \quad \text{ almost surely } \\
    \|F\|_{L^{\infty}\big(\Om;\cnt\big([0,T]; \W^{m,p}(\dom)\big)\big)} &\le e^{\const_{N}T}\|F_0\|_{\W^{m,p}}.
\end{align*}
\end{cor}
\begin{proof}
    Since $a$ and $N$ are bounded functions, For $C, H \in L^p(\Om;\cnt(T;\R))$, $a, N$ are elements of $L^{\infty}(\Om;\cnt(T;\R))$ i.e. $a, N$ are almost surely bounded continuous functions. Thus, pointwise ($\omega$-wise) continuous paths exists, and thus we can apply the above Corollary to obtain the claim.
\end{proof}

\begin{lem}
\label{lem:FC_lip}
Let \( N(x,y) \) and \( a(x,y) \) be functions as per Assumption \ref{asm:coef_cond} and let \(F_0 \in \W^{2, p}(\dom) \) be Lipschitz continuous with respect to $x$, so that:
\begin{align*}
    \|F_0\|_{\W^{2, p}} &\leq \const_{F_0} \text{ \rm  and } |F_0(x) - F_0(y)| \leq \const_{F_0} |x - y| \text{ \rm  for some constant } \const_{F_0} > 0.
\end{align*}
Then, the function
\[
F_t(C,H) = e^{\int_0^t N(C(s),H(s))\,\ds} F_0 \left(x - \int_0^t a(C(s),H(s))\, \ds\right)
\]
is {Lipschitz} with respect to \( C, H \), that is, there exists a constant \( \const_{F} > 0 \) such that for any two pairs of functions \( C^1, H^1 \) and \( C^2, H^2 \)
\begin{align}
\label{eq:lip_psd}
\|F(C^1,H^1) - F(C^2,H^2)\|_{\cnt_T(\W^{2,p})} \leq \const_{F} \big(\|C^1 - C^2\|_{\cnt_T(\R)} + \|H^1 - H^2\|_{\cnt_T(\R)} \big).
\end{align}
\end{lem}
\begin{proof}
Let \( C^1, H^1 \) and \( C^2, H^2 \) be two pairs of functions. Define:
\[
I_N^1 = \int_0^t N(C^1(s),H^1(s))\,\ds, \quad I_N^2 = \int_0^t N(C^2(s),H^2(s))\,\ds,
\]
\[
I_a^1 = \int_0^t a(C^1(s),H^1(s))\,\ds, \quad I_a^2 = \int_0^t a(C^2(s),H^2(s))\,\ds.
\]
Then, the difference \( F(C^1,H^1) - F(C^2,H^2) \) can be written as:
\begin{align*}
F(C^1,H^1) - F(C^2,H^2) &= e^{I_N^1} F_0(x - I_a^1) - e^{I_N^2} F_0(x - I_a^2) \\
\implies \|F(C^1,H^1) - F(C^2,H^2)\|_{\cnt_T(\W^{2,p})} &\leq \|e^{I_N^1} F_0(x - I_a^1) - e^{I_N^2} F_0(x - I_a^1)\|_{\cnt_T(\W^{2,p})} \\
                                                &\hspace*{1cm} + \|e^{I_N^2} F_0(x - I_a^1) - e^{I_N^2} F_0(x - I_a^2)\|_{\cnt_T(\W^{2,p})}.
\end{align*}
Considering the first summand we have:
\begin{align*}
\|e^{I_N^1} F_0(x - I_a^1) &- e^{I_N^2} F_0(x - I_a^1)\|_{\cnt_T({\W^{2,p}(\dom)})} = \|F_0(x - I_a^1)\|_{\W^{2,p}(\dom)} \cdot \|e^{I_N^1} - e^{I_N^2}\|_{\cnt_T(\W^{2,p}(\dom))} \\
    &\leq \|F_0(x - I_a^1)\|_{\W^{2,p}(\dom)} ~ \left \| e^{\max(I_N^1, I_N^2)} \right\|_{\cnt_T(\W^{2,p}(\dom))} ~ \|I_N^1 - I_N^2\|_{\cnt_T(\W^{2,p}(\dom))} \\
    &\leq \const_{F_0} e^{\const_{N} T} \const_{N} \left( \|C^1 - C^2\|_{\cnt_T(\R)} + \|H^1 - H^2\|_{\cnt_T(\R)} \right).
\end{align*}
Since \( F_0\) and \( N(C(s),H(s)) \) are bounded and \( e^x \) is Lipschitz on bounded domains 
$$  \|F_0(x - I_a^1)\|_{\W^{2,p}} \leq \const_{F_0}, \quad \|I_N^1 - I_N^2\|_{\cnt_T(\R)} \leq \const_{N} \left( \|C^1 - C^2\|_{\cnt_T(\R)} + \|H^1 - H^2\|_{\cnt_T(\R)} \right).$$
Considering the second summand we have:
\begin{align*}   
\left\|e^{I_N^2} F_0(x - I_a^1) - e^{I_N^2} F_0(x - I_a^2)\right\|_{\cnt_T(\W^{2,p}(\dom))} &= \left\||e^{I_N^2}| \cdot |F_0(x - I_a^1) - F_0(x - I_a^2)| \right\|_{\cnt_T(\W^{2,p}(\dom))} \\
&\leq \|e^{I_N^2}\|_{\cnt_T(\R)} \|F_0(x - I_a^1) - F_0(x - I_a^2)\|_{\W^{2,p}(\dom)} \\
&\leq \|e^{I_N^2}\|_{\cnt_T(\R)} \const_{F_0} \|I_a^1 - I_a^2\|_{\cnt_T(\R)} \\
&\leq e^{\const_{N} T} \const_{F_0} \const_a \left( \|C^1 - C^2\|_{\cnt_T(\R)} + \|H^1 - H^2\|_{\cnt_T(\R)} \right)
\end{align*}
since \( F_0 \in \W^{2,p} \) and $\W^{2,p} \hookrightarrow C^{1,{1-{1\over p}}}$, $F_0$ is Lipschitz 
\begin{align*}
    \|F_0(x - I_a^1) - F_0(x - I_a^2)\|_{\cnt_T(\W^{2,p})} &\leq \const_{F_0} |\dom| \|I_a^1 - I_a^2\|_{\cnt_T(\R)} \\
    \|I_a^1 - I_a^2\|_{\cnt_T(\R)} &\leq \const_a |\dom| \left( \|C^1 - C^2\|_{\cnt_T(\R)} + \|H^1 - H^2\|_{\cnt_T(\R)} \right).
\end{align*}
Combining the two we get
\begin{align*}
\|F(C^1,H^1) - F(C^2,H^2)\|_{\W^{1,p}} &\leq \const_{F_0} e^{\const_{N} T} \const_{N}|\dom| \left( \|C^1 - C^2\|_{\cnt_T(\R)} + \|H^1 - H^2\|_{\cnt_T(\R)} \right) \\
    &\hspace*{2cm} + e^{\const_{N} T} \const_{F_0} \const_a |\dom| \left( \|C^1 - C^2\|_{\cnt_T(\R)} + \|H^1 - H^2\|_{\cnt_T(\R)} \right) \\
    &= \left( \const_{F_0} e^{\const_{N} T} \const_{N} |\dom| + e^{\const_{N} T} \const_{F_0}\const_a |\dom| \right) \left( \|C^1 - C^2\|_{\cnt_T(\R)} + \|H^1 - H^2\|_{\cnt_T(\R)} \right) \\
    &= \const_{F} \left(\|C^1 - C^2\|_{\cnt_T(\R)} + \|H^1 - H^2\|_{\cnt_T(\R)}\right), \quad \const_{F} := \const_{F_0} e^{\const_{N} T}|\dom|(\const_{N} + \const_a).
\end{align*}
\end{proof}

\begin{lem}
\label{lem:mom_lip}
For $C,H \in L^p\left(\Om; \cnt([0,T];\W^{1,p})\right)$ and $F_0 \in \W^{1,p}(\dom)$ given, let $F \in L^{\infty}\left(\Om; \cnt([0,T];\W^{1,p}(\dom;\R))\right)$ be the solution obtained as per Corollary \ref{cor:iso_psd_sol_lp}. Then the mapping $F \mapsto S^F$, defined as per \eqref{eq:F_mom}, satisfies the following properties:
\begin{align}
\label{eq:FC_mom}
\begin{aligned}
S^{C,H} := S[F^{C,H}] &:= \int_{\dom} x^2 F^{C,H}(\cdot,x) \dx \\
\|S^{C,H}\|_{L^{\infty}_{\F}(\cnt(\R))}  &\le \const_{F_0} |\dom| e^{\const_{N}}T (\const_a |\dom| + T)^2 \\
\|S^1 - S^2\|^{2p}_{L^p_{\F}(\cnt(\R))} &\le \const^p_{S} \left(\|C^1 - C^2\|^p_{\cnt_T(\R)} + \|H^1 - H^2\|^p_{\cnt_T(\R)}\right).
\end{aligned}
\end{align}
\end{lem}
\begin{proof} Since $0 \le F_0 \in \W^{1,p}(\dom)$ and $L^p \hookrightarrow L^1(\dom)$ for $p > 1$ we have that $ \Upsilon^0_0 := \int_{\dom} F_0(x) \dx < \infty$. 
Similarly, due to compactness of $\dom \subset \R$, we also have that
    \begin{align*}
        \Upsilon^n_0 &= \int_{\dom} x^n F_0(x) \dx 
    \implies \Upsilon^n_0 \le |\dom|^{n} \int_{\dom} F_0(x) ~ \dx = \const^{n}_{\Upsilon} \const_{\Upsilon^0_0}
    \end{align*}
    where $\const_{\Upsilon} := |\dom|$ and $\const_{\Upsilon^0_0} := \Upsilon_0$. Thus, applying Lemma \ref{lem:F_mom_bd} we get that
    \begin{align*}
        S^{C,H} \le \const_{\Upsilon^0_0}  e^{\const_{N}T} (\const_a \const_{\Upsilon} + T)^2.
    \end{align*}
    Since the above inequality holds almost surely, we consequently also get that
    \begin{align*}
        \|S^{C,H}\|_{L^{\infty}(\Om)} \le \const_{\Upsilon^0_0} e^{\const_{N}T} (\const_a \const_{\Upsilon} + T)^2.
    \end{align*}
Next, let $C^1, C^2, H^1, H^2 \in L^p(\Om;C([0,T];\R))$ be arbitrary functions, then denoting $S^1 := S^{C^1, H^1}$, $S^2 := S^{C^2, H^2}$, $F^1 := F^{C^1, H^1}$, $F^2 := F^{C^2, H^2}$ we get
\begin{align*}
    S^1 - S^2 &= \int_{\dom} x^2 F^1 - \int_{\dom} x^2 F^2 = \int_{\dom} x^2 (F^1 - F^2) \\
    \implies \|S^1 - S^2\|_{\cnt(\R)} &\le \|F^1 - F^2\|_{\cnt(\R)} \int_{\dom} x^2 \\
    &\le |\dom|^3 \|F^1 - F^2\|_{\cnt(\R)}  \\
    &\le \const_{\W} \const_F |\dom|^3 \left(\|C^1 - C^2\|_{\cnt(\R)} + \|H^1 - H^2\|_{\cnt(\R)}\right), \quad (\text{using Lemma \ref{lem:FC_lip}})\\
\implies \E\left[\|S^1 - S^2\|^{2p}_{\cnt(\R)}\right] &= \E\left[\|S^1 - S^2\|^{p}_{\cnt(\R)} \|S^1 - S^2\|^{p}_{\cnt(\R)}\right] \\
&\le 2^p \E\left[\|S^{C,H}\|^{p}_{\cnt(\R)}\right] \E\left[\|S^1 - S^2\|^{p}_{\cnt(\R)}\right] \\
&\le 2^p \E\left[ \left(\const_{\Upsilon^0_0} e^{\const_{N}T} (\const_a \const_{\Upsilon} + T)^2\right)^{p} \right] \E\left[\|S^1 - S^2\|^{p}_{\cnt(\R)}\right] \\
&\le \left(2\const_{\Upsilon^0_0} e^{\const_{N}T} (\const_a \const_{\Upsilon} + T)^2\right)^{p} \E\left[\|S^1 - S^2\|^{p}_{\cnt(\R)}\right] \\
&\le \const^p_{S} \left( \|C^1 - C^2\|^p_{L^p_{\F}(\cnt(\R))} + \|H^1 - H^2\|^p_{L^p_{\F}(\cnt(\R))} \right)
\end{align*}
where, $\const^p_{S} := \left(2\const_{\Upsilon^0_0} e^{\const_{N}T} (\const_a \const_{\Upsilon} + T)^2\right) \left(\const_{\W} \const_F |\dom|^3 \right)^p.$ 
$$\implies \|S^1 - S^2\|^{2p}_{L^{2p}_{\F}(\cnt(\R))} \le \const^p_{S} \left( \|C^1 - C^2\|^p_{L^p_{\F}(\cnt(\R))} + \|H^1 - H^2\|^p_{L^p_{\F}(\cnt(\R))} \right).$$
\end{proof}

\subsection{Analysis of the noisy kinetic system}
Consider the reaction kinetic system in isolation, i.e., given $F$, we consider the following SDE system
\begin{subequations}
\label{eq:noisy_kin_sys}
\begin{align}
dC(t) &= \big[ r(t) - C(t) \tilde{k}_v(t)  - b^F(t)\big]dt + \sigma_C C(t) dW^1_t, \quad t > 0 \label{eq:C_sde}\\ 
dQ(t) &= -[K_{\rm sp} P(H) U_r + \tilde{k}_v(t)] Q(t) \dt + \sigma_Q Q(t) dW^2_t, \quad  t > 0 \label{eq:Q_sde}\\
dH(t) &= k_H U_H ~ \dt + \sigma_H H(t) dW^3_t, \quad  t > 0 \label{eq:pH_sde}\\
\dot{R} &= k_v(t), \quad C(0) = C_0, \quad R(0) = R_0, \quad H(0) = H_0, \quad Q(0)=Q_0  \\
b^F(t) &:= {\tilde \rho ~ k_v(t) S^F_t a(t)} + {\tilde \rho~v_{\rm nuc} N(t)}.
\end{align}
\end{subequations}
Firstly, we establish some properties of $b^F$.
\begin{lem} 
\label{lem:decay_term_bound}
Let $b^F$ be as defined in \eqref{eq:noisy_kin_sys} with $a(t), N(t), R(t)$ and $S^F(t)$ defined as per \eqref{eq:psd_trans} and \eqref{eq:F_mom}. If $C,H \in L^{p}(\Om; \cnt([0,T]))$, then $b^F$ satisfies the following properties:
\begin{align}
\begin{aligned}
b^F \ge 0  \text{ \rm with } b^F(\R^-) = 0 ~~ a.s. \quad &\text{ \rm and } \quad b^F_t \to 0  \text{ almost surely, as } C_t \to 0 \\
b^F \in L^{2p}(\Om;\cnt(T;\R)) \quad &\text{ \rm with } \quad \|b^F\|_{L^{2p}_{\F}(\cnt(\R))} \le \const_b \\
\|b^{h}(C^1,H^1) - b^{h}(C^2,H^2)\|^{2p}_{L^{2p}_{\F}(\cnt(\R))} &\le \const_b^p\left( \|C^1 - C^2\|^p_{L^p_{\F}(\cnt(\R))} + \|H^1 - H^2\|^p_{L^p_{\F}(\cnt(\R))} \right).
\end{aligned}
\end{align}
\end{lem}

\begin{proof}
Since $a(C,H) = \const_g \tanh(\tilde C^2)$ and $N(C,H) = \exp(-{\delta \over \tilde{C}^2})$ based on \eqref{eq:growth_nuc_lip} in 
Lemma $\ref{lem:growth_lip}$ we obtain that
\begin{align*}
    |a(C^1,H^1) - a(C^2,H^2)| &\le \const_{L_a} (|C^1 - C^2| + |H^1-H^2| )\\
    \implies \sup_{t\in[0,T]} |a(C^1,H^1) - a(C^2,H^2)| &\le (\const_{L_a})\sup_{t\in[0,T]} (|C^1 - C^2| + |H^1-H^2| )
\end{align*}
where $\const_{L_a} := (L_a)^{1 \over p}$.
\begin{align*}
    \implies \E\left[\left(\sup_{t\in[0,T]} |a(C^1,H^1) - a(C^2,H^2)|\right)^{2p}\right] &= \E\Bigg[\left(\sup_{t\in[0,T]} |a(C^1,H^1) - a(C^2,H^2)|\right)^{p}   \\
    &\hspace{3cm} \times \left(\sup_{t\in[0,T]} |a(C^1,H^1) - a(C^2,H^2)|\right)^{p}\Bigg]\\
    &\le \const_{L_a}^{2p} \E\left[ 2^p \left(\sup_{t\in[0,T]} |C^1 - C^2|\right)^{p} + 2^p \left(\sup_{t\in[0,T]} |H^1 - H^2|\right)^{p}\right] \\
    \implies \|a(C^1,H^1) - a(C^2,H^2)\|^{2p}_{L^{2p}_{\F}(\cnt(\R))} &\le 2^p \const_{L_a}^{2p} \Big( \|C^1 - C^2\|^p_{L^{p}_{\F}(\cnt(\R))} + \|H^1 - H^2\|^p_{L^{p}_{\F}(\cnt(\R))} \Big)
\end{align*}
Similarly, we also obtain that 
\begin{align*}
\|N(C^1,H^1) - N(C^2,H^2)\|^{2p}_{L^{2p}_{\F}(\cnt(\R))} &\le 2 \const^{2p}_{L_N} \Big( \|C^1 - C^2\|^{p}_{L^{p}_{\F}(\cnt(\R))} + \|H^1 - H^2\|^{p}_{L^{p}_{\F}(\cnt(\R))} \Big)
\end{align*}
Next, considering the expression for $b^F$, using the definitions of $\tilde \rho $ and $ S^F_t $, we have
\begin{align}
\label{eq:f_definition}
\begin{aligned}
    b^F(C,R) &= { \tilde \rho k_v(t) a(C,H) S^F_t} + {\tilde \rho v_{\rm nuc} N(t)} \\ 
    \implies \E[\|b^F\|^{2p}_{\cnt}] &\le \const_{\tilde\rho} \|S^F\|_{\cnt} \E[\|a\|^{2p}_{\cnt}]   + \const_{\tilde\rho} \|v_{\rm nuc}\|_{\cnt} \E[\|N\|^{2p}_{\cnt}] \\
    &\le \const_{\tilde \rho} (\const_{a} \const_{S^F} + \const_{v_{\rm nuc}} \const_{N}) \\
    &= \tilde \const_{b}
\end{aligned}
\end{align}
Now, letting $b^1 := b^{h}(C^1,H^1)$ and $b^2 := b^{h}(C^2,H^2)$, consider the difference $|b^1 - b^2|$, 
\begin{equation}
\|b^1 - b^2\|^{2p}_{L^{2p}_{\F}(\cnt(\R))} = \underbrace{\left\| {\tilde \rho k_v(t) a^1 S^1} - {\tilde \rho k_v(t) a^2 S^2} \right\|^{2p}_{L^{2p}_{\F}(\cnt(\R))}}_{\text{Term 1}} + \underbrace{\left\| {\tilde  \rho v_{\rm nuc} N^1} - {\tilde \rho v_{\rm nuc} N^2)} \right\|^{2p}_{L^{2p}_{\F}(\cnt(\R))}}_{\text{Term 2}}
\label{eq:main_difference}
\end{equation}
Term1 can be estimated as
\begin{align}
\label{eq:coefficient_diff}
\text{Term 1} &= (\tilde \rho k_v)^{2p}\left\| \left( {a^1} - {a^{2p}} \right) S^1 \right\|^{2p}_{L^{2p}_{\F}(\cnt(\R))} + (\tilde \rho k_v)^{2p}\left\| \left({S^1} - {S^{2p}} \right) a^{2p} \right\|^{2p}_{L^{2p}_{\F}(\cnt(\R))}\notag\\
&\le (\tilde \rho k_v)^{2p}  \|S\|^{2p}_{\cnt} \left\|{a^1 - a^{p}} \right\|^{2p}_{L^{2p}_{\F}(\cnt(\R))} + (\tilde \rho k_v)^{2p}  \|a^{2p}\|^{2p}_{\cnt} \left\|{S^1 - S^{p}} \right\|^{2p}_{L^{2p}_{\F}(\cnt(\R))}\notag\\
&\le (\tilde \rho k_v)^{2p}  \|S\|^{2p}_{\cnt} \const^{2p}_{L_a} 2^p \Big( \left\|C^1 - C^2\right\|^p_{L^{p}_{\F}(\cnt(\R))} + \left\|H^1 - H^{p}\right\|^p_{L^{p}_{\F}(\cnt(\R))} \Big) \\
&\hspace*{2cm} + (\tilde \rho k_v)^{2p}  \|a^{2p}\|^{2p}_{\cnt}  \const^p_S \Big( \left\|C^1 - C^2\right\|^p_{L^{p}_{\F}(\cnt(\R))} + \left\|H^1 - H^{p}\right\|^p_{L^{p}_{\F}(\cnt(\R))} \Big) \notag\\
&\le (\const_{\tilde \rho} \const_v \const_{L_a})^{2p} \const^p_{S} (2^p \const^p_{S} + 1) \Big( \left\|C^1 - C^2\right\|^p_{L^{p}_{\F}(\cnt(\R))} + \left\|H^1 - H^{p}\right\|^p_{L^{p}_{\F}(\cnt(\R))} \Big)
\end{align}
Similarly, Term2 can be estimated as:
\begin{align}
\text{Term 2} &= \left\| {\tilde \rho v_{\rm nuc} N^1} - {\tilde \rho v_{\rm nuc} N^2)} \right\|^{2p}_{L^{2p}} \\
&\le \left| {\tilde \rho v_{\rm nuc}}\right|^{2p}   \left\|N^1 - N^2 \right\|^p_{L^{2p}_{\F}(\cnt(\R))} \notag\\
&\le \left| {\tilde \rho v_{\rm nuc}}  \right|^{2p} \const^{2p}_{L_N} 2^p \Big( \left\|C^1 - C^2\right\|^p_{L^{p}_{\F}(\cnt(\R))} + \left\|H^1 - H^2\right\|^p_{L^{p}_{\F}(\cnt(\R))} \Big) \\
&\le \const^{2p}_{\tilde \rho} \const^{2p}_{v_{\rm nuc}} \const^{2p}_{L_{N}} 2^p \Big( \left\|C^1 - C^2\right\|^p_{L^{p}_{\F}(\cnt(\R))} + \left\|H^1 - H^2\right\|^p_{L^{p}_{\F}(\cnt(\R))} \Big)
\label{eq:integral_bound}
\end{align}
Combining the two, i.e. substituting \eqref{eq:coefficient_diff} and \eqref{eq:integral_bound} into \eqref{eq:main_difference}, we get
\[ \|b^F(C^1,H^1) - b^F(C^2,H^2)\|^{2p}_{L^{2p}_{\F}(\cnt(\R))} \leq \const^{p}_b \Big( \left\|C^1 - C^2\right\|^p_{L^{p}_{\F}(\cnt(\R))} + \left\|H^1 - H^2\right\|^p_{L^{p}_{\F}(\cnt(\R))} \Big) \]
where $\const^p_b := \const^{2p}_{\tilde \rho}   \Big( \const^{2p}_v \const^{2p}_{L_a} \const^p_{S} (2^p \const^p_{S} + 1) + \const^{2p}_{v_{\rm nuc}} \const^{2p}_{L_{N}} 2^p\Big).$
\end{proof}

\subsubsection{Estimates of geometric BM}
\begin{lem}
Let $r(t)$, $P(H)$ and $U_t$ be some given functions of time, then the time varying processes
\begin{align}
\label{eq:sg_kin_eqns}
\begin{aligned}
\Psi^C_{s,t} &:= \Big(e^{- {1\over 2} \int_s^t 2\tilde{k}_v(\tau) + \sigma^2_C (\tau) \dtau }\Big) \Big(e^{\int_s^t \sigma_C(\tau)~{\rm dW}^1_{\tau}}\Big) \\ 
\Psi^Q_{s,t} &:= \Big(e^{- {1\over 2} \int_s^t 2 \tilde P(\tau) + \sigma^2_Q(\tau) \dtau }\Big) \Big(e^{\int_s^t \sigma_Q(\tau) ~ {\rm dW}^2_{\tau}}\Big), \quad  \tilde P := \left(\Ksp P U_r + \tilde{k}_v(t)\right) \\
\Psi^H_{s,t} &:= \Big(e^{- {1\over 2}\int_s^t  \sigma^2_H(\tau) \dtau }\Big) \Big(e^{\int_s^t \sigma_H(\tau) ~ {\rm dW}^3_{\tau}}\Big)
\end{aligned}
\end{align}
satisfy the following, corresponding, SDEs
\begin{align}
\label{eq:sg_sde}
\begin{aligned}
    d\Psi^C_{0,t} &= - \tilde{k}_v(t) \Psi^C_{0,t} \dt + \Psi^C_{0,t} \sigma_C(t) \dWt{1} \\
    d\Psi^Q_{0,t} &=  - \tilde{P}(t) \Psi^Q_{0,t} \dt + \Psi^Q_{0,t} \sigma_Q(t) \dWt{2} \\
    d\Psi^H_{0,t} &= \Psi^C_{0,t} \sigma_H(t) \dWt{3}.
\end{aligned}
\end{align}
Furthermore, their first two moments have the following expression
\begin{align}
\label{eq:sg_mom_est}
\begin{aligned}
\E[\psic_{s,t}] &=\Big(e^{-\int_s^t \tilde k_v(\tau) \dtau }\Big), \hspace*{-0.04cm} & \E[(\psic_{s,t})^{2p}] &= \Big(e^{\int_s^t p(2p-1)\sigma^2_C(\tau) - 2p \tilde k (\tau) \dtau }\Big), \\
\E[\psiq_{s,t}] &=\Big(e^{-\int_s^t \tilde P(\tau) \dtau }\Big), \hspace*{-0.04cm} & \E[(\psiq_{s,t})^{2p}] &= \Big(e^{\int_s^t p(2p-1)\sigma^2_Q(\tau) - 2p \tilde P(\tau) \dtau }\Big), \\
\E[\psih_{s,t}] &= 1, \hspace*{-0.04cm} & \E[(\psih_{s,t})^{2p}] &= \Big(e^{\int_s^t p(2p-1)\sigma^2_H(\tau) \dtau }\Big). 
\end{aligned}
\end{align}
\end{lem}

\begin{proof} Note that the expressions in \eqref{eq:sg_kin_eqns} have the form of an exponential Brownian motion. Thus, based on the formulas \eqref{eq:Exp_OU_proc} and \eqref{eq:Exp_OU_proc_solved} the SDE relations \eqref{eq:sg_sde} follows immediately. 
Consequently, the moment estimates \eqref{eq:sg_mom_est} follows immediately from the independence of $W_t$ and $r(t),  P(H)$, Ito calculus and the fact that $\E[e^{W_t}] = e^{{1 \over 2} t}$.
\end{proof}
\begin{cor} Let $r(t)$, $P(H)$ and $U_H(t)$ be bounded functions of time. Similarly, let the covariance functions $\sigma_C$, $\sigma_H$, $\sigma_Q$ be also bounded in time such that
\begin{align}
\label{eq:sg_coefs_bound}
\begin{array}{llll}
\hspace*{-0.75cm}\|\tilde k_v\|_{L^\infty([0,T])} \le \const_{k_v},& \|U_r\|_{L^\infty([0,T])} \le \const_{U_r},& \|P\|_{L^\infty([0,T])} \le \const_{P},& \|k_H U_H\|_{L^\infty([0,T])} \le \const_{U_H}, \\
\hspace*{-0.75cm}\|\tilde P\|_{L^\infty([0,T])} \le \const_{\tilde P},& \|\sigma^2_C\|_{L^\infty([0,T])} \le \const_{\sigma_C},& \|\sigma^2_Q\|_{L^\infty([0,T])} \le \const_{\sigma_Q},& \|\sigma^2_H\|_{L^\infty([0,T])} \le \const_{\sigma_H},
\end{array}
\end{align}
with $\const_{\tilde P} := \Ksp \const_{P} \const_{U_r}$. Then the moments of $\psic$, $\psih$ and $\psiq$ satisfy the following, corresponding, inequalities
\begin{align}
\label{eq:sg_mom_est_upp}
    \begin{aligned}
        \E[(\psih_{s,t})^{2p}] &\le e^{(t-s)\const_{p,\psi^h}}, \quad \E \left [\int_0^T (\psih_{s,t})^{2p} \ds \right] \le T~\const_{p,\psi^h} \exp(T~\const_{p,\psi^h}) \\
        \E[(\psiq_{s,t})^{2p}] &\le e^{(t-s)\const_{p,\psi^q}}, \quad \E \left [\int_0^T (\psiq_{s,t})^{2p} \ds \right] \le T~\const_{p,\psi^q} \exp(T~\const_{p,\psi^q}) \\
        \E[(\psic_{s,t})^{2p}] &\le e^{(t-s)\const_{p,\psi^c}}, \quad \E \left [\int_0^T (\psic_{s,t})^{2p} \ds \right] \le T~\const_{p,\psi^c} \exp(T~\const_{p,\psi^c})
    \end{aligned}
\end{align}
where,
\begin{align*}
    \const_{p,\psi^h} &:= 2p(2p-1)\const_{\sigma_H}, \quad \const_{p,\psi^q} := 2p(2p-1)\const_{\sigma_Q} + 2p \const_{\tilde P}, \quad 
    \const_{p,\psi^c} := 2p(2p-1)\const_{\sigma_C}. 
\end{align*}
Furthermore, we also have the following uniform upper bounds
\begin{align}
\label{eq:sg_sup_lp_est}
    \begin{aligned}
        \E\left[\left(\sup_{s \le t, s,t \in [0,T]}\psih_{s,t}\right)^{2p}\right] &\le \left({2p \over 2p -1}\right)^{2p} e^{T\const_{p,\psi^h}}, \\
        \E\left[\left(\sup_{s \le t, s,t \in [0,T]}\psiq_{s,t}\right)^{2p}\right] &\le \left({2p \over 2p -1}\right)^{2p} e^{T\const_{p,\psi^q}}, \\
        \E\left[\left(\sup_{s \le t, s,t \in [0,T]}\psic_{s,t}\right)^{2p}\right] &\le \left({2p \over 2p -1}\right)^{2p} e^{T\const_{p,\psi^c}}.
    \end{aligned}
\end{align}
\end{cor}
\begin{proof}
    The claim \eqref{eq:sg_coefs_bound} follows directly by applying the \eqref{eq:sg_coefs_bound} to \eqref{eq:sg_mom_est}. For \eqref{eq:sg_sup_lp_est} we proceed as follows.
  From \eqref{eq:sg_kin_eqns}, since $\Psi^H_{s,t}$ is a $\F_t$ adapted martingale. Thus, we can apply Doob's $L^{^p}$ inequality to obtain
  \begin{align*}
      \sup_{s} \Psi^H_{s,t} &\le \sup_{s \in [0,T]} \Psi^H_{s,T} = \sup_{t \in [0,T]} (\Psi^H_{0,t}) \\
    \implies \E\left[\left(\sup_{t \in [0,T]}\psih_{s,t}\right)^{2p}\right] &\le \left({2p \over 2p -1}\right)^{2p} \E\left[\left(\psih_{0,T}\right)^{2p}\right] \le \left({2p \over 2p -1}\right)^{2p} e^{T\const_{p,\psi^h}}.
  \end{align*}
  Similarly, for $\Psi^Q_{s,t}$ we observe that
  \begin{align*}
  \sup_{t \in [0,T]} \Psi^Q_{s,t} \le \sup_{t \in [0,T]}\Big(e^{\int_0^t \tilde P(\tau) \dtau }\Big) \sup_{t \in [0,T]} \left[\Big(e^{-{1\over 2} \int_0^t \sigma^2_Q(\tau) \dtau }\Big) \Big(e^{\int_0^t \sigma_Q(\tau) dW_{\tau}}\Big)\right]
  \end{align*}
  Since the 2nd term on the RHS is a $\F_t$ adapted martingale, thus applying Doob's $L^{^p}$ inequality we obtain
  \begin{align*}
      \E\left[\left(\sup_{t \in [0,T]}\psiq_{s,t}\right)^{2p}\right] &\le \sup_{t \in [0,T]}\Big(e^{2p\int_0^t \tilde P(\tau) \dtau }\Big) \left({2p \over 2p -1}\right)^{2p} \E\left[\left(\psiq_{0,T}\right)^{2p}\right] \\
      &\le e^{2p \const_{\tilde P}} \left({2p \over 2p -1}\right)^{2p} e^{(p(2p-1)\const_{\sigma_Q})T}\\
      &= \left({2p \over 2p -1}\right)^{2p} e^{T\const_{p,\psi^q}}.
  \end{align*}
  Following, exactly like above for $\psic_{s,t}$ and observing that 
  $$\Big(e^{-{1\over 2} \int_0^t \sigma^2_C(\tau) \dtau + \int_0^t \sigma_C(\tau) dW_{\tau}}\Big)_{t\ge0}$$ is $\F_t$ adapted Martingale, we arrive at the conclusion that
  \begin{align*}
      \E\left[\left(\sup_{t \in [0,T]}\psiq_{s,t}\right)^{2p}\right] &\le \left({2p \over 2p -1}\right)^{2p} e^{T\const_{p,\psi^c}}.
  \end{align*}
\end{proof}

\subsubsection{Analysis of pH dynamics}
\begin{lem}
\label{lem:ph_dyn_wlpd}
    Let $U_H \in L^{\infty}(\Om;L^{\infty}([0,T];\R))$ be a given bounded $\F_t$ measurable input process, let $\sigma_H \in L^{\infty}([0,T];\R)$ and $H_0 \in L^{p}(\Om;\R)$, then the $\rm pH$ dynamics specified by 
    \begin{align}
    \label{eq:pH_sde_rep}
    \begin{aligned}
        dH(t) &= k_H U_H ~ \dt + \sigma_H H(t) dW^3_t, \quad  t > 0 \\
        H(0) &= H_0
    \end{aligned}
    \end{align}
    has a unique strong solution $H \in L^{p}(\Om;L^{\infty}([0,T];\R))$, with $H = (H_t)_{t\ge0}$ and 
    \begin{align}
    \label{eq:pH_sde_ineq}
    \begin{aligned}
        H_t &= \psih_{0,t} H_0 + \int_0^t \psih_{s,t} k_H U_H(s) \ds \\
        \E[|H_t|^{p}] &\le  e^{\const_{p,\psi^h} t} \left( \const^p_{H_{_0}} + \const_{p,\psi^h} \const^p_{U_H} ~ t \right)\\
        \E\left[\left(\sup_{t\in[0,T]}|H_t|\right)^{p}\right] &\le \left({p \over p -1}\right)^{p} e^{T\const_{p,\psi^h}} \Big(\const^p_{H_{_0}} +  \const^p_{U_H} T \Big) \\
        \E(|H_t - H_{s}|^{2p}) &\le 2^{2p-1} (t-s)^{p} \Bigg( k_H \const_u + \const_{p,H_0,\sigma_H, U_H} e^{p(p-1) \const_{\sigma_H}(t-s)} \Bigg) \\
        \E \left[\int_0^T |H_t|^{p} \dt \right]  &\le T~\const_{p,\psi^h} \exp(T~\const_{p,\psi^h}) \left( \const_{H_{_0}} + \const_{U_H^{2p}} ~ T \right)
    \end{aligned}
    \end{align}
    where,
    \begin{align}
    \begin{array}{lll}
       \E[\|\const_H U_H\|^{p}_{L^\infty_T(\R)}] \le \const^{p}_{U_H}, & \E[|H_0|^{p}] \le \const^p_{H_0}, &  \|\sigma_H\|^{p}_{L^\infty_T(\R)} \le \const^{p}_{\sigma_H}.
    \end{array}
    \end{align}
    Furthermore, $H$ has a $\vartheta$-H\"older continuous modification for $\vartheta < 1/4$. 
\end{lem}
\begin{proof}
    For a given bounded input process $U_H \in L^{2p}(\Om;L^{\infty}([0,T];\R))$, using \eqref{eq:Exp_OU_proc} and \eqref{eq:Exp_OU_proc_solved} we get that solution to \eqref{eq:pH_sde} is uniquely given by 
    $$\Phi^{H}_t H_0 := H_t = \psih_{0,t} H_0 + \int_0^t \psih_{s,t} k_H U_H(s) \ds.$$ 
    Thus,
    \begin{align}
    \label{eq:htp_ineq}
        |H_t|^{p} &\le | \psih_{0,t} H_0 |^{p} + \int_0^t \left |\psih_{s,t}k_H  U_H(s) \right |^{p} \ds
    \end{align}
    \begin{align}
    \implies \E \left[ |H_t|^{p} \right] &\le \E[ (\psih_{0,t})^{p}] \E[(H_0)^{p}] +  \int_0^t\E[\left (\psih_{s,t} \right)^{p}] \E[\left | k_H U_H(s) \right |^{p}] \ds \\
     &\le \E[(\psih_{0,t})^{p}] \E[(H_0)^{p}] + \E \left [\left \|k_H  U_H \right \|^{p}_{\infty} \right] \int_0^t \E[\left (\psih_{s,t} \right)^{p}]  \ds \\
     &\le e^{p(p-1)\const_{\sigma_H}t} ~ \const^p_{H_{_0}} + \const^p_{U_H}  [e^{p(p-1)\const_{\sigma_H}t} - 1] \\
     &\le e^{p(p-1)\const_{\sigma_H}t} ~ \left( \const^p_{H_{_0}} + p(p-1)\const_{\sigma_H}\const^p_{U_H} ~ t \right)\\
     &= e^{\const_{p,\psi^h} t} \left( \const^p_{H_{_0}} + \const_{p,\psi^h} \const^p_{U_H} ~ t \right).
    \end{align}
    Taking supremum of \eqref{eq:htp_ineq} we obtain
    \begin{align*}
        \left(\sup_{t\in[0,T]}|H_t|\right)^{p} &\le \left(\sup_{t\in[0,T]} | \psih_{0,t} H_0 |\right)^{p} + \left(\sup_{t\in[0,T]} \int_0^t \left |\psih_{s,t} k_H  U_H(s) \right | \ds\right)^p \\
        &\le \left(\sup_{t\in[0,T]} | \psih_{0,t} H_0 |\right)^{p} +   \left(\int_0^t \sup_{t\in[0,T]}\left |\psih_{s,t} k_H U_H(s) \right | \ds \right)^p \\
        &\le \left(\sup_{t\in[0,T]} | \psih_{0,t} H_0 |\right)^{p} +  \left(\int_0^t \sup_{t\in[0,T]}\left |\psih_{0,t}\right | \sup_{s\in[0,T]}\left| k_H U_H(s) \right | \ds \right)^p\\
        &\le \left(\sup_{t\in[0,T]} | \psih_{0,t} H_0 |\right)^{p} + \int_0^t  \left(\sup_{t\in[0,T]}\left |\psih_{0,t}\right | \right)^p \left(\sup_{s\in[0,T]}\left|U_H(s) \right |\right)^p \ds \\
        &\le \left(\sup_{t\in[0,T]} | \psih_{0,t}|\right)^{p} |H_0 |^{p} +   \int_0^t \left(\sup_{t\in[0,T]} | \psih_{0,t}|\right)^{p} \const_{U_H}^p \ds
    \end{align*}
    \begin{align*}
        \implies \E\left[\left(\sup_{t\in[0,T]}|H_t|\right)^{p}\right] &\le \E \left[\left(\sup_{t\in[0,T]} | \psih_{0,t}|\right)^{p}\right] \E\left[|H_0 |^{p}\right] +  \const_{U_H}^p \int_0^t \E\left[\left(\sup_{t\in[0,T]} | \psih_{0,t}|\right)^{p}\right]  \ds \\
        &=\E \left[\left(\sup_{t\in[0,T]} | \psih_{0,t}|\right)^{p}\right] \Big(\E\left[|H_0 |^{p}\right] +  \const_{U_H}^p  T \Big) \\
        &\le \E \left[\left(\sup_{t\in[0,T]} | \psih_{0,t}|\right)^{p}\right] \Big(\const^p_{H_{_0}} + \const_{U_H}^p T \Big) \\
        &\le \left({p \over p -1}\right)^{p} e^{T \const_{p,\psi^h}} \Big(\const^p_{H_{_0}} +  \const_{U_H}^p T \Big).
    \end{align*}
    
    Integrating \eqref{eq:htp_ineq} we get
    \begin{align}
    \int_0^T |H_t|^{p} \dt &\le \int_0^T | \psih_{0,t} H_0 |^{p} \dt + \int_0^T \int_0^t \left |\psih_{s,t} k_H U_H(s) \right |^{p} \ds \dt \\
    \implies \E \left[\int_0^T |H_t|^{p} \dt \right]    &\le \E \left[\int_0^T | \psih_{0,t} H_0 |^{p} \dt \right] +  \E \left[\int_0^T \int_0^t \left |\psih_{s,t} k_H U_H(s) \right |^{p} \ds \dt \right] \\
    &= \int_0^T \E [ | \psih_{0,t} H_0 |^{p}] \dt  +  \int_0^T \int_0^t \E[ \left |\psih_{s,t} k_H U_H(s) \right |^{p}] \ds \dt \\
    &\le \int_0^T \E[ (\psih_{0,t})^{p}] \E[(H_0)^{p}] +  \int_0^T \int_0^t\E[\left (\psih_{s,t} \right)^{p}] \left\|\left ( k_H U_H(s) \right )^{p}\right\|_{\infty} \ds \\
     &\le \int_0^T \E[(\psih_{0,t})^{p}] ~ \E[(H_0)^{p}] + \E \left [\left \| k_H U_H \right \|^{p}_{\infty} \right] \int_0^T \int_0^T \E \left[\left (\psih_{s,t} \right)^{2p} \right]^{1\over2}  \ds \\
     &\le T~\const_{p,\psi^h} \exp(T~\const_{p,\psi^h}) \left( \const^p_{H_{_0}} + \const^p_{U_H} ~ T \right).
    \end{align}
    Now using SDE relation \eqref{eq:pH_sde}, $(a+b)^{2p} \le 2^{2p-1}(a^{2p} + b^{2p})$  and H\"older inequality
    $$ \E\left[\left|\int_s^t f(s) dW_s\right|^{2p}\right] \le (t-s)^{p-1} (p(2p-1))^p \int_s^t \E[|f(s)|^{2p}] ds $$
    we observe
    \begin{align*}
    \E(|H_t - H_{s}|^{2p}) &= \E\left[\left|\int_{s}^t k_H U_H(\tau) \dtau + \int_{s}^t \sigma_H(\tau) H_{\tau} dW_\tau\right|^{2p}\right] \\
        &\leq 2^{2p-1} \left\{ \E\left[\left|\int_{s}^t k_H U_H(\tau) \dtau\right|^{2p}\right] + \E\left[\left|\int_{s}^t \sigma_H(\tau) H_{\tau} dW_\tau\right|^{2p}\right] \right\} \\
        &\leq 2^{2p-1} (t-s)^{2p-1} \int_{s}^t \E[|k_H U_H(\tau)|^{2p}] \dtau + 2^{2p-1} (t-s)^{p-1} \left(p(2p-1)\right)^{p} \int_{s}^t \E[|\sigma_H(\tau) H_{\tau}|^{2p}] \dtau \\
        &= 2^{2p-1} (t-s)^{p-1} \Bigg(\int_{s}^t T^p \E[|k_H U_H(\tau)|^{2p}] \dtau + \left(p(2p-1)\right)^{p} \int_{s}^t \E[|\sigma_H(\tau) H_{\tau}|^{2p}] \dtau \Bigg) \\
        &\le 2^{2p-1} (t-s)^{p-1} \Bigg( T^p \E[\|k_H U_H\|^{2p}_{\infty}] (t-s) + \left(p(2p-1)\right)^{p} \|\sigma_H\|^{2p}_{\infty} \int_{s}^t \E[|H_{\tau}|^{2p}] \dtau \Bigg) \\
        &\le 2^{2p-1} (t-s)^{p-1} \Bigg( T^p \E[\|k_H U_H\|^{2p}_{\infty}] (t-s) + \left(p(2p-1)\right)^{p} \const^{2p}_{\sigma_H} \int_{s}^t \E[|H_{\tau}|^{2p}] \dtau \Bigg) \\[-6ex]
    \end{align*}
    \begin{align*}
        &\le 2^{2p-1} (t-s)^{p-1} \Bigg( T^p\const^{2p}_{U_H} (t-s) + \left(p(2p-1)\right)^{p} \const^{2p}_{\sigma_H} \int_{s}^t e^{\const_{p,\psi^h} \tau} \bigg( \const^{2p}_{H_{_0}}  + \const_{p,\psi^h} \const^{p}_{U_H} ~ \tau \bigg) \dtau \Bigg) \\
        &\le 2^{2p-1} (t-s)^{p-1} \Bigg( T^p\const^{2p}_{U_H} (t-s) + \left(p(2p-1)\right)^{p} \const^{2p}_{\sigma_H} \bigg(\const^{2p}_{H_{_0}} + \const_{p,\psi^h}\const^{2p}_{U_H}(t-s) \bigg) \bigg(e^{\const_{p,\psi^h}(t-s) } - 1 \bigg) \Bigg) \\
        &\le 2^{2p-1} (t-s)^{p-1} \Bigg( T^p\const^{2p}_{U_H} (t-s) + \left(p(2p-1)\right)^{p} \const^{2p}_{\sigma_H} \bigg( \const^{2p}_{H_{_0}} + \const_{p,\psi^h} \const^{2p}_{U_H}(t-s) \bigg) \bigg(e^{\const_{p,\psi^h} (t-s)}\bigg)  (t-s)  \Bigg) \\
        &\le 2^{2p-1} (t-s)^{p} \Bigg( T^p\const^{2p}_{U_H} + \const^{2p}_{p,H_0,\sigma_H, U_H} e^{\const_{p,\psi^h}(t-s)} \Bigg), \const^{2p}_{p,H_0,\sigma_H, U_H} := \left(p(2p-1)\right)^{p} \const^{2p}_{\sigma_H} \bigg( \const^{2p}_{H_{_0}} + \const_{p,\psi^h} \const^{2p}_{U_H}(t-s) \bigg).
    \end{align*}
    Thus, for any $p > 1$, by applying Kolmogorov-Chentsov's law, we get that $H_t$ has a $\vartheta$-H\"older continuous version with $\vartheta < (p-1)/2p$.
\end{proof}
\subsubsection{Analysis of $[\text{Ca}^{2+}]$ dynamics}
\begin{lem}
\label{lem:Q_dyn_wlpd}
    Let $H \in L^{p}(\Om;L^{p}([0,T];\R))$ be a given $\F_t$ measurable process that solves \eqref{eq:pH_sde} and let $\tilde P := K_{\rm sp} P(H) U_r + \tilde k_v$, then the dynamics of $Q$ specified by 
    \begin{align}
    \label{eq:Q_sde_rep}
    \begin{aligned}
        dQ(t) &= -\tilde P(t) Q(t) \dt + \sigma_Q Q(t) dW^2_t, \quad  t > 0 \\
        Q(0) &= Q_0
    \end{aligned}
    \end{align}
    has a unique strong solution $Q \in L^{p}(\Om;L^{p}([0,T];\R))$ with $Q = (Q_t)_{t\ge0}$ such that
    \begin{align}
    \label{eq:Q_sde_ineq}
    \begin{aligned}
        Q_t &= \psiq_{0,t} Q_0, \text{ and } ~~~ \E\left[ |Q_t|^{p} \right] \le e^{\const_{p,\psi^q} t} \const^p_{Q_{_0}} 
    \end{aligned}
    \end{align}
    where $\const_{p,\psi^q} := p(2p-1)\const_{\sigma_Q} + 2p \const_{\tilde P}$  \text{ as per \eqref{eq:sg_mom_est_upp}} and $\E[\|Q_0\|^{p}] \le \const^p_{Q_{_0}}$. Furthermore, $Q$ has a $\vartheta$-H\"older continuous modification for $\vartheta < 1/4$. 
\end{lem}
\begin{proof}
    Given $H \in L^p(\Om;L^p([0,T];\R))$, $(\Psi^{Q}_{s,t})_{0\le s \le t}$ is fully determined as per \eqref{eq:sg_kin_eqns}. Thus, using \eqref{eq:Exp_OU_proc} and \eqref{eq:Exp_OU_proc_solved} we get that solution to \eqref{eq:Q_sde} is uniquely given by
    $$\Phi^{Q}_t Q_0 := Q_t = \psiq_{0,t} Q_0.$$ 
    Thus,
    \begin{align}
    \begin{aligned}
        |Q_t|^{p} &\le | \psiq_{0,t} Q_0 |^{p} \\
        \implies \E\left[ |Q_t|^{p} \right]   &\le  \E \left[ | \psiq_{0,t} Q_0 |^{p} \right ] \le  \E[|\psiq_{0,t}|^{p}] ~  \E[|Q_0|^{p}] \\
        &\le e^{\left(p(2p-1)\const_{\sigma_Q} + 2p \const_{\tilde P} \right)t} \const^p_{Q_{_0}} = \const^p_{Q_{_0}} e^{\const_{p,\psi^q} t} 
    \end{aligned}
    \end{align}
    \noindent Taking supremum of $|Q_t|$ we obtain
    \begin{align}
    \begin{aligned}
        \left(\sup_{t\in[0,T]}|Q_t|\right)^{p} &\le \left(\sup_{t\in[0,T]} | \psiq_{0,t} Q_0 |\right)^{p} \le \left(\sup_{t\in[0,T]} | \psiq_{0,t}|\right)^{p} |Q_0 |^{p}\\
        \implies \E\left[\left(\sup_{t\in[0,T]}|Q_t|\right)^{p}\right] &\le \E\left[\left(\sup_{t\in[0,T]} |\psiq_{0,t}| \right)^{p}\right] \E\left[|Q_0 |^{p}\right] \le \const_{Q_{_0}} \left({p \over p -1}\right)^{p} e^{T\const_{p,\psi^q}}. 
    \end{aligned}
    \end{align}
    Next, by integrating $|Q_t|^p$ over $[0,T]$ we get
    \begin{align}
        \begin{aligned}
            \int_0^T |Q_t|^{p} \dt &\le  \int_0^T | \psiq_{0,t} Q_0 |^{p} \dt \\
        \implies \E\left[ \int_0^T |Q_t|^{p} \dt \right]   &\le  \E \left[ \int_0^T | \psiq_{0,t} Q_0 |^{p} \dt \right ] \le  \int_0^T \E[|\psiq_{0,t}|^{p}] ~ \E[|Q_0|^{p}]  \\
         &\le \int_0^T \E[|\psiq_{0,t}|^{p}] ~  \E[|Q_0|^{p}]  \le  T~\const_{p,\psi^q}~\const^p_{Q_{_0}}~\exp(T~\const_{p,\psi^q})  
        \end{aligned}
    \end{align}
    Now using SDE relation \eqref{eq:Q_sde}, we observe
    \begin{align*}
    \E(|Q_t - Q_{s}|^{2p}) &= \E\left[\left|\int_{s}^t -[\tilde P(\tau)] Q(\tau) \dtau + \int_{s}^t \sigma_Q(\tau) Q_{\tau} dW_\tau\right|^{2p}\right] \\
        &\leq 2^{2p-1} \left\{ \E\left[\left|\int_{s}^t -[\tilde P(\tau)] Q(\tau) \dtau\right|^{2p}\right] + \E\left[\left|\int_{s}^t \sigma_Q(\tau) Q_{\tau} dW_\tau\right|^{2p}\right] \right\} \\
        &\leq 2^{2p-1} (t-s)^{2p-1} \int_{s}^t \E[|[\tilde P(\tau)] Q(\tau)|^{p}] \dtau + 2^{2p-1} (t-s)^{p-1} \left(p(2p-1)\right)^{p} \int_{s}^t \E[|\sigma_Q(\tau) Q_{\tau}|^{p}] \dtau \\
        &\le 2^{2p-1} (t-s)^{p-1} \Bigg( \|[\tilde P(\tau)]\|^p_{\infty} \int_{s}^t \E[|Q_{\tau}|^{p}] + \left(p(2p-1)\right)^{p} \|\sigma_Q\|^p_{\infty} \int_{s}^t \E[|Q_{\tau}|^{p}] \dtau \Bigg) \\
        &\le 2^{2p-1} (t-s)^{p-1} \Bigg( T^p\const^p_{\tilde P}  + \left(p(2p-1)\right)^{p} \const^p_Q \Bigg) \int_{s}^t \E[|Q_{\tau}|^{p}] \\
        &\le 2^{2p-1} (t-s)^{p-1} \Bigg( T^p\const^p_{\tilde P}  + \left(p(2p-1)\right)^{p} \const^p_Q \Bigg) \const^p_{Q_{_0}} \int_{s}^t e^{\const_{p,\psi^q}\tau}   \\
        &\le 2^{2p-1} (t-s)^{p-1} \Bigg( T^p\const^p_{\tilde P}  + \left(p(2p-1)\right)^{p} \const^p_Q \Bigg) \const^p_{Q_{_0}} \left [ e^{\const_{p,\psi^q} (t-s)} - 1 \right ]  \\
        &\le 2^{2p-1} (t-s)^{p} \Bigg( T^p\const^p_{\tilde P}  + \left(p(2p-1)\right)^{p} \const^p_Q \Bigg) \const^p_{Q_{_0}} \const_{p,\psi^q} \left [ e^{\const_{p,\psi^q} (t-s)} \right ]. 
    \end{align*}
    Thus, for any $p > 1$, by applying Kolmogorov-Chentsov's law, we get that $H_t$ has a $\vartheta$-H\"older continuous version with $\vartheta < (p-1)/2p$.
\end{proof}
\subsubsection{Analysis of CaCO$_3$ dynamics}
\begin{lem}
\label{lem:C_dyn_wlpd}%
    Let $F \in \cnt([0,T];\W^{1,p}(\dom;\R))$ as per Lemma \ref{lem:mom_lin_op} and let $b^F$ be as defined in Lemma \ref{lem:decay_term_bound}. Let $Q,H \in L^{2p}(\Om;L^{2p}([0,T];\R))$ be $\F_t$ measurable process solving \eqref{eq:Q_sde} and \eqref{eq:pH_sde} respectively, and let $\sigma_C \in L^{\infty}([0,T];\R)$ and $C_0 \in L^{2p}(\Om;\R)$ then the concentration dynamics specified by 
    \begin{align}
    \label{eq:C_sde_rep}
    \begin{aligned}
        dC(t) &= [- \tilde k_v(t) C(t) + r(t) - b^F(t)]dt + \sigma_C C(t) dW^1_t, \quad t > 0 \\
        C(0) &= C_0
    \end{aligned}
    \end{align}
    has a unique strong solution $C \in L^{p}(\Om;L^{2p}([0,T];\R))$ with $C = (C_t)_{t\ge0}$ such that 
    \begin{align}
    \label{eq:pC_sde_ineq}
    \begin{aligned}
        C_t &= \psic_{0,t} C_0 + \int_0^t \psic_{s,t} \left(r(s) - b^F(s) \right) \ds \\
        \E[|C_t|^{p}] &\le  e^{\const_{p,\psi^c} t}  \left( \const_{C_{_0}} + \const_{r,b^F} ~ t \right)\\
        \E(|C_t - C_{s}|^{2p}) &\le 2^{2p-1} (t-s)^{p} \Bigg( k_C \const_u + \const_{p,\psi^c} e^{\const_{p,\psi^c}(t-s)} \Bigg)
    \end{aligned}
    \end{align}
    where
    \begin{align}
    \begin{array}{lll}
        \const_{\Psi_{C,p}} := \left(p(p-1)\const_{\sigma_C} + p \const_{\tilde k_v} \right), & \E[\|r(t)\|^{p}] \le \const^p_r & \E[\|b^F(s)\|^{p}] \le \const_{b^F}, \\[1ex]
       \E[\|U_r\|^{p}_{L^\infty_T(\R)}] \le \const^p_{U_r},~~ \E[\|U_H\|^{p}_{L^\infty_T(\R)}] \le \const^p_{U_H}, & \E[|C_0|^{p}] \le \const^p_{C_{_0}}, &  \|\sigma_C\|^{p}_{L^\infty_T(\R)} \le \const^p_{\sigma_C}.
    \end{array}
    \end{align}
    Furthermore, $C$ has a $\vartheta$-H\"older continuous modification for $\vartheta < 1/4$. 
\end{lem}
\begin{proof}
    Firstly, note that given $F$, since $b^F$ is a function of the solution $C$, \eqref{eq:C_sde} is a nonlinear SDE. To show existence we proceed by the method of successive approximation. To this end, consider the following closed subspace of $L^p(\Om; \cnt(T;\R))$
    \begin{align}
        \mathcal{Y} := \Big\{f \in L^{p}\big(\Om; \cnt(T;\R)\big): \|f_0\|_{L^{p}_{\F}(\R)} < K_2, \|f\|_{L^{p}_{\F}(\cnt(\R))} \le K_1 \Big\}
    \end{align}
    where, $K_2 < K_1$ are some arbitrary positive constants. 
    
    Now consider an arbitrary element $C^n \in \mathcal{Y}$, $n \in \N_0$ with $C^0 := C_0$, where $C_0$ is the given initial condition. Then, letting $b_n^F = b^F(C^n)$ with the convention that $b_0^F = b^F(C^0) := b^F(C_0)$, we can apply Lemma \ref{lem:decay_term_bound} to obtain that $b^F_n \in \mathcal{X}$. Consequently, given $H$ and $Q$ in $L^{p}(\Om; \cnt(T;\R))$, since $b_n^h \in \mathcal{X}$ is independent of the solution for arbitrary $C^n \in \mathcal{Y}$, we get that \eqref{eq:C_sde} is a non-autonomous and inhomogeneous SDE whose unique solution is given as
    \begin{align}
    \label{eq:nth_approx_C}
        C^{n+1} = \psic_{0,t} C_0 + \int_0^t \psic_{s,t} \left(r_n(s) - b^F_{n}(s)\right) \ds.
    \end{align}
    Thus, given $H,Q \in L^{p}_{\F}\Big(\Om;\cnt(T;\R)\Big)$, the solution \eqref{eq:nth_approx_C} defines a mapping $C^n \mapsto \Phi^C C^n := C^{n+1}$.
    Next, we show that $C^{n+1} \in \mathcal{Y}$, which then establishes that the solution operator $\Phi^C$ maps elements of $\mathcal{Y}$ into $\mathcal{Y}$ itself. 
    To this end, firstly we note that since $U_H\mapsto U_r$ and $(H,Q,U_H) \mapsto r$ are continuous functions, $U_r, r \in L^{p}(\Om;\cnt(T;\R))$. For the sake of ease of notation, in the subsequent we shall ignore the denotation of $k$. Accordingly, let 
    $$C_t = \psic_{0,t} C_0  + \int_0^t \psic_{s,t} \left(r_n(s) - b^F(s)\right) \ds, $$ 
    then since $b^F > 0$ a.s. we obtain that
    \begin{align*}
         C_t &\le \psic_{0,t} C_0 + \int_0^t \psic_{s,t} r(s) \ds, \quad \text{a.s.} \\
         |C_t|^{p} &\le |\psic_{0,t} C_0|^{p} + \int_0^t |\psic_{s,t} r(s) |^{p} \ds, \quad \text{a.s.} \\
         \implies \E \left[ |C_t|^{p}  \right ] &\le \E \left[  |\psic_{0,t} C_0|^{p} \right ] + \E \left[\int_0^t |\psic_{s,t} r_s |^{p} ~\ds \right ], \\
         &\le  \E \left[(\psic_{0,t})^{p}\right] \E \left[ |C_0|^{p}\right]  + \left[ \left \|r(t) \right\|_{\infty}\right] \int_0^t \E \left[ (\psic_{s,t})^{p}\right] ~\ds,
    \end{align*}
    For the last inequality, the first term follows from the independence of $C_0$ from $W_t$ while the second term follows from almost sure boundedness of $\|r(t)\|_{\infty}$. Accordingly, we get 
    \begin{align*}
    \E \left[ |C_t|^{p}  \right ]  &\le  \E \left[(\psic_{0,t})^{p}\right] \const^p_{C_{_0}}  + \const^p_r \int_0^t \E \left[ (\psic_{s,t})^{p}\right] ~\ds, \le  \E \left[(\psic_{0,t})^{p}\right] \const^p_{C_{_0}} + \const^p_r \int_0^t \E \left[ (\psic_{s,t})^{p}\right] ~\ds, \\
    &\le  e^{\const_{p,\psi^c} ~ t} \const^p_{C_{_0}} + \const^p_r \int_0^t e^{\const_{p,\psi^c} ~ (t-s)} ~\ds \le  e^{\const_{p,\psi^c} ~ T} \left(\const^p_{C_{_0}} + \const^p_r T\right).
    \end{align*}
    Taking supremum of $|C_t|$ we obtain
    \begin{align*}
        \left(\sup_{t\in[0,T]}|C_t|\right)^{p} &\le \left(\sup_{t\in[0,T]} | \psic_{0,t} C_0 |\right)^{p} + \left(\sup_{t\in[0,T]} \int_0^t \left |\psic_{s,t} r(s) \right | \ds\right)^{p} \\%
        &\le \left(\sup_{t\in[0,T]} | \psic_{0,t} C_0 |\right)^{p} +   \left(\int_0^t \sup_{t\in[0,T]}\left |\psic_{s,t} r(s) \right | \ds \right)^{p} \\
        &\le \left(\sup_{t\in[0,T]} | \psic_{0,t} C_0 |\right)^{p} +  \left(\int_0^t \sup_{t\in[0,T]}\left |\psic_{0,t}\right | \sup_{s\in[0,T]}\left|r(s) \right | \ds \right)^{p}\\
        &\le \left(\sup_{t\in[0,T]} | \psic_{0,t} C_0 |\right)^{p} + \int_0^t  \left(\sup_{t\in[0,T]}\left |\psic_{0,t}\right | \right)^{p} \left(\sup_{s\in[0,T]}\left|r(s) \right |\right)^{p} \ds \\
        &\le \left(\sup_{t\in[0,T]} | \psic_{0,t}|\right)^{p} |C_0 |^{p} +   \int_0^t \left(\sup_{t\in[0,T]} | \psic_{0,t}|\right)^{p} \left\|r(s) \right \|^{p}_{\infty} \ds
    \end{align*}
    \begin{align*}
        \implies \E\left[\left(\sup_{t\in[0,T]}|C_t|\right)^{p}\right] &\le \E \left[\left(\sup_{t\in[0,T]} | \psic_{0,t}|\right)^{p}\right] \E\left[|C_0 |^{p}\right] +  \int_0^t \E\left[\left(\sup_{t\in[0,T]} | \psic_{0,t}|\right)^{p}\right] \left[\left\|r(s) \right \|^{p}_{\infty}\right]^{1\over2} \ds \\
        &=\E \left[\left(\sup_{t\in[0,T]} | \psic_{0,t}|\right)^{p}\right] \Big(\E\left[|C_0 |^{p}\right] +  \left[\left\|r(s) \right \|^{p}_{\infty}\right] T \Big) \\
        &\le \E \left[\left(\sup_{t\in[0,T]} | \psic_{0,t}|\right)^{p}\right] \Big(\const^p_{C_{_0}} +  \const^p_r T \Big) \\
        &\le \left({p \over p -1}\right)^{p} e^{T \const_{p,\psi^c}} \Big(\const^p_{C_{_0}} +  \const^p_r T \Big).
    \end{align*}
    Thus, for $\E\left[|C_0 |^{p}\right] \le \const^p_{C_{_0}} < (K_2)^{p}$, and $K^p_1 > \left( {p-1 \over p} \right)^p (K^p_2 + \const^p_r)$, we can take $T \in (0, T_{\max})$ 
    \begin{align*}
        T_{\max} = {1 \over \const_{p,\psi^c}} \ln \left[{K^p_1 \over (K^p_2 + \const^p_r)} \left( {p-1 \over p} \right)^p \right] > 0
    \end{align*}
    such that 
    \begin{align}
        \E\left[\left(\sup_{t\in[0,T]}|C_t|\right)^{p}\right] \le (K_1)^{p} \implies \|C\|_{L^{p}_{\F}(\cnt(\R))} \le K_1.
    \end{align}
    Thus, we obtain that the process $(C_t)_{t\ge0}$ is an element of $\mathcal{Y}$, which in turn implies that $\Phi^C: \mathcal{Y} \to \mathcal{Y}$. In particular, we have that $\mathcal{Y} \ni C^{n} \mapsto C^{n+1} := \Phi^C C^{n} \in \mathcal{Y}.$
    Now consider a sequence of uniformly bounded process $(C^n)_{k\in\N} \in \mathcal{Y}$ where $C^n$ is obtained via $\Phi^C$. Then, we have that
    \begin{align*}
        \|C^{n+1} - C^{n}\|^p_{L^p_{\F}(\cnt(\R))} &= \left \| \int_0^t \psic_{s,t} [b^F_n(s) - b^F_{n-1}(s)] \ds \right\|^p_{L^p_{\F}(\cnt(\R))} \le \int_0^t \left \| \psic_{s,t} [b^F_n(s) - b^F_{n-1}(s)]\right\|^p_{L^p_{\F}(\cnt(\R))} \ds \\
        &\le \int_0^t \| \psic_{s,t}\|^p_{L^{2p}} \left \|b^F_n(s) - b^F_{n-1}(s)\right\|^p_{L^{2p}} \ds \\
        &\le \| \psic_{0,T_{\max}}\|^p_{L^{2p}} \int_0^{T_{\max}} \left \|b^F_n(s) - b^F_{n-1}(s)\right\|^p_{L^{2p}} \ds \\
        &\le \| \psic_{0,T_{\max}}\|^p_{L^{2p}} \left \|b^F_n - b^F_{n-1} \right\|^p_{L^{2p}} T_{\max} \\
        &\le \| \psic_{0,T_{\max}}\|^{2p\over2}_{L^{2p}} \const_b T_{\max} \|C^n - C^{n-1}\|^{p\over2}_{L^{p}} \\
        &\le \sqrt{\left({2p \over 2p -1}\right)} e^{{T_{\max} \const_{p,\psi^c} \over 2}} ~ \const_b ~ T_{\max} ~ \|C^n - C^{n-1}\|^{p\over2}_{L^{p}}
    \end{align*}
    \noindent Thus for 
    $$T_{\max} < { \sqrt{\left(1 - {1 \over 2p}\right)} e^{-{\const_{p,\psi^c} \over 2}} \over \const_b}$$
    we get that $\Phi^C$ is a contraction in $L^p(\Om)$ with 
    $$ \|C^{n+1} - C^{n}\|_{L^p_{\F}(\cnt(\R))} \le \eta \|C^{n} - C^{n-1}_t\|^{1\over2}_{L^p_{\F}(\cnt(\R))}, \quad \eta < 1. $$
    Thus, 
    \begin{align*}
        \|C^{n+1} - C^{n}\|_{L^p_{\F}(\cnt(\R))}  &\le \eta^n \left(\|C^{1}_t - C^{0}_t\|_{L^p_{\F}(\cnt(\R))}\right)^{1\over 2n}\\
        &\le \eta^n (2\const_{C})^{1\over 2n}.
    \end{align*}
    Since $(2\const_{C})^{1\over 2n} <1$ uniformly in $n$ and $\eta^n \to 0$ as $n \to \infty$ we get that $C^n \to C$ in $L^p_{\F}(\cnt(\R))$.     
    Consequently, we obtain that given $Q$ and $H$ there is a unique solution in $\mathcal{Y}$ given by
    $$(Q,H) \mapsto \Phi^C_t C_0 := C_t = \psic_{0,t} C_0 + \int_0^t \psic_{s,t} \left(r(s) -  b^F(s)\right)  \ds.$$
\end{proof}

\noindent Altogether, we combine the Lemmas \ref{lem:ph_dyn_wlpd}, \ref{lem:Q_dyn_wlpd} and \ref{lem:ph_dyn_wlpd} and formulate the following result. 
\newcommand{\phikin}{\Phi^U_H}
\begin{thm} 
\label{thm:exist_uniq_psd}
Let $U_H \in \mathcal{U}$ be a given $\F_t$ measurable process. Let $\sigma_H, \sigma_Q, \sigma_C \in L^{\infty}([0,T];\R)$ and $H_0, Q_0, C_0 \in L^{p}(\Om;\R)$, $F_0 \in \W^{2,p}(\dom)$ be independent of $\F_0$. Then there exists a finite time interval $[0,T]$ for which the system \eqref{eq:noisy_kin_sys} has a unique $\F_t$ measurable solution $\Big((H_t, Q_t, C_t, F_t)\Big)_{0\le t \le T}$ given as
\begin{subequations}
\label{eq:sol_noisy_sys}
\begin{align}
U_H \mapsto \Phi^{H}_t H_0 := H_t &= \psih_{0,t} H_0 + \int_0^t \psih_{s,t} U_H(s) \ds \label{eq:sol_ph}\\
H \mapsto \Phi^{Q}_t Q_0 := Q_t &= \psiq_{0,t} Q_0 \label{eq:sol_ca}\\
H \mapsto (\Phi^C_t C_0, \Psi^F_t F_0) &:= (C_t, F_t) \label{eq:sol_caco} \\
C_t &= \psic_{0,t} C_0 + \int_0^t \psic_{s,t} (r(s) - b^F(s)) \ds \\ 
F_t &= \psif_{0,t} F_0 
\end{align}
\end{subequations}
where $\psic$, $\psiq$ and $\psih$ are as in \eqref{eq:sg_kin_eqns} and $\psif$ is as per \eqref{eq:sol_psd_sg}.  In particular, for fixed time interval $[0,T]$, there exists a solution mapping 
\begin{align*}
\phikin&: \mathcal{U} \to \mathbf{X}, \quad U_H \mapsto \phikin U_H := \Big((H_t, Q_t, C_t, F_t)\Big)_{0\le t \le T} \\
\mathbf{X} := \Big(\mathbf{X}_1 \times \mathbf{X}_1 \times \mathbf{X}_1 \times &\mathbf{X}_2\Big), \quad \mathbf{X}_1 := L^p_{\F}\Big(\Om; \cnt([0,T];\R)\Big), \quad \mathbf{X}_2 := L^p_{\F}\big(\Om; \cnt([0,T];W^{2,p}(\dom))\big)
\end{align*}
\end{thm}
\begin{proof}
    From Lemma \ref{lem:C_dyn_wlpd} we have that the sequence $(C^n)_{n \in \N}$ as per $\eqref{eq:nth_approx_C}$ converges in $L^p_{\F}(\cnt(\R))$. Since the solution to \eqref{eq:psd_iso} is dependent on $C$ and $H$, we can construct a sequence $(F^n)_{n\in \N}$ such that
    $F^n := \Psi^{F(C^n,H)}_{0,t} F_0$, then due to \eqref{eq:lip_psd} we have that $F^n \to F := \Psi^{F(C,H)}_{0,t} F_0$ in $L^p_{\F}(\cnt(T;\W^{2,p}))$ for some $\tilde T_{\max} < T_{\max}$. Thus taking $T = \min(\tilde T_{\max}, T_{\max})$ we get the existence of the solution 
    $$\Big((H_t, Q_t, C_t, F_t)\Big)_{0\le t \le T} \text{ in } \Big(\mathbf{X}_1 \times \mathbf{X}_1 \times \mathbf{X}_1 \times \mathbf{X}_2\Big).$$
\end{proof}

\noindent Next we provide sufficient conditions for the stationarity of the solution. To this end, we note that the above local existence result can be extended by the method of pathwise continuation to any finite time interval $[0,T] \subset (0,\infty)$. 

\begin{thm} 
\label{thm:stable_soln}
Let $k_v(t)$ be a positive function such that for some $0 < t^* <  \infty$, $k_v(t) = 0$ for $t \ge t^*$, let $U_H(t)$ be an input function such that $U_H(t) \to 0$ as $t^* \le t \to \infty$ and let the mapping $U_H \mapsto U_r$ be such that $U_r \ge 0$ for $C < \csat$ and $U_r \to U_r^*$ as $U_H \to 0$, with $U_r^* > 0$, then the solution to the $X_t = (F_t, C_t, Q_t, H_t)^{\top}$ satisfies the following properties 
\begin{enumerate}
    \item $Q_t \to 0$ almost surely as $t \to \infty$ for $\sigma^2_{Q} < K_{\rm sp} P^* U_r^*$
    \item $C_t \to \csat^*$ as $t \to \infty$ in a pathwise manner when $Q_0 \gg 0$ large enough, where $\csat^* = C_{\rm sat}(H^*)$ and $H^* = \Phi^H_{t^*} H_0$. If $Q_0$ is large enough then $C_t \to C^*$ as $t \to \infty$ with $C^* \le \csat^*$.
    \item $F_t \to F_{t^*}$ as $C_t \to \csat^*$ in a pathwise manner.
\end{enumerate}
\end{thm}
\begin{proof}
    Since $U_H \to 0$ as $ t \to \infty$, from \eqref{eq:sol_ph}, we have that $H_t \to H^*$ a.s., where
    $$ H^* = H_{t^*} = \psih_{0,t^*} H_0 + \int_0^{t^*} \psih_{s,t} U_H(s) ~ \ds, \quad t^* = \min\{t > 0: U_H(t) = 0\}.$$
    Next, recalling that $ \tilde P := \Ksp P U_r + \tilde k_v$, we observe the following. Since $P > 0$ for all $H \ge 0$, $\tilde k_v(t) \to 0$ and $U_r \to U^*_r$ as $t \to \infty$, with $U_r^* > 0$ a.s., we obtain that $\tilde P \to \tilde P^*$ with $\tilde P^* > 0$ almost surely. 
    \noindent Consequently, we have that $\E[\psiq_t] \to 0$ as $ t \to \infty.$ Furthermore, if $\tilde P^* > \sigma^2_Q$ we also get that $\psiq_t \to 0$ almost surely. Consequently, we also get that $Q_t \to 0$ and $r \to 0$ almost surely as $t \to \infty$. In particular, there exists $t^* \le t^{**} < \infty$ such that $r(t^{**}) = 0$ almost surely.\\
    \noindent Now, let 
    $$I_{r_t} := \int_0^t \psic_{s,t} r(s) \ds  \quad \text{ and } \quad I_{b^F_t} :=  \int_0^t \psic_{s,t} b^F(s) \ds.$$
    Since, $b^F_t = 0$ for $C(t) \le C_{\rm sat}(t)$, without loss of generality we can only consider the case $C_0 = 0$, i.e. the edge case when $C_t$ starts with the lowest possible value, as since $C_t \ge 0$ almost surely. 
    Furthermore, since $C \in \cnt([0,T],\R)$, it is sufficient to consider pathwise analysis.
    
    \noindent Accordingly, since $U_r \ge 0$ for $C < \csat$ implies $r \ge 0$, we obtain that $C \ge 0$ for $t > 0$. 
    Thus, for $Q_0 \gg 0$ (so that $Q_t > 0$ for $t>0$), there exists $t > 0$ such that $I_{r_t} \ge \csat$. In particular, as $r \to 0$ for $t \to t^{**}$ we have that 
    $$ I_{r_t} \to I^* \: \: { \rm a.s.}  \text{ as }  t \to t^{**}, \text{\rm with } I^* \ge C^*_{\rm sat}.$$
    Consequently, $b^{F(C_t)} > 0$ for all $t > t^{**}$ and  (since $C_0 = 0$) 
    \begin{align*}
        C_t = I^* - I_{b^F_t}  \ge 0  \: \: \text{ \rm for } t > t^{**}.
    \end{align*}
    As a consequence, we get that $C_t \to C^*_{\rm sat}$ as $t \to \infty$ in a pathwise manner. Indeed, since $C_t$ keeps decreasing until $0 < b^{F(C_t)}$ tends to $0$ as $0 < C_t \to C^*_{\rm sat}$, where $C^*_{\rm sat} = C_{\rm sat}(H^*)$ then satisfies the implicit relation
    $$ C^*_{\rm sat} = I^* - I_{b^{F(C^*_{\rm sat})}}, \quad I_{b^{F(C^*_{\rm sat})}} :=  \int_0^{t^{***}} \psic_{s,t} b^F(s) \ds,$$
    where $t^{***} > t^{**}$ is such that $b^F(t) = 0$ for $t > t^{***}$. In case $Q_0 > 0$ is not large enough so that there is enough time for $C_t$ to get larger than $\csat$, then $b_t^F = 0$ and $r_t \to 0$ as $Q_t \to 0$ so $C_t \to C^*$ where $C^*$ is not necessarily equal to $C^*_{\rm sat}.$

    \noindent Finally, since $b^F \to 0$ as $C_t \to C^*_{\rm sat}$, it in turn implies $a_t \to 0$ and $N_t \to 0$ as $C_t \to C^*_{\rm sat}$, a consequence of which we immediately obtain that $F_t \to F^*$ as $C_t \to C^*_{\rm sat}.$
\end{proof}

\begin{rem}
\begin{enumerate} 
    \item The assumption about $k_v(t)$ implies that there is only an addition of volume to the initial solution which happens for a limited time i.e. $[0,t^*]$ after which there is no change in the volume. This simply says that total volume during the experimentation is bounded and does not keep increasing. 
    \item The assumption about $U_H(t)$ implies that the changes to the solution pH also ceases shortly after the changes to the volume become zero. This is simply because pH changes are induced by volume addition so if the latter stops the former has to stop as well.
    \item The assumption about $U_r \ge 0$ for $C < C_{\rm sat}$ simply states that when the system is undersaturated the kinetics should proceed in the direction of $\caco$(aq) formation which implies that $U_r$ needs to be non-negative

    \item The assumption about $Q_0$ being sufficiently large enough means that there is sufficient initial $\caion$ so that upon conversion to $\caco$(aq) the system becomes oversaturted which then promotes solid formation.

    \item Under the above nominal oversaturation condition, the theorem states that the stationary point $U^*_r$ of $U_r$ needs to be positive, in particular greater that $\sigma^2_{Q}/(\Ksp P^*)$ for the solution to fully exhaust the free calcium ions. 
\end{enumerate}
\end{rem}

\subsection{Optimal Control Formulation \label{sec:ocf}}
In this section we are concerned with the formulation of an optimal control problem for the task of model identification as well as for autonomously controlling the precipitation process. Thus, in this section the control variable is in general denoted as $U \in \{U_H, U_r\}$. For the case of model identification $U$ takes the role of a design variable $U_r$ while for the case of autonomous control, $U$ takes the role of a control variable $U_H$. To this end, we first recall the precipitation model as given in \eqref{eq:noisy_psd}. 
Accordingly, the system of state equations consisting of a deterministic PDE for the PSD \( F(t,x) \) coupled with the SDEs for the kinetic reactions is abstractly written as:
\begin{align}
\label{eq:abs_psd_spde}
\begin{aligned}
&\dd X_t = \mu(X_t,U_t;\theta) \dt + \Sigma_t \dWt{}\\
&X(0) = X_0, \\
&\begin{array}{l l}
& X_t := [H_t, Q_t, C_t, F_t]^{\top},\quad \Sigma_t := \mat{\bs{\sigma}(t) & {\bs{0}} \\ {\bs{0}} & {\bs{0}}}, \\[3ex] 
& \mu(X_t, U_t;\theta) := A(X_t;\theta) X_t + G(X_t,U_t;\theta), \\[2ex]
&A{(X_t;\theta)} := \smat{0 & 0 & 0 & 0 \\
0 & -\tilde P & 0 & 0\\
0 & \tilde P & {-\tilde k_v(t)} & 0 \\
0 & 0 & 0 & \left(-a_t \partial_x + N_t\right)
}\\[3.5ex]
&G(X_t,U_t;\theta) := \mat{k_H U_t, ~0, -\left({\tilde \rho(t) a_t S_t} + {\tilde \rho(t) v_{\rm nuc} N_t} \right), 0}^{\top}, \\[3ex]
&R_t = R_0 + \int_0^t U_s \: \ds,\: \bs{\sigma}_t = {\diag}([ H_t \sigma_H, Q_t \sigma_Q, C_t \sigma_C])   \\[3ex]
&\theta := [k_v(t), k_g, {\rho}, \delta, K_1, K_2, K_{\ce{[CO_2]}}]^{\top}
\end{array}
\end{aligned}
\end{align}

\noindent As mentioned in Section \ref{sec:intro} the core objective of the work is to autonomously adapt the pH of the solution so that the precipitation of calcium carbonate of specified size is maximized.  This application requirement is mathematically specified via the following cost/objective functional:
\begin{align}
\label{eq:cost_func}
\begin{aligned}
    J(U) &=  \E  \left[\int_0^T f(X_t, U_t) \dt + g(X_T) \Big| X(0) = x_0 \right],
\end{aligned}
\end{align}
where \(\alpha > 0\) is a regularization parameter and $X = [F,C,Q,H]^{\top}$ is the state vector function.
Based on this the finite horizon stochastic optimal control problem is given as
\begin{align}
\label{eq:socp}
\begin{aligned}
                U^*  &= \arg \min _{U \in \mathcal{U}} J(U) \\
\rm{s.t. }~~\dd X_t &= \mu(X_t,U_t;\theta) \dt + \Sigma_t \dWt{}\\
            X(0) &= X_0.
\end{aligned}
\end{align}
In the case of model identification 
\begin{align}
\label{eq:sys_id_feedback}
\begin{aligned}
    f(X_t, U_t) := {1 \over 2} \big( | Q_t - \bar Q_t |^2  + \alpha |U_t - \bar U_t|^2 \big) \\
    g(X_T, U_T) := {1 \over 2} |Q_T - \bar Q_T|^2.
\end{aligned}
\end{align}
In the case of autonomous control
\begin{align}
\label{eq:auto_con_feedback}
\begin{aligned}
    f(X_t, U_t) := {\alpha \over 2}  \big( |U_t - \bar U_t|^2 \big) \\
    g(X_T, U_T) := {1 \over 2} |F_T - \bar F_T|^2.
\end{aligned}
\end{align}
Since this work is mainly concerned with process modeling and process identification, the available data is mainly for the state variable $Q$ (i.e. $\caion$ ions) we shall consider the running cost and terminal cost as provided in \eqref{eq:sys_id_feedback}. Based on this we shall look for conditions that will guarantee the existence of an optimal control $U^* \in \mathcal{U}$ which then eventually enables us to solve the above stochastic optimal control problem (SOCP) \eqref{eq:socp}.
To this end we first need to define the adjoint variables and the Hamiltonian.
Let $\blam = (\lam_H, \lam_Q, \lam_C, \lam_F) \in \mathbf{X}$ be the vector functions representing the adjoint process of the state $X = (H, Q, C, F)$. Correspondingly, let $\bvsig = (\vsig_H, \vsig_Q, \vsig_C, 0) \in \mathbf{Z}$ be the covariances associated to the adjoint $\blam$, then the Hamiltonian function $\mathcal{H} : [0,T] \times \mathcal{X}_p \times \mathcal{U} \times \mathcal{X}_p \times \R^{4} \to \R$ is defined as:
\begin{equation}
\label{eq:Ham}
\begin{aligned}
    \mathcal{H}(t, X_t, U_t, \blam_t,\bvsig_t) &= \big \langle \mu(X_t,U_t), \blam_t \big \rangle_{\mathcal{X}} + \Tr[\Sigma_t ~ {\diag(\bvsig_t)}] - f(X_t,U_t)\\    
    &= \int_{\dom} a(t) \partial_x F(t,x) \lam_F(t,x) ~ \dx + N(C,H) \int_\dom F(t,x) \lambda_F(t,x) ~ \dx \\
    &\quad + \lambda_C(t) \left(r(t) -  {\tilde k_v(t) C(t) } - {\tilde{\rho}(t) S(t) a(C,H)} - {\tilde{\rho} v_{\rm nuc} N(C,H)} \right)  \\
    &\quad - \lambda_Q(t) \left(r(t) + {\tilde k_v(t) Q(t) } \right) + \lambda_H(t) k_H U    \\
    &\quad + C(t) \sigma_C \vsig_C + Q(t) \sigma_Q \vsig_Q + H(t) \sigma_H \vsig_H  - \frac{\alpha}{2} |U|^2
\end{aligned}
\end{equation}

\noindent Associated to the system \eqref{eq:abs_psd_spde}, design cost \eqref{eq:sys_id_feedback} and the Hamiltonian \eqref{eq:Ham} we define the following adjoint system for the (adjoint) variables $(\blam, \bvsig) \in \mathbf{X} \times \mathbf{Z}$
\begin{equation}
\label{eq:adj_sys}    
\begin{aligned}
 {\partial_t \lambda_F(t,x)} &= -a(t) {\partial_x \lambda_F(t,x)} - N(C_t) \lam_F + {\tilde \rho k_v S^*(1) a} \lam_C  \\
    \dd\lambda_C(t) &= -\bigg[ \int_{\dom} \lam_F(t,x) \Big( \partial_C a ~ \partial_x F(t,x) + \partial_C N ~ F(t,x) \Big) \dx \\
    &\quad + \lam_C(t) \Big( -\tilde{k}_v(t) - {\tilde \rho(t) S_t} \partial_C a - {\tilde \rho(t) v_{\rm nuc}} \partial_C N \Big)  + \sigma_C \vsig_C(t) \bigg] \dt + \vsig_C(t) \dWt{1} \\
    \dd\lambda_Q(t) &= -\left[ \left(-\lam_Q + \lam_C  \right) \partial_Q r - \lam_Q \tilde{k}_v(t) - \alpha (Q - \bar Q) + \sigma_Q \vsig_Q(t) \right] \dt + \vsig_Q(t) \dWt{2},  \\
    \dd\lambda_H(t) &= -\bigg[ \int_{\dom} \lam_F(t,x) \Big( \partial_H a ~ \partial_x F(t,x) + \partial_H N ~ F(t,x) \Big) \dx \\
    &\hspace*{2cm} + \lam_C(t) \Big( -\partial_H r - {\tilde \rho(t) S_t} \partial_H a - {\tilde \rho(t) v_{\rm nuc}} \partial_H N \Big) \\
    &\hspace*{3cm} + \left(-\lam_Q + \lam_C  \right) \partial_H r + \sigma_H \vsig_H(t) \bigg] \dt + \vsig_H(t) \dWt{3},  \\
    \lambda_F(t, 0) &= 0, \quad \blam_T = - \nabla g(X_T). 
\end{aligned}
\end{equation}

\noindent Following this the existence of an optimal control is given by the following sufficiency theorem.
\begin{thm}
\label{thm:suf_opt_con}
Let $X^*_t$, $(\blam^*_t$, $\bvsig^*_t)$ be the solutions to \eqref{eq:abs_psd_spde} and \eqref{eq:abs_psd_adj_spde}, respectively, corresponding to a given $U^* \in \mathcal{U}$. Furthermore, let (i) the function $x \mapsto g(x)$, i.e. the terminal cost function be a convex function and (ii) the function $x \mapsto \hat \Ham(x) := \underset{v \in \mathcal{U}}{\rm{sup}} \, \mathbb{E}[\Ham(t, x, v, \blam^*_t, \bvsig^*_t)]$ be a concave function for all $t \in [0,T]$.
If the Hamiltonian fulfills the following maximality condition
\begin{align}
\label{eq:ham_max_cond}
\underset{v \in \mathcal{U}}{\rm{sup}} \, \mathbb{E}[\Ham(t, X^*_t, v, \blam^*_t, \bvsig^*_t] = \mathbb{E}[\Ham(t, X^*_t, U^*_t, \blam^*_t, \bvsig^*_t)] \quad t \in [0, T]
\end{align}
then $U^*$ is an optimal solution for the problem \eqref{eq:socp}
\end{thm}
\begin{proof}
    Let \(U \in \mathcal{U}\) be an admissible control with corresponding state process \(X_t = X^{U}_t\). Then

\begin{equation*}
J(U^*) - J(U) = \E\left[{\alpha \over 2}\int_{0}^{T} |U^*|^2 - |U|^2 \dt +  \underbrace{{1 \over 2} \left(\left\|F^*_T - \bar{F} \right\|_{\dom}^2 - \left\|F_T - \bar{F} \right\|_{\dom}^2 \right)}_{\big(g(X^*_T) - g(X_T)\big)}\right].
\end{equation*}
Since \(g\) is convex, using the fact that $$g(x) - g(y) \le \nabla g(x)^{\top} (x-y) $$
\begin{align*}
\E[g(X^*_T) &- g(X_T)] \leq \E[(X^*_T - X_T)^{\top} \nabla g(X^*_T)] = -\E[(X^*_T - X_T)^{\top} \blam^*(T)] \\
\intertext{using integration by parts we get}
&= -\E\left[\int_{0}^{T} (X^*_t - X_t)^{\top} d\blam^*_t + \int_{0}^{T} (\blam^*_t)^{\top} (dX^*_t - dX_t) \dt + \int_{0}^{T} \Tr[\{\Sigma(X^*_t) - \Sigma(X_t)\} ~ \diag(\bvsig^*_t)] \dt\right] \\
&= -\E\left[\int_{0}^{T} (X^*_t - X_t)^{\top} \Big(-\nabla_{X}\Ham(t, X^*_t, U^*_t, \blam^*_t, \bvsig^*_t)\Big) \dt + \int_{0}^{T} (\blam^*_t)^{\top} \Big \{\mu(t, X^*_t, U^*_t) - \mu(t,X_t, U_t) \Big \} \dt\right. \\
&\hspace*{1cm} \left. + \int_{0}^{T} \Tr\left[\{\Sigma(X^*_t) - \Sigma(X_t)\} ~ \diag(\bvsig^*_t)\right] \dt\right].
\end{align*}
\noindent By the definition of \(\Ham\), we find
\begin{align*}
{\alpha \over 2} \E\left[\int_{0}^{T} \Big(|U^*|^2 - |U|^2 \Big) \dt\right] &= \E\left[\int_{0}^{T} \{\Ham(t, X_t, U_t, \blam^*_t, \bvsig^*_t) - \Ham(t, X^*_t, U^*_t, \blam^*_t, \bvsig^*_t)\} \dt\right. + \\
&\quad \left. \int_{0}^{T} \{\mu(t, X^*_t, U^*_t) - \mu(t,X_t, U_t)\} \blam^*_t \dt + \int_{0}^{T} \Tr[\{\Sigma(X^*_t) - \Sigma(X_t)\}^{\top} \diag(\bvsig^*_t)] \dt\right].
\end{align*}
Thus adding the two we obtain
\begin{align*}
J(U^*) - J(U) &\leq \E\bigg[\int_{0}^{T} \{\Ham(t, X_t, U_t,  \blam^*_t, \bvsig^*_t) - \Ham(t, X^*_t, U^*_t, \blam^*_t, \bvsig^*_t)  + (X^*_t - X_t)^{\top} \nabla_{X}\Ham(t, X^*_t, U^*_t, \blam^*_t, \bvsig^*_t)\} \dt\bigg] \\
&\leq \bigg[\int_{0}^{T} \E[\Ham(t, X_t, U_t,  \blam^*_t, \bvsig^*_t)] - \E[\Ham(t, X^*_t, U^*_t, \blam^*_t, \bvsig^*_t)]  + (X^*_t - X_t)^{\top} \nabla_{X}\Ham(t, X^*_t, U^*_t, \blam^*_t, \bvsig^*_t)\} \dt\bigg] \\
&= \E\left[\int_{0}^{T} \Ham(X_t) - \Ham(X^*_t) - (X_t - X^*_t)^{\top} \nabla_{X}\Ham(X^*_t) \dt \right]
\end{align*}

Since $\hat \Ham$ is concave, there exists some $p \in \partial_x \hat \Ham$ such that 
    $$h(x) := \hat \Ham(x) - \hat \Ham(X^*_t) - (x - X^*_t)^{\top} p \le 0$$ 
with $h(X^*_t) = 0$. This in turn implies that, due to \eqref{eq:ham_max_cond}, $\partial_x \Ham(t,X^*_t, U^*_t) = \partial_x \hat \Ham(X^*_t) = p$.
Consequently, we obtain
\begin{align*}
    0 \ge \E[h(X_t)] &= \E[\hat \Ham(X_t)] - \E[\hat \Ham(X^*_t)] - \E[(X_t - X^*_t)^{\top} p] \\
            &= \E[\hat \Ham(X_t)] - \E[\hat \Ham(X^*_t)] - \E[(X_t - X^*_t)^{\top} \nabla_{X} \hat \Ham(X^*_t) ] \\
            &= \E[\hat \Ham(X_t)] - \E[\hat \Ham(X^*_t)] - \E[(X_t - X^*_t)^{\top} \nabla_{X} \Ham(X^*_t) ] \\
            &= \E\left[ \sup_{v \in \mathcal{U}} \Ham(X_t, v) - \sup_{v \in \mathcal{U}} \Ham(X^*_t, v)\right] - \E[(X_t - X^*_t)^{\top} \nabla_{X} \Ham(X^*_t) ] \\
            &= \E\left[ \sup_{v \in \mathcal{U}} (\Ham(X_t, v)] - \Ham(X^*_t, v))\right] - \E[(X_t - X^*_t)^{\top} \nabla_{X} \Ham(X^*_t) ] \\
            &\ge \E[\Ham(X_t) - \Ham(X^*_t)] - \E[(X_t - X^*_t)^{\top} \nabla_{X} \Ham(X^*_t) ] \\
            &\ge \E[\Ham(X_t)] - \E[\Ham(X^*_t)] - \E[(X_t - X^*_t)^{\top} \nabla_{X} \Ham(X^*_t) ] \\
    \implies \int_0^T 0 \ge \int_0^T \E[h(X_t)] &\ge \int_0^T \E[\Ham(X_t)] - \E[\Ham(X^*_t)] - \E[(X_t - X^*_t)^{\top} \nabla_{X} \Ham(X^*_t) ]
\end{align*}
Thus we obtain that 
\begin{align*}
    J(U^*) - J(U) &\le \E\left[\int_{0}^{T} \Ham(X_t) - \Ham(X^*_t) - (X_t - X^*_t)^{\top} \nabla_{X}\Ham(X^*_t) \dt \right] \\
     &\le \int_{0}^{T} \E\left[\Ham(X_t) - \Ham(X^*_t) - (X_t - X^*_t)^{\top} \nabla_{X}\Ham(X^*_t) \dt \right] \\
     &\le 0 \\
\implies J(U^*) \le J(U).
\end{align*}
\end{proof}

\noindent Theorem \eqref{thm:suf_opt_con} emphasizes the importance of the Hamiltonian and the associated adjoint system. However, to numerical compute the optimal control $U^*$ we need a rather constructive method. Consequently, it involves solving the adjoint system \eqref{eq:adj_sys} which, along with other coupled relations results in the following necessary condition for the existence of an optimal control $U^*$. 



\noindent Firstly, similar to \eqref{eq:abs_psd_spde}, let us reformulate the adjoint system \eqref{eq:adj_sys} in the following abstract form
\begin{align}
\label{eq:abs_psd_adj_spde}
\begin{aligned}
\dd \blam_t &= -\Gamma(X_t, U_t, \blam_t, \bvsig_t;\theta) \dt + \diag(\bvsig_t) \dWt{}\\
\blam(T) &= \blam_T, \\
\blam_t &:= [\lam_{H,t}, \lam_{Q,t}, \lam_{C,t}, \lam_{F,t}]^{\top},\quad \bvsig = [\vsig_H, \vsig_Q, \vsig_C, 0]^{\top}, \\[3ex]
\blam_T &:= [0, -(Q_T - \bar{Q}), 0, 0]^{\top}, \quad \lam_F(t,0) = 0, \\[2ex]
\Gamma(X_t, U_t;\theta) & := A^*(X_t,\blam_t;\theta) \blam_t + G^*(X_t,U_t;\theta), \\[2ex]
A^*{(X_t;\theta)} & := \smat{0 & 0 & -{\tilde \rho(t) S^*(1) a} & 0 \\
0 & -\partial_Q r - {\tilde k_v(t)} & \partial_Q r & 0\\
0 & 0 & -{\tilde k_v(t)} - {\tilde \rho(t) S_t} \partial_C a - {\tilde \rho(t) v_{\rm nuc}} \partial_C N  & 0 \\
0 & 0 & 0 & \left(a_t \partial_x + N_t\right)
}, \\
G^*(X_t,U_t;\theta) &:= \mat{\Gamma_1(t), -\alpha(Q - \bar Q) + \sigma_Q\vsig_Q, ~\int_{\dom} \lam_F \Big( \partial_C a ~ \partial_x F + \partial_C N ~ F \Big) \dx + \sigma_C \vsig_C, ~0}^{\top}, \\[3ex]
\Gamma_1(t) &:= \Bigg[ \int_{\dom} \lam_F \Big( \partial_H a ~ \partial_x F + \partial_H N ~ F \Big) \dx + \lam_C \Big( -\partial_H r -{\tilde \rho(t) S_t} \partial_H a - {\tilde \rho(t) v_{\rm nuc}} \partial_H N \Big) \\
    &\hspace*{2cm}  + \left(-\lam_Q + \lam_C  \right) \partial_H r + \sigma_H \vsig_H(t) \Bigg]
\end{aligned}
\end{align}

Accordingly, we shall now formulate the necessary condition based on which the finite horizon stochastic optimal control problem \eqref{eq:socp} can be iteratively solved.

\begin{thm}[Necessary maximum principle] Let $U \in \mathcal{U}$ and let $\widetilde{U}$ be a perturbation of $U$ such that $\widetilde{U} := U + \alpha \beta$ belongs to $\mathcal{U}$ for $\alpha \in (-1,1)$. Let $\Xalp := X^{U+\alpha \beta}$, denote the process driven by $\alpha$ perturbed control $U + \alpha \beta$ and let $\x(t):=\frac{d}{d\alpha}\Xalp(t)|_{\alpha=0}$ be a derivative process whose time differential $d\x$ is an element of $L^p(\Om)$, $p \ge 2$. \\
Then 
\begin{align}
\label{eq:nec_cond}
\begin{aligned}
 &\frac{d}{d\alpha} J(U+\alpha\beta)|_{\alpha=0} = 0  \quad \forall \beta \in \mathcal{U} \text{ with } (U + \alpha\beta) \in \mathcal{U}, \\
 \hspace*{-3cm} \text{if and only if} \hspace*{2.5cm}  & \\
 &\mathbb{E} \left[ {\partial_U \Ham}(t) \right] = 0 \text{ for all } t \in [0, T].    
\end{aligned}
\end{align}
\end{thm}
\begin{proof}
Firstly, we observe that for any given $U, \beta_0 \in \mathcal{U}$ we can always choose a perturbation 
\begin{align}
\label{eq:cont_perturb}
    \beta := \delta(t) \beta_0(t) \text{ with } \delta(t):= \left(\frac{1}{2K}\text{dist}(U_t, \partial \mathcal{U}) \wedge 1\right) > 0,
\end{align}
such that $\widetilde{U} := U + \alpha \beta$ is still in $\mathcal{U}$ for all $\alpha \in (-1,1)$. Furthermore, the time differential $d\x_t$ of the derivative process $\x$ satisfies the following SDE:
\begin{align}
\label{eq:der_proc_sde}
\begin{aligned}
d\x(t) &= \Big\{ 
{\partial_X \mu(t)} \cdot {\x(t)}+ {\partial_U \mu(t)}\beta(t) 
\Big\} \dt + \Big\{ 
{\partial_X \Sigma(t)} ~ {\diag(\x(t))} + {\partial_U \Sigma(t)}\beta(t) 
\Big\} \dWt{}, \\
\x(0) &= 0,
\end{aligned} 
\end{align}
where, $\diag(\x(t)) \in \R^{d \times d}$ denotes the diagonalized matrix of the vector $\x(t) \in \R^d$, $t \in [0, T]$.
\noindent Now, let us define a sequence of stopping times $\tau_{n}; n=1,2,\dots$ as follows: 
\begin{align*}
\tau_{n} = \inf \left \{t>0;\int_{0}^{t} \Tr(\blam_s \blam^{\top}_s) \Tr(\Sigma_s \Sigma^{\top}_s) + \Tr(X_sX^{\top}_s) \Tr(\bvsig_s \bvsig^{\top}_s) \ds \ge n \right \}\wedge T.
\end{align*}
Then it is clear that 
$\tau_{n} \rightarrow T \text{ as } n \rightarrow \infty$.
Due to Martingale property, we also have that:
\begin{align*}
\mathbb{E} \left[ \int_{0}^{\tau_{n}} {\diag(\blam_t)~\Sigma_t} ~ dW(t) \right] & = \mathbb{E} \left[ \int_{0}^{\tau_{n}} {\diag(\x_t) ~ \diag(\bvsig_t)} ~ dW(t) \right] = 0 \text{ for all } n.
\end{align*}
We can write 
\[
\frac{d}{d\alpha} J(U + \alpha \beta) \big|_{\alpha=0} = I_1 + I_2,
\]
where 
\[
I_1 = \frac{d}{d\alpha} \mathbb{E} \int_0^T f \left( t, X^{U+\alpha \beta}(t), U_t + \alpha \beta(t) \right) \dt \big|_{\alpha=0},
\]
and
\[
I_2 = \frac{d}{d\alpha} \mathbb{E}[g(X^{U+\alpha \beta}(T))] \big|_{\alpha=0}.
\]
By our assumptions on \(f\) and \(g\), we have
\[
I_1 = \mathbb{E} \int_0^T \left[ \frac{\partial f}{\partial_x}(t) x(t) + \frac{\partial f}{\partial U}(t) \beta(t) \right] \dt,
\]
\[
I_2 = \mathbb{E}[g'(X_T) \cdot \x(T)] = -\mathbb{E}[\blam(T) \cdot \x(T)].
\]
By the Ito formula and \eqref{eq:der_proc_sde},
\begin{align*}
I_2 = -\mathbb{E}[\blam(T) &\cdot \x(T)] =  -\lim_{n \to \infty} \E[ \blam(\tau_n) \cdot \x(\tau_n)] \\
&= -\lim_{n \to \infty} \mathbb{E} \left[ \int_0^{\tau_n}  \blam(t)\cdot d\x(t) + \int_0^{\tau_n} \x(t) \cdot d\blam(t) + \int_0^{\tau_n} d[\blam,\x](t) \right]  \\
&= -\lim_{n \to \infty} \mathbb{E} \left[ \int_0^{\tau_n} \blam(t) \cdot \Big\{ {{\partial_X \mu(t)} \cdot \x(t)} + {\partial_U \mu(t)} \beta(t) \Big\} \dt + \int_0^{\tau_n} \left(-\partial_X \Ham \right) \cdot \x(t)  ~ \dt \right . \\
&\qquad \qquad + \left . \int_0^{\tau_n} \Tr\Big( \diag(\bvsig(t)) \Big \{ \partial_X \Sigma(t) {\diag(\x(t))} + \partial_U \Sigma(t) \beta(t)  \Big \} \Big) \dt \right]  \\
&= -\lim_{n \to \infty} \mathbb{E} \left[ \int_0^{\tau_n}  \x(t) \cdot \Big( {\partial_X \mu(t)} \cdot \blam(t) + \partial_X \Sigma \bvsig(t)  - \partial_X \Ham \Big) \dt  \right. \\
&\qquad \qquad \qquad + \left. \int_0^{\tau_n} \beta(t) \Big( \partial_U \mu(t) \cdot \blam(t) + \Tr \big(\partial_U \Sigma_t ~ \diag(\bvsig_t)\big) \Big) \dt \right]\\
&= -\lim_{n \to \infty} \mathbb{E} \left[ \int_0^{\tau_n}  \partial_X f \cdot \x(t) ~ \dt   + \int_0^{\tau_n} \beta(t) \Big( \partial_U \Ham + \partial_U f \Big) ~ \dt \right]\\
&= -I_1 - \mathbb{E} \left[ \int_0^T \partial_U \Ham ~ \beta(t) ~ \dt \right].
\end{align*}
Thus, we get:
\[
\frac{d}{d\alpha} J(U + \alpha \beta) \big|_{\alpha=0} = -\mathbb{E} \int_0^T {\partial_U \mathcal{H}}(t) \beta(t) \dt.
\]
Consequently, we conclude that
\[
\frac{d}{d\alpha} J(U + \alpha \beta) \big|_{\alpha=0} = 0
\]
if and only if 
\[
\mathbb{E} \left[ \int_0^T {\partial_U \mathcal{H}}(t) \beta(t) \dt \right] = 0, \quad \text{for all bounded } \beta \in \mathcal{U} \text{ of the form \eqref{eq:cont_perturb}}.
\]
In particular, applying this to \(\beta(t) = \indc_{[0,T]}\), we get that this is again equivalent to
\[
\mathbb{E} \left[ {\partial_U \mathcal{H}}(t) \right] = 0 \quad \text{for all } t \in [0, T].
\]
\end{proof}

The above theorem says that if there exists $U^* \in \mathcal{U}$ such that $J^* = J(U^*) = \min_{U \in \mathcal{U}} J(U)$ then $U^*$ is a stationary point of $J$ only if $\mathbb{E} \left[ {\partial_U \mathcal{H}} \right] |_{U^*} = 0.$ Thus, set of solutions to the equation $\mathbb{E} \left[ {\partial_U \mathcal{H}}\right]  = 0$ or even $\partial_U \mathcal{H}  = 0$ provides the set of candidate solutions to the problem \eqref{eq:socp}. Based on this we can develop an algorithm that uses the state-system \eqref{eq:abs_psd_spde}, the adjoint-system \eqref{eq:adj_sys} and the condition \eqref{eq:nec_cond} to iteratively obtain a sequence of controls $(U^n)_{n \in \N}$ along which $J$ is decreasing and consequently converges to a stationary value $J^*$ that corresponds to a (sub) optimal value $U^*$. 

\section{Numerical investigations \label{sec:numerics}}
\subsection{Experimental setup}
The pH control is not only necessary to induce supersaturation for CaCO$_3$ precipitation but also to steer the process so that only CaCO$_3$ is selectively precipitated from a solution with different ions dissolved in it. A semi-batch mode of operation enhances the pH control by allowing the gradual addition of CO$_2$ gas and NaOH for pH regulation, ensuring uniform mixing and avoiding localized supersaturation that can lead to undesired by-products. 

\subsection{Setup Description}
\begin{figure}[h!]
    \centering
    \includegraphics[width=0.8\textwidth]{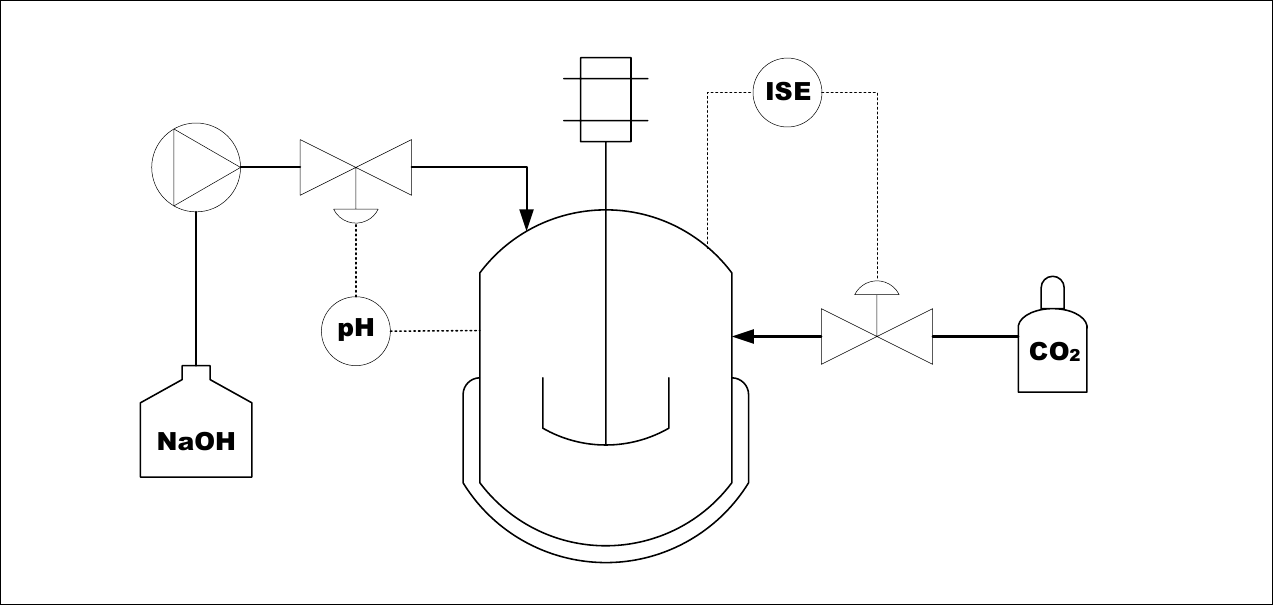}
    \caption{Semi-batch experimental setup with provision for controlled addition of NaOH and CO$_2$, monitoring pH and [Ca$^{2+}$] concentration.}
    \label{fig:setup}
\end{figure}

An experimental setup in a semi-batch mode of operation was designed, as shown in Figure \ref{fig:setup}. Provision was made for the controlled addition of NaOH solution and bubbling CO$_2$ gas while continuously monitoring the pH and [Ca$^{2+}$] concentration. CO$_2$ gas and NaOH solution were added in sequence, achieving several pH-swing cycles.

For the carbonation and selective precipitation experiments, 300~mL of ultrapure Milli-Q water at neutral pH was used. To this, Calcium and Magnesium salts (CaCl$_2$ and MgCl$_2$) were added in varying concentrations to simulate mine tailing extracts. While dissolving CaCl$_2$ kept the solution neutral, MgCl$_2$ increased the pH to 8-9, depending on its concentration.

To this ionic solution, CO$_2$ was bubbled at a steady rate of 4~L/h until saturation was reached (CO$_2$(g) Vol \% = 100\%) while continuously monitoring the pH and [Ca$^{2+}$] ion concentration. The pH of the CO$_2$-saturated solution dropped to approximately 4 regardless of dissolved Ca and Mg concentrations. Samples taken for IC measurements revealed that the concentrations of Ca and Mg remained the same after CO$_2$ addition due to negligible carbonate ion ([CO$_3^{2-}$]) formation at low pH, preventing CaCO$_3$ precipitation.

To precipitate CaCO$_3$, the pH of the solution was increased by adding 1~M NaOH at 3~mL/min until the pH reached 9, where the solubility of CaCO$_3$ is very low ($\sim10^{-8}$~mol/L). As the pH rose, the Ca$^{2+}$ concentration decreased due to the quick formation and precipitation of CaCO$_3$. In one cycle of the pH swing, only a small amount of dissolved [Ca$^{2+}$] in the solution could be precipitated (around 0.03~M) due to the limited amount of CO$_2$ that can be dissolved in the solution (0.034~M at 1~bar and 25~$^{\circ}$C). 

To completely extract most of the dissolved [Ca$^{2+}$] ions as CaCO$_3$ precipitate, the pH-swing cycle was implemented several times by alternately bubbling CO$_2$ until saturation and then supplying NaOH until pH 9. After each CO$_2$ bubbling and NaOH addition, IC samples were taken for ion concentration measurements, while continuously monitoring the pH and [Ca$^{2+}$] concentration with ion sensors. After the precipitation of the carbonates, the precipitated sample was dried and measured using XRD to determine if pure CaCO$_3$ was precipitated selectively.

\begin{figure}[!h]
\begin{center}
    \includegraphics[scale=.4]{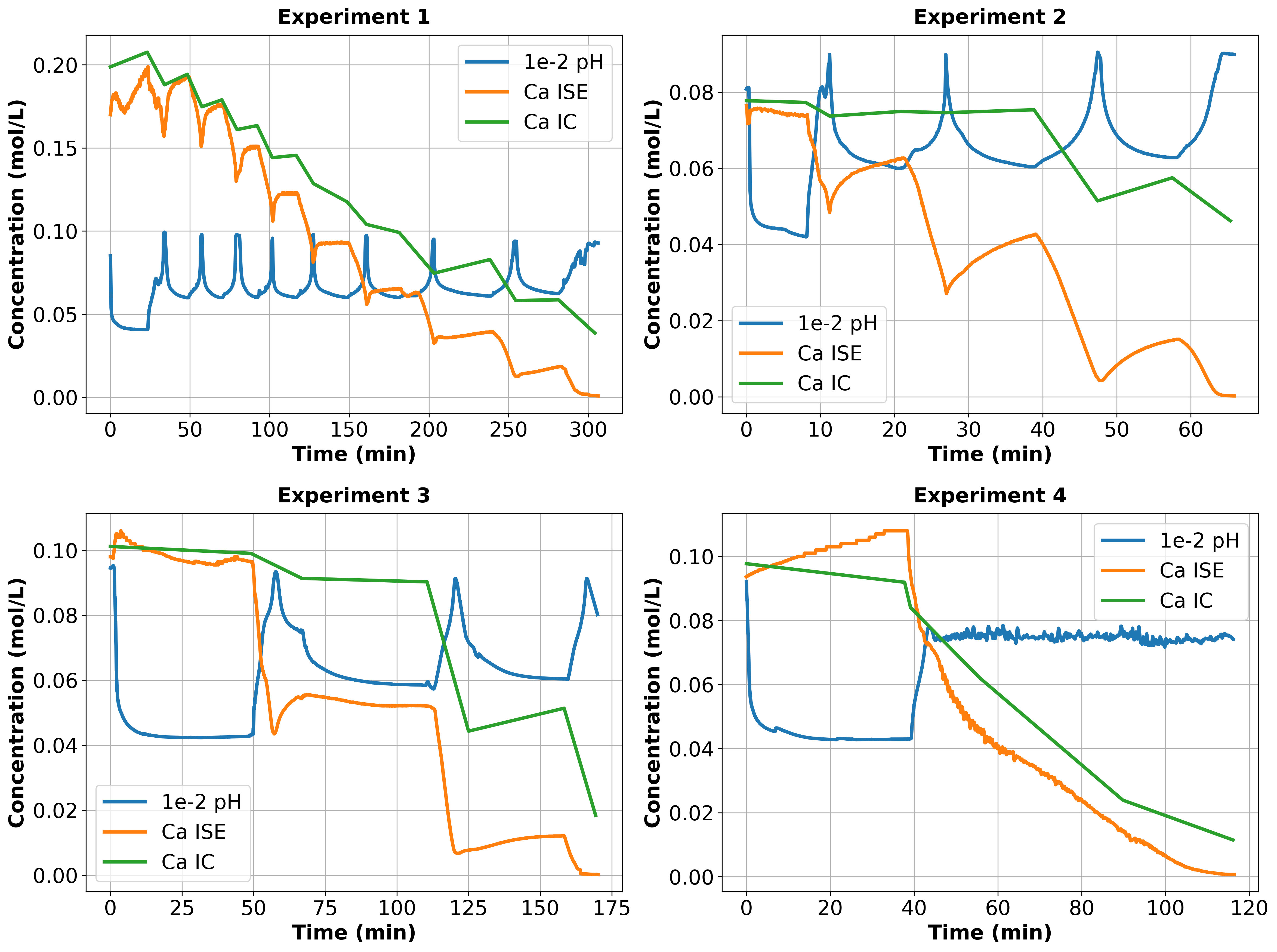}
\end{center}
\caption{Experimental measurements of $\caion$ measured using two different sensor devices, IC (green) and ISE (orange). The measurements are observed in response to the pH changes induced either by adding NaOH solution or CO$_2$ gas. The resulting  pH profile (after scaling by a factor of $0.01$) is shown in blue. \label{fig:exp_data}}
\end{figure}
\section{Simulation results \label{sec:sim_results}}
In this section we consider the simulation of the precipitation model \eqref{eq:abs_psd_spde} and compare the results with the experimentally obtained data. For this purpose, the system of equations \eqref{eq:abs_psd_spde} is solved numerically using the standard discretization method which we shall briefly describe below. 
Accordingly, the finite time domain $[0,T]$ is discretized uniformly via the stepping parameter $\tau$ in the following way:
$$t_j = t_0 + j \tau, \quad j = 0 \ldots N, \quad h = \frac{T}{M}.$$
Similarly, the size domain $[0,L]$ is discretized uniformly in the $x$-direction based on the step parameter $h$ as:
$$x_i = x_0 + i h, \quad i = 0 \ldots M, \quad h = \frac{L}{M}.$$
Accordingly, the numerical approximation of the solution component $F$ at time point $t_j$ and size value $x_i$, i.e. $F(t_j, x_i)$,  is shortly denoted as  $F_i^j$. Analogously, the numerical approximations of $C, Q, H$ at $t_j$ and $x_i$ are denoted as $C_i^j, Q_i^j$ and $H_i^j$ respectively. Since, the PSD dynamics is modeled via an advection equation, we employ the Lax-Wendroff method \cite{Lax1964}, a second-order accurate finite difference scheme for solving hyperbolic partial differential equation. Due to its stability and consistency properties, it is particularly well-suited for advection problems. The method utilizes a Taylor series expansion in time and incorporates partial derivatives with respect to the 2nd variable to achieve higher-order accuracy. Using the Lax-Wendroff scheme the discrete version $D^j_i$ of the partial differentiation operator $\partial_x$ at $(t_j,x_i)$ is given as \\[-6ex]
\begin{align*}
D_i^j F &= -{1 \over h} \left(D_{i-\frac{1}{2}}^{j+\frac{1}{2}} - D_{i+\frac{1}{2}}^{j+\frac{1}{2}}\right) F = -{1 \over h} \left(F_{i-\frac{1}{2}}^{j+\frac{1}{2}} - F_{i+\frac{1}{2}}^{j+\frac{1}{2}}\right) \\
D_{i-\frac{1}{2}}^{j+\frac{1}{2}} F = F_{i-\frac{1}{2}}^{j+\frac{1}{2}} &= \frac{1}{2} \left( F_i^j + F_{i-1}^j \right) - \frac{\tau}{2 h} \left( F_i^j - F_{i-1}^j \right), \\
D_{i+\frac{1}{2}}^{j+\frac{1}{2}} F = F_{i+\frac{1}{2}}^{j+\frac{1}{2}} &= \frac{1}{2} \left( F_i^j + F_{i+1}^j \right) - \frac{\tau}{2 h} \left( F_{i+1}^j - F_i^j \right).
\end{align*}
These equations represent the core of the Lax-Wendroff method, where the first two equations compute the solution at intermediate time and size values, and the third equation updates the solution at the next time step. The method's stability is subject to the Courant-Friedrichs-Lewy (CFL) condition, which relates the stepping parameters $\tau$, $h$ and the advection speed. In particular the CFL condition states that ${k_g \tau \over 2h} < 1,$ where $k_g$ denotes the advection speed.

Subsequently, the remaining set of equations that are described by SDEs, a simple Euler-Maruyama \cite{Gikhman2007} scheme is used. Thus, based on the Lax-Wendroff (a discretized advection) operator and Euler-Maruyama discretization, the forward simulation of the precipitation process is performed using the following scheme
\begin{align}
\label{eq:abs_disc_scheme}
\begin{aligned}
    X^{j+1}_i   &= X^{j}_i + \tau \mu^{j}_i + \sqrt{\tau} ~ \Sigma^j_i ~ {\rm z^j}, \quad \\
    X^0_i       &= X_{0,i}, \\
    {\rm z^j}   &= [z^j_1, z^j_2, z^j_3, z^j_4], \: z^j_{\ell} \in \mathcal{N}(0,1), \: \ell \in \{1,\cdots,4\} \\[-3ex]
\end{aligned}
\end{align}
where,
\begin{align}
\label{eq:disc_scheme}
\begin{aligned}
\begin{array}{l l}
    &X^j_i := [H^j, Q^j, C^j, F^j_i]^{\top},\quad \Sigma^j := \mat{\bs{\sigma}^j & {\bs{0}} \\ {\bs{0}} & \bs{0}},\quad U^j := U(t_j),\\
    &\mu^j_i   = A(X^j_i) X^j_i  + G(X^j_i, U^j_i), \\[1ex]
    &A{(X^j_i;\theta)} := \smat{0 & 0 & 0 & 0 \\
                        0 & -\tilde P^j & 0 & 0\\
                        0 & \tilde P^j & {-\tilde k_v^j} & 0 \\
                        0 & 0 & 0 & \left(-a^j D^j_i + N^j\right) \\
                        },\\[3.5ex]
    &G(X^j_j,U^j) := \mat{k_H U^j, 0, -\left(\frac{\rho k_v(t) a^j S^j}{R^j} + \frac{\rho v N^j}{R^j} \right), 0}^{\top}, \\[3ex]
    &\bs{\sigma}^j = {\diag}([H^j \sigma_H, Q^j \sigma_Q, C^j \sigma_C]),   \\[3ex] 
    &R^j = R^{j-1} + \tau k_v^j.
\end{array}
\end{aligned}
\end{align}
To validate the proposed model we aim to compare the simulation results with that of the measured data (Figure \ref{fig:exp_data}). For this purpose, we set the model parameters $\theta$, in accordance to the experimental setup, as per the values in Table \ref{tab:model_params}. Based on this, the rate of change of pH, i.e. $\const_H U_t$, is taken as a known function of time  and the rate modulation function $U_r$ is the unknown which needs to be determined. 

\begin{table}[!h]
    \vspace*{0.1cm}
    \centering
    \begin{tabular}{|c|c|c||c|c|c|}
    \hline
Constant &  Value  &  Unit & Constant &  Value  &  Unit\\\hline
$K_{sp}$ & $2.8 \times 10^{-9}$ & 1 & $K_{a_1}$ & $10^{-6}$ &  1\\\hline
$K_{a_2}$ & $10^{-10}$ & 1 & $d$ & $1 \times 10^{-6}$ &  1$\mu m$ \\\hline
$\tilde K_{\ce{[CO_2]}}$ & $100$ & 1 & $\tilde K_{a_1}$ & 0.45 & 1 \\\hline $\tilde K_{a_2}$ & $6.6$ & 1 & $\tilde K_{a_3}$ & $13.5$ & 1 \\\hline
$K_{1, \rm sat}$ & $.05$ & 1 & $K_{2, \rm sat}$ & $1.2 \cdot 10^6$ & 1 \\\hline
$v_{nuc}$ & $\frac{\pi}{6} \cdot d^3$ & $m^3$ & $K_{\rm co2}$ & $1 \times 10^{-3}$ & ${\rm mol/\ell}$ \\\hline
$\nu$ & $1.0$ & 1 & $C_{0}$ & $0$ & ${\rm mol/\ell}$\\\hline
$R_0$ & $5.18$ & $\ell$ & $\dot{R}$ & $1 \times 10^{-2}$ & $1/s$\\\hline
$k_{g}$ & $0.459$ & ${\mu m / s}$ & $\rho$ & $0.315$ & $g/cm^3$\\\hline
    \end{tabular}
    \vspace*{0.01cm}
    \caption{Fixed parameters of the model chosen as per the experimental setup.}
    \label{tab:model_params}
    \vspace*{-.3cm}
\end{table}
In the following, we consider three different approaches to obtain the most suitable function $U_r$ which eventually is able to show that the proposed model is able to correspond to the experimental data. Since $U_H$ is a known input function of time, without loss of generality, $\const_H$ is set to 1.

\subsection{Manual fitting of the function}
\begin{figure}[!h]
\begin{center}
    \includegraphics[scale=.24,trim = 0 0 0 0,clip]{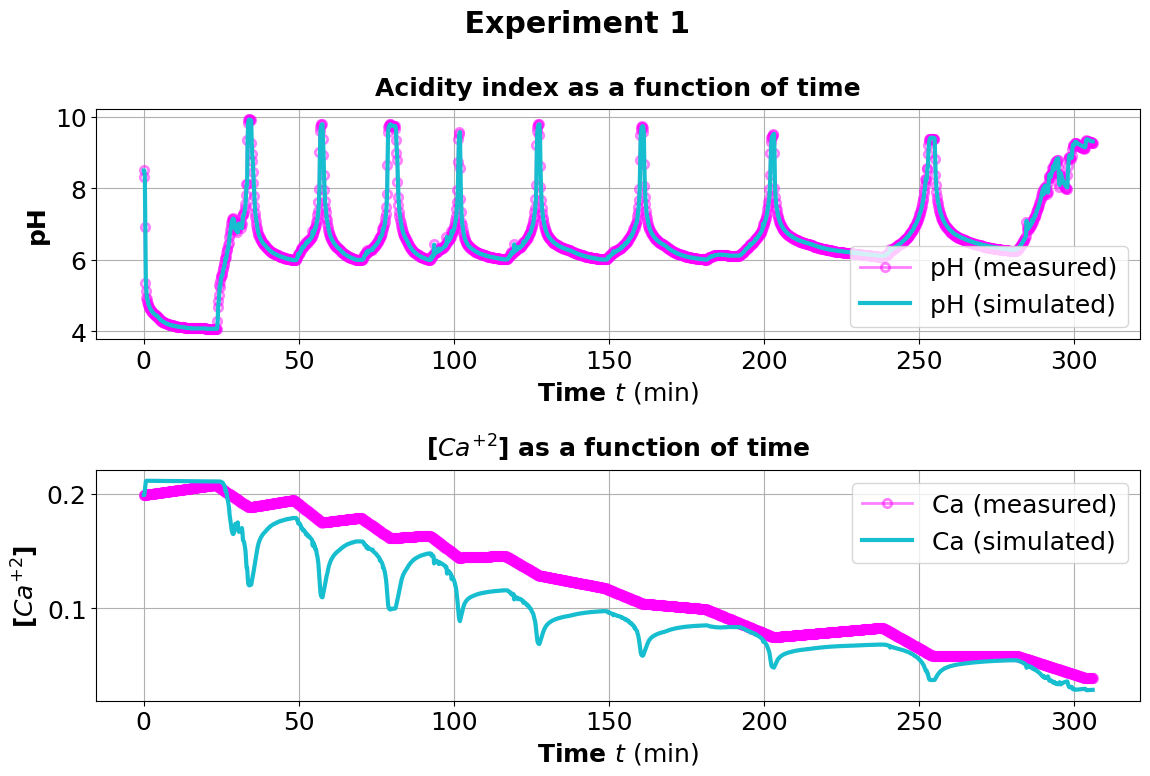} \qquad
    \includegraphics[scale=.24,trim = 0 0 0 0,clip]{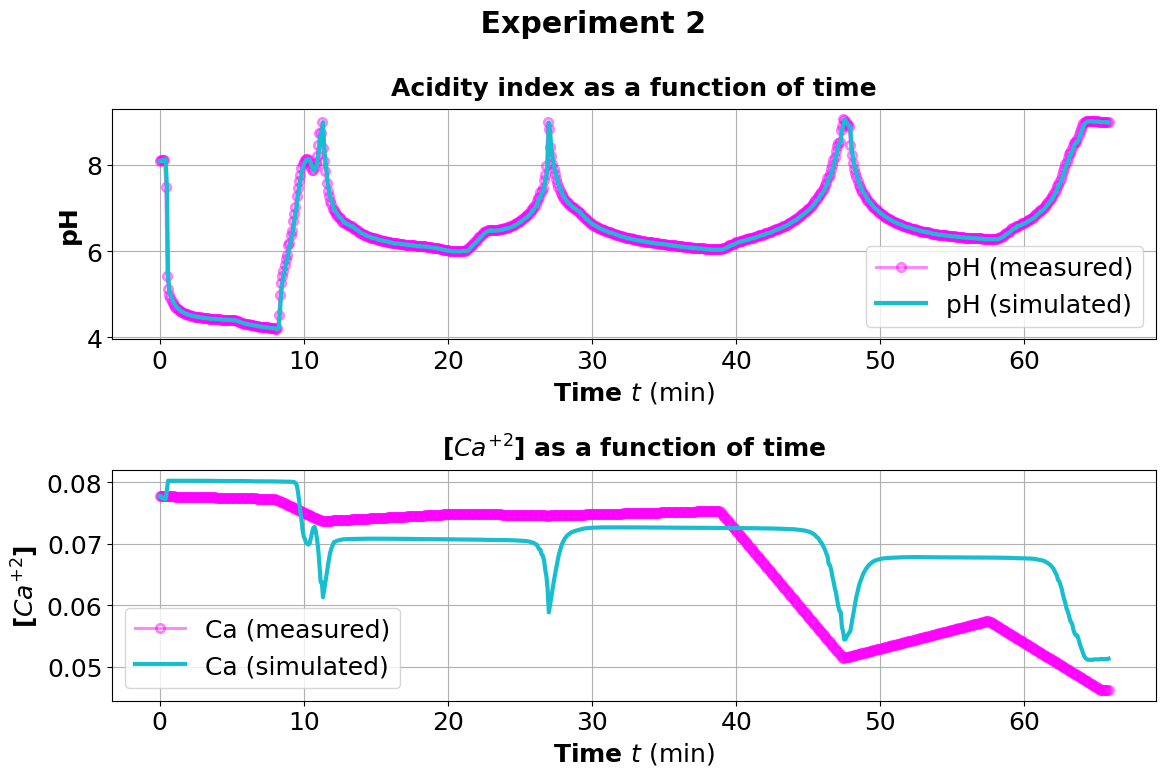} \\
    \includegraphics[scale=.24,trim = 0 0 0 0,clip]{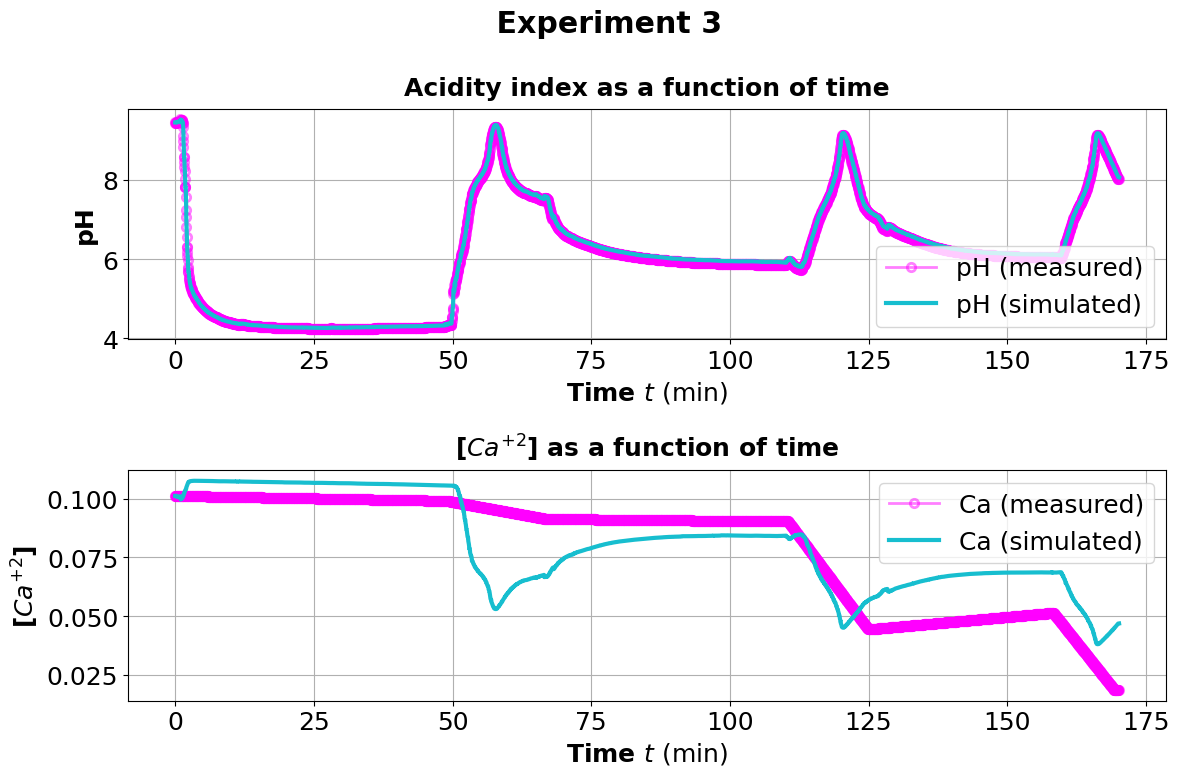} \qquad
    \includegraphics[scale=.24,trim = 0 0 0 0,clip]{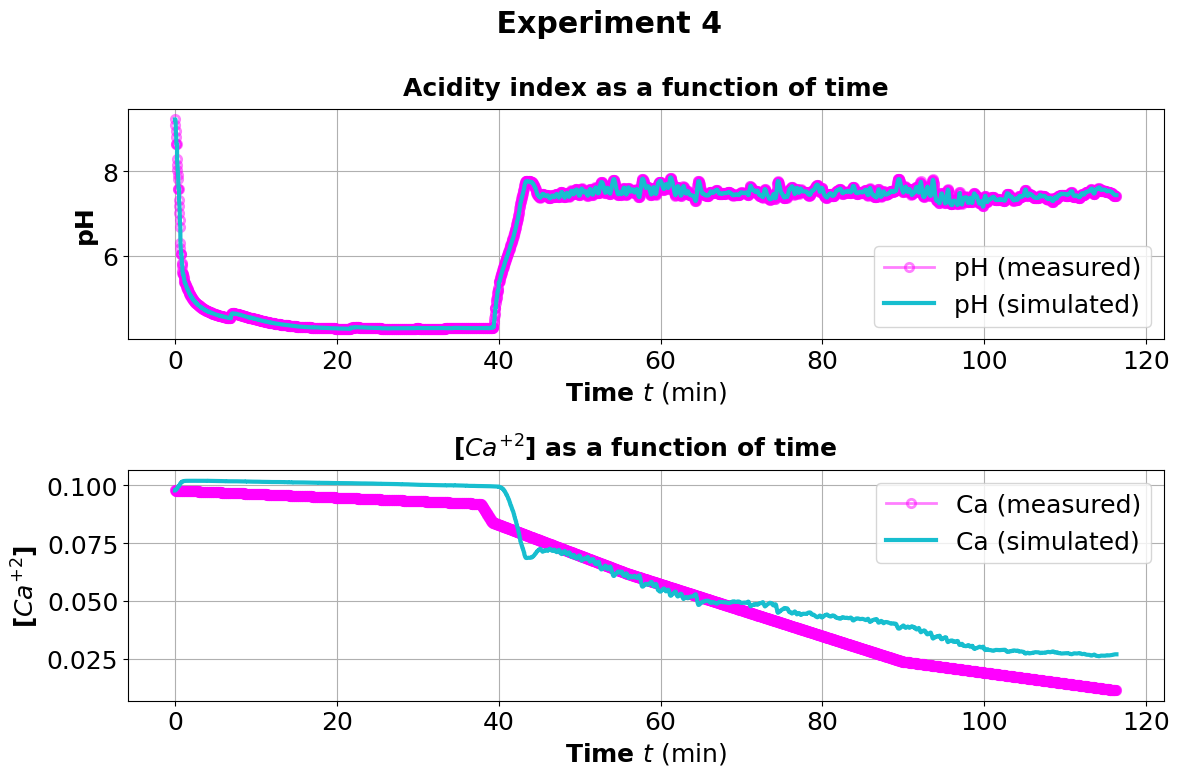} 
\end{center}
\caption{\label{fig:man_fit} Concentration of calcium ions as per the manually determined function $U_r$ with experiment specific constants as mentioned in Table \ref{tab:man_fit_par}.}
\end{figure}

As a first approach we follow a brute force and intuitive method to model the function $U_r$. Based on the measured observations (cf. Figure \ref{fig:exp_data}) we see that whenever there is a spike in the pH, i.e. a sudden increase or decrease in the pH, the $\caion$ ions also exhibits a corresponding change in its trend i.e. a sudden decreasing or increasing trend respectively.  Motivated from this we represent $U_r$ to be an explicitly dependent function of $U_H$. In particular, we define $U_r$ as 
\begin{align}
\label{eq:man_ur}
\begin{aligned}
    U_r(t) = \begin{cases}
            k^+_{uc}        & \text{if }  \quad k_{rc} U_H > k^+_{uc} \\
            k^-_{uc}        & \text{if }  \quad k_{rc} U_H < k^-_{uc} \\
            k_{rc} U_H      & \text{if }  \quad k_{rc} U_H \in [k^-_{uc}, k^+_{uc}] 
        \end{cases}\\
\end{aligned}
\end{align}
The parameters of the function $U_r$ are assembled as the vector $\kappa = [k_{rc}, k^-_{uc}, k^+_{uc}]^{\top}$. The values of the components of the vector $\kappa$ are determined by trail and error guided by dynamical intuition. Accordingly, thus obtained values are as shown in Table \ref{tab:man_fit_par}.
\begin{table}[!h]
    \centering
    \begin{tabular}{|c|c|c|c|c|}
         \hline
         Exps&  $k_{rc}$ & $k^{-}_{uc}$& $k^+_{uc}$ \\\hline
         1&  0.175&  -0.02&  .042 \\\hline
         2&  0.15&  -0.02&  .06 \\\hline
         3&  0.19&  -0.003&  .19 \\\hline
         4&  0.19&  -0.002&  .09 \\\hline
    \end{tabular}
    \caption{Manually chosen constants for the function $U_r$.}
    \label{tab:man_fit_par}
\end{table}
The resulting simulated dynamics of the calcium concentration, in all four cases, are as shown in Figure \ref{fig:man_fit}. Based on the obtained results we see that the chosen constants (as per Table \ref{tab:man_fit_par}) is only able to qualitatively mimic the trend of the measurement data. In particular, the inflection point of pH closely matches the inflection point of [$\caion$] (i.e. Q). This is to say that increase in alkalinity of the solution corresponds to decrease in $\caion$ concentration, while increase in acidity (or decrease in alkalinity) results in increase in $\caion$ concentration. The time point of these events match closely which can be attributed to the fact that $U_r \propto U_H$. Lastly, from the plot for experiment 4, we infer that constant pH profile results in linear decrease in $\caion$ ions which is suitably matched by the chosen constants. 
Despite the qualitative agreement, the quantitative deviation is too prominent to be ignored. Thus, it this method is not convincing enough to claim that the model is suitable to mimic the observed data. This motivates us to adopt a more mathematical approach for determining the unknown model function $U_r$.

\subsection{Forward backward adjoint SDE method}
With the aim of algorithmically determining the function $U_r$ we use two different methods. In the first method, we use a classical numerical scheme for solving the stochastic optimal control problem. For this purpose, as already mentioned above, the function $U_H$ is considered as a given or a known input function while the function $U_r$ is taken as the design variable. Now, applying the necessary and sufficient condition for the existence of an optimal control in Section \ref{sec:ocf}, we formulate the following optimality system \eqref{eq:opt_sys} with respect to the design variable $U_r$:
\begin{align}
\label{eq:opt_sys}
\begin{aligned}
    {\partial_{U_r} \mathcal{H}} &= 0 \\
    \partial_t X_t &= \partial_{\blam} \Ham(X, U_r, \blam) \\
    \partial_t \blam(t) &= -\partial_{\x} \Ham(X, U_r, \blam) \\
\end{aligned}
\end{align}
Since $\Ham$ is differentiable in $U_r$, the variation of \(\mathcal{H}\) with respect to \(U_r\) gives the optimality condition:
\begin{align}
\label{eq:opt_cond}
    -\alpha U_r + \Big(\lambda^k_C - \lambda^k_{Q} \Big) K_{sp} P Q \in {\partial_u \mathcal{H}}  
\end{align}
We can use \eqref{eq:opt_cond} to generate a sequence $(U^k_r)_{k\in\N}$ that maximizes the Hamiltonian in the following manner: 
\begin{align*}
    U^{k+1}_r &= (1 - \eta \alpha)U^{k} + \eta \Big(\lambda^k_C - \lambda^k_{Q} \Big) K_{sp} P^k Q^k\\
    \partial_t X^{k+1}(t) &= \partial_{\blam} \Ham(X^{k+1}, U^{k+1}_r, \blam^k) \\
    \partial_t \blam^{k+1}(t) &= -\partial_{\x} \Ham(X^{k+1}, U^{k+1}_r, \blam^{k+1}) \\
\end{align*}
Thus, for a sequence of $U^k_r$ such that $U^k_r \to U^*_r$ as $k \to \infty$ we see that the optimal design function $U^*_r$ is given as
\begin{align}
    U^*_r = {1 \over \alpha} \Big(\lambda^*_C - \lambda^*_{Q} \Big) K_{sp} P^* Q^*
\end{align}
Solving the above system involves solving the forward and backward adjoint system sequentially for which the above-mentioned discretization scheme \eqref{eq:disc_scheme} is used.
Furthermore, we make two vital simplifications that greatly improves the rate of convergence. Firstly, since $\caion$ is the only measurement that is used for determining the optimal design variable $U_r$, the role of adjoint variables $\lam_C, \lam_H, \lam_F$ is superfluous and can be set to 0 and omitted from the adjoint system.
Secondly, since [$\caion$] and [CO$_3^{2-}$] are small in magnitude and could potentially get close to zero, we find it convenient to choose the stepping parameter $\eta$ and $\alpha$ such that the products $\tilde \eta = \eta \alpha$ and $\hat \eta = \eta P Q $ are constants. Accordingly, the sequence of approximate functions $U^k_r$ is obtained using the following simpler relation
$$U^{k+1}_r = (1 - \tilde \eta)U^{k}_r - \hat \eta \lambda^k_{Q} $$
Overall the above described iteration scheme for obtaining an optimal design function $U^*_r$ is outlined in Algorithm \ref{algo:pgd} which we call as Forward Backward Stochastic Sweeping Method (FBSSM).

\SetKwFor{KwFor}{for}{do}{end}
\begin{algorithm}[!htb]
\caption{FBSSM \label{algo:pgd}}
 \KwData{$\epsilon > 0$, $\tau, h > 0$, $N_{\tau} \in \N$, $T > 0$}
 $U^0 := 0$ \\
 \KwFor{$k = 1, \dots$}
 {
    $X^0 = X_0$\\
    \KwFor{$j = 1, \dots, N$}
    {
        $X^{j+1}_i   = X^{j}_i + \tau \mu^{j}_i + \sqrt{\tau} ~ \Sigma^j_i ~ {\rm z}^j,$  \\[1ex]
    }
    
    $\blam^N = \blam_T$ \\
    
    \KwFor{$j = N-1, \dots, 1$}
    {
        $\bvsig^{j}   = {1 \over \tau} \E[{\rm z}^j \blam^{j+1} \Big| \F_{t_j}]$ \\
        $\blam^{j}_i  = \blam^{j+1}_i + \tau \Gamma^{j+1}_i - \sqrt{\tau} ~ \bvsig^{j} ~ {\rm z}^j,$
    }
    $U_r^{k+1} = (1 - \tilde \eta)U^k_r - \hat \eta \blam_{Q}$ \\[.5ex]
     
 }
\end{algorithm}

\begin{table}[h!]
    \centering
    \begin{tabular}{lll}
        \toprule
        Parameter & Description & Value \\
        \midrule
        $\tau$ & Time step size & 0.01 \\
        $h$ & Spatial step size & 0.1 \\
        $N_x$ & Number of spatial grid points & 64 \\
        $\alpha$ & Regularization parameter for control & 0.1 \\
        $\sigma_C$ & Noise standard deviation for C & 0.025 \\
        $\sigma_{Q}$ & Noise standard deviation for Q & 0.0005 \\
        $\sigma_H$ & Noise standard deviation for H & 0.001 \\
        $\tilde \eta$ & Step size for $U_r$ update  term & 0.0001 \\
        $\hat \eta$ & Step size for $\lam_Q$ update term & 0.005 \\
        \bottomrule
    \end{tabular}
    \caption{\label{tab:num_par_fbssm} Numerical parameters used in the FBSSM algorithm. \ref{algo:pgd}}
\end{table}

\begin{figure}[!h]
\begin{center}
    \includegraphics[scale=.195]{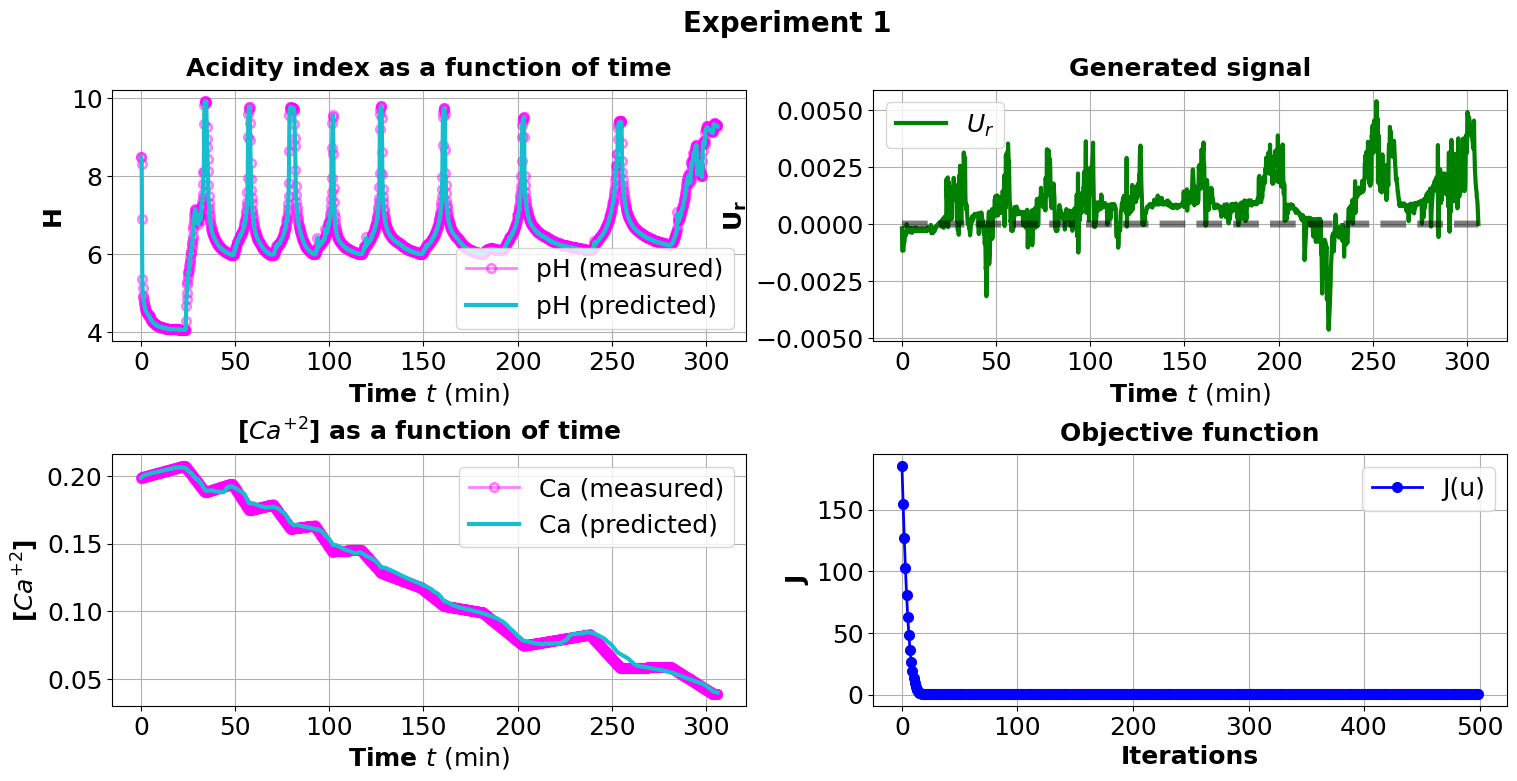} 
    \includegraphics[scale=.195]{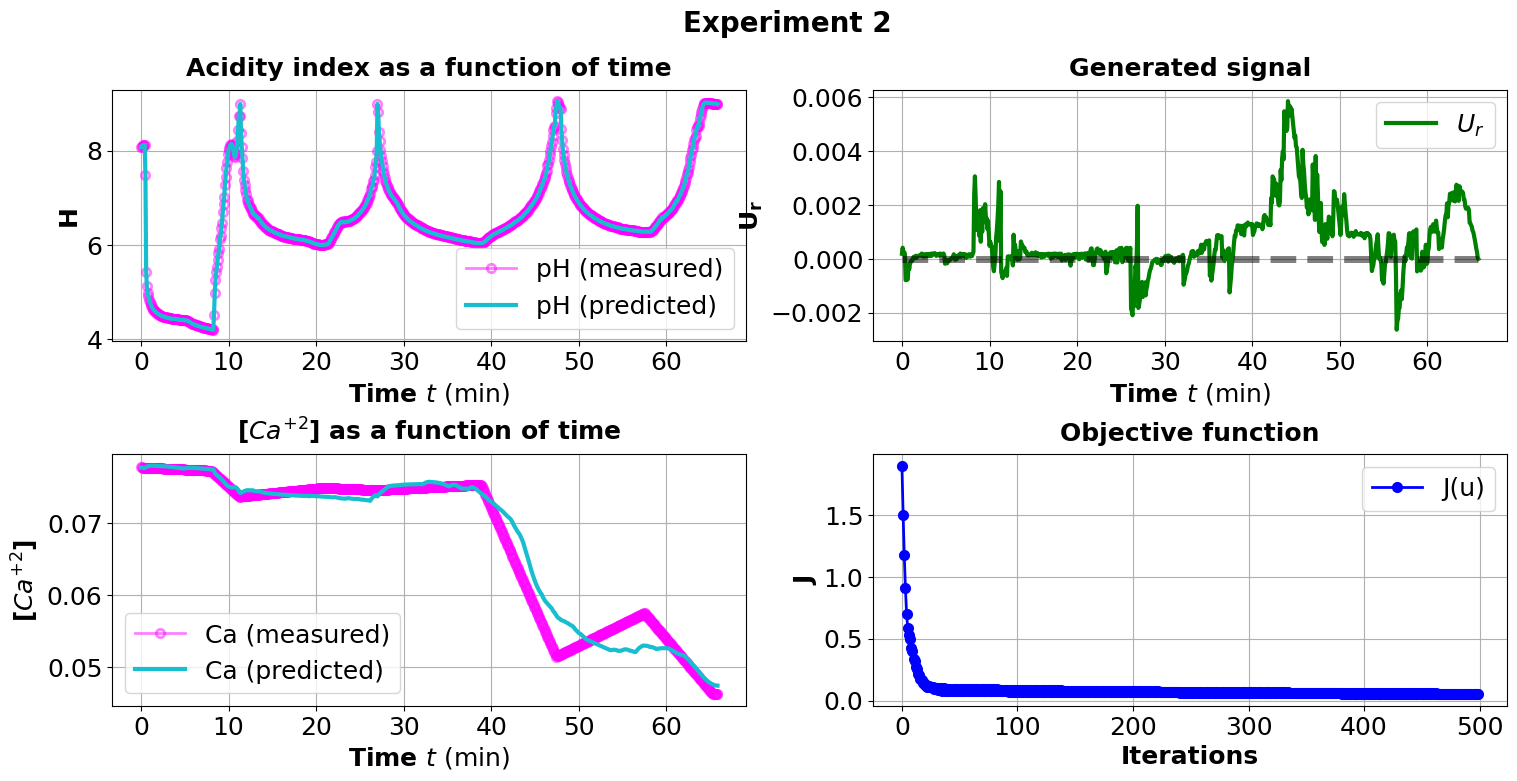} \\
    \includegraphics[scale=.195]{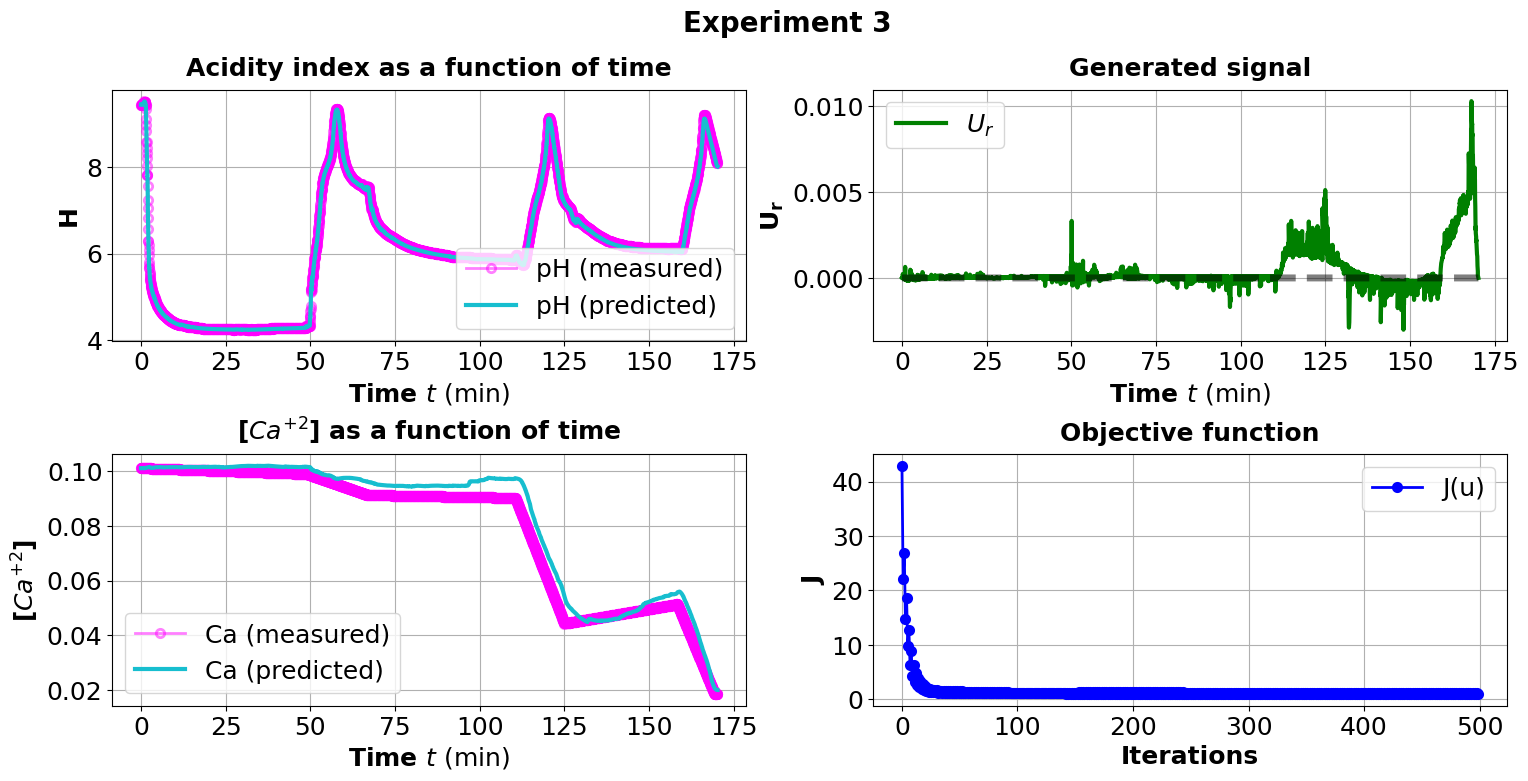}
    \includegraphics[scale=.195]{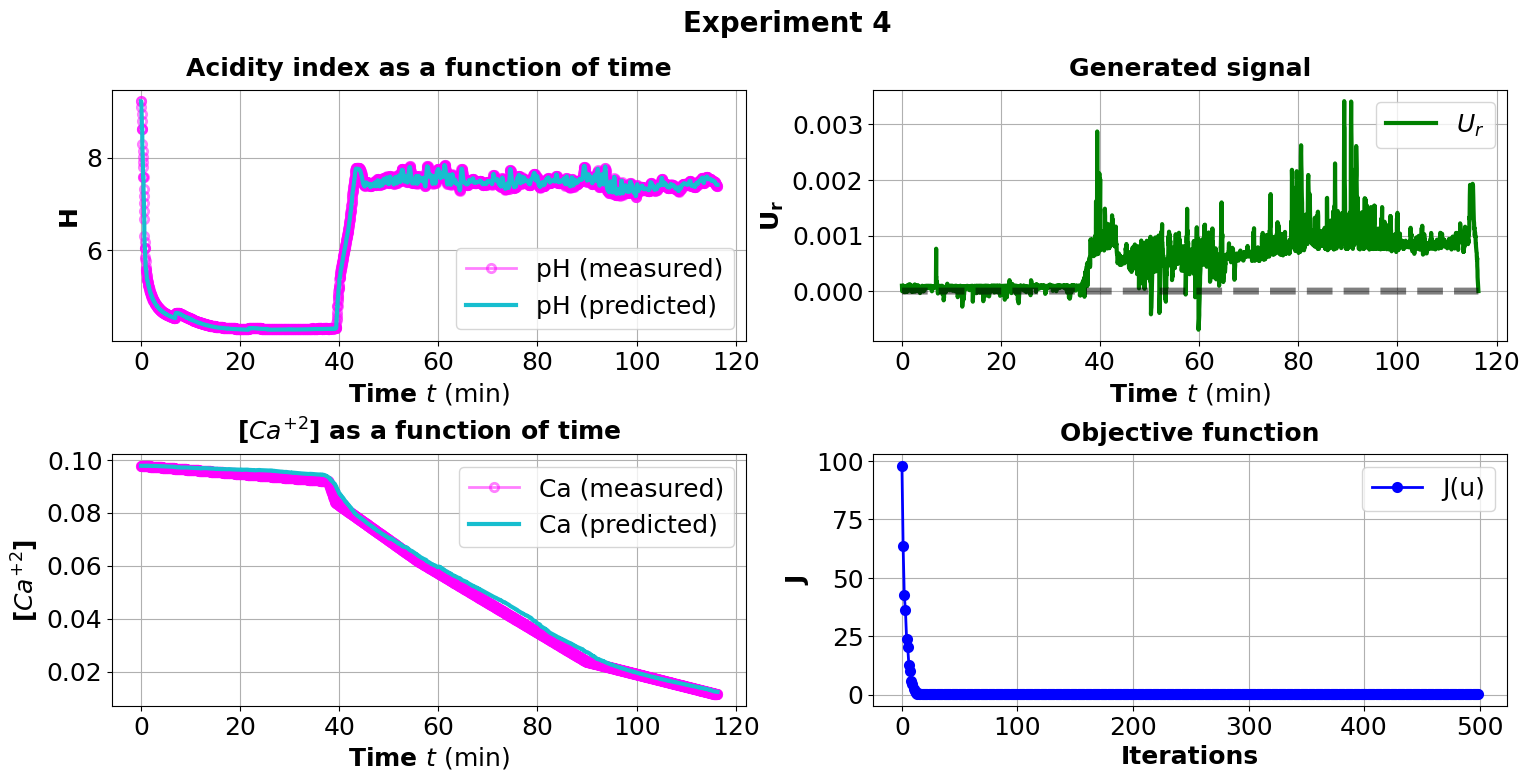} \vspace*{-.1cm}
\end{center}
\caption{\label{fig:FBSSM} Optimal design function $U_r$ obtained after running the Algorithm \ref{algo:pgd} for 500 iterations for each of the four experimental data trajectory separately. For each experimental plot, we see 4 subplots which are arranged from top-left to bottom-right in the following manner: the input function pH, the approximated or generated signal $U_r$, observable function Q ([$\caion$]) and the objective function $J(U_r)$. In each of the plots the raw experimental data trace is shown in magenta color.}
\end{figure}

The constants used for the implementation of the FBSSM are provided in Table \ref{tab:num_par_fbssm}. The algorithm was run for 500 iterations and the obtained results are as shown in Figure \ref{fig:FBSSM}. The figures indicate that the proposed method is able to monotonically reduce the cost $J$ and the approximated solution $Q^k$ corresponding to $U^k_r$ slowly converges to the target profile $\bar Q$. This is evident in the plots for Experiments 1, 3 and 4. In there we can observe that the approximated $\caion$, $\hat Q$ (shown in cyan) closely matches the experimental trace $\bar Q$. However, for Experiment 2, the approximated $\caion$ deviates quite a bit from the experimental data. This can be improved by simply running (or continuing) the algorithm for more number of iterations. Unfortunately, this is indeed a major disadvantage of this method where the order of convergence is $\mathcal{O}({1/2})$ \cite{Zeller1956, Chassagneux2017} and thus the time required for convergence can be quite large, up to 3-4 hours in particular for a long-time experiment trajectory.
In order to address this, we adopt a more modern approach of using DNN to solve the stochastic optimization problem \eqref{eq:socp}.

\subsection{Neural network based method}
In order to address the shortcomings of the FBSSM algorithm, in this section we investigate the neural network based method to estimate the function $U_r$ that fits the data. For this purpose we investigated two distinct neural network architectures: an Artificial Neural Network (ANN) and a Gated Recurrent Unit (GRUN) based network (GRUN). 

\subsubsection{Artificial Neural Network (ANN)}

The ANN comprised a series of linear layers and activation functions designed to process input data without inherent temporal dependencies. The network consisted of three linear layers (Linear-1, Linear-3, Linear-4) and two activation functions (ReLU-2, Tanh-5). The input, with a shape of $(32, 128)$, was initially transformed to a shape of $(32, 128)$ via Linear-1, which employed 16,512 trainable parameters. This output was then passed through a Rectified Linear Unit (ReLU-2) activation function. A subsequent linear layer (Linear-3), utilizing 16,512 parameters, further processed the data to a shape of [32, 128]. Finally, Linear-4, also with 16,512 parameters, mapped the output to [32, 128], followed by a hyperbolic tangent activation function (Tanh-5). The ANN, with a total of 49,536 trainable parameters, processed the input to generate the final predictions.

\subsubsection{Gated Recurrent Unit Network (GRUN)}
The GRUN was designed to capture temporal dependencies within sequential data. It is composed of a three-layered bidirectional GRUN (GRU-1) with 32 hidden units per direction. The input sequence, with a shape of $(32, 64, 2)$, was fed into the GRUN, along with a context vector of shape $(6, 32, 32)$. The GRUN output, $(32, 64, 64)$, was then passed through a linear layer (Linear-2) with 130 parameters, transforming the output to $(32, 64, 2)$. Subsequently, a hyperbolic tangent activation function (Tanh-3) was applied, maintaining the output shape. The GRUN, with a total of 44674 trainable parameters, was responsible for extracting relevant temporal features from the input sequence and generating predictions.
Both networks are designed to take as input the experimental data $\bar U_H$ and $\bar Q$, namely the system's input and observation signal respectively, and generates as output the input signal $\hat U_H$ and the design function $\hat U_r$. The latter is obtained via an explicit form $\hat U_r = \hat U_H \tilde U_r$, where $\tilde U_r$ is one of the direct output (the other being $\hat U_H$) of the network. We note here that, the explicit form of $\hat U_r$ is not a necessity but is simply motivated from \eqref{eq:man_ur} and due to the convenience network based optimization method provides. This choice also allows use to indirectly validate our intuition of considering the explicit dependence of $U_r$ on $U_H$. Alternatively, one can set $\hat U_r$ to be equal to $\tilde U_r$, i.e. as one of the direct output of the network itself.

\begin{figure}[!h]
\begin{center}
    \includegraphics[scale=.195]{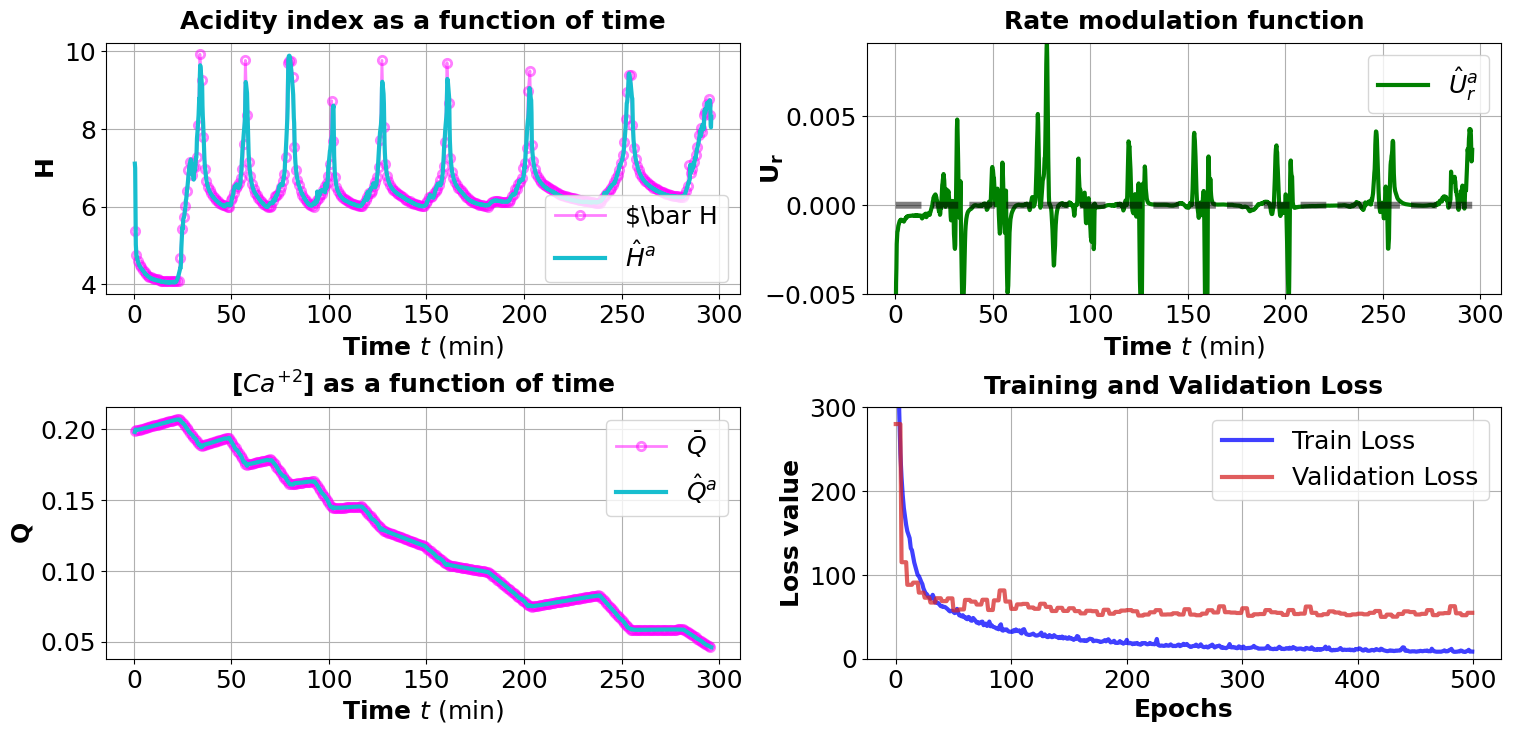} 
    \includegraphics[scale=.195]{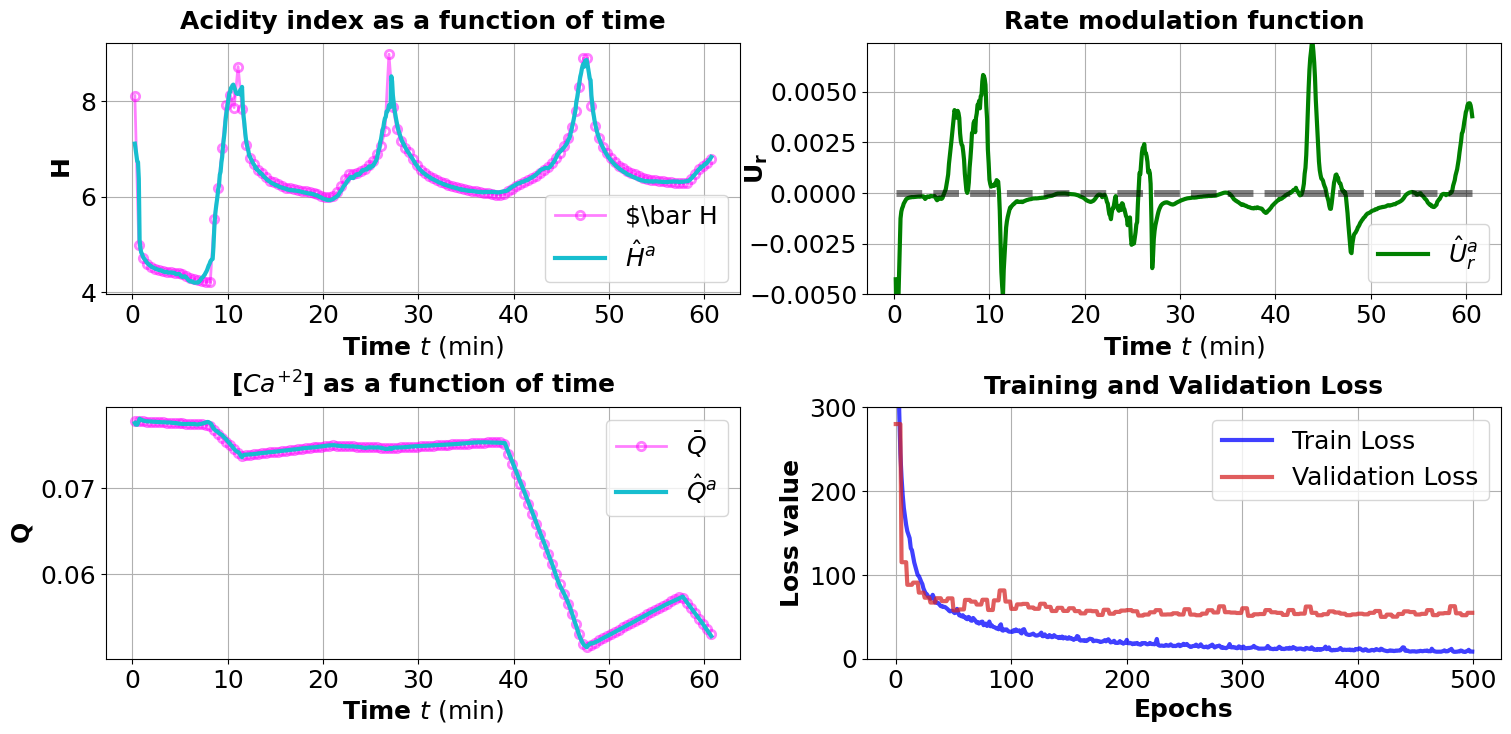} \\
    \includegraphics[scale=.195]{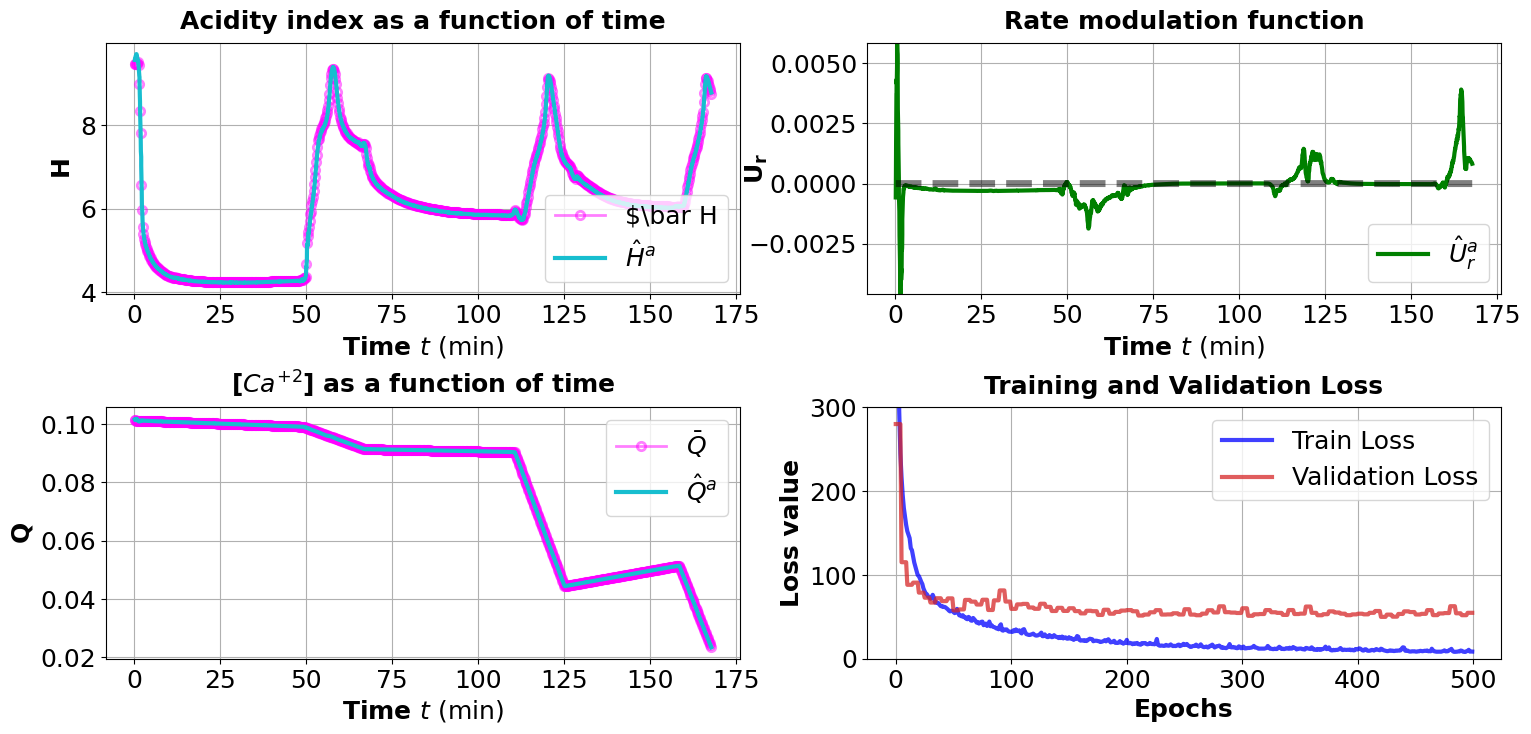} 
    \includegraphics[scale=.195]{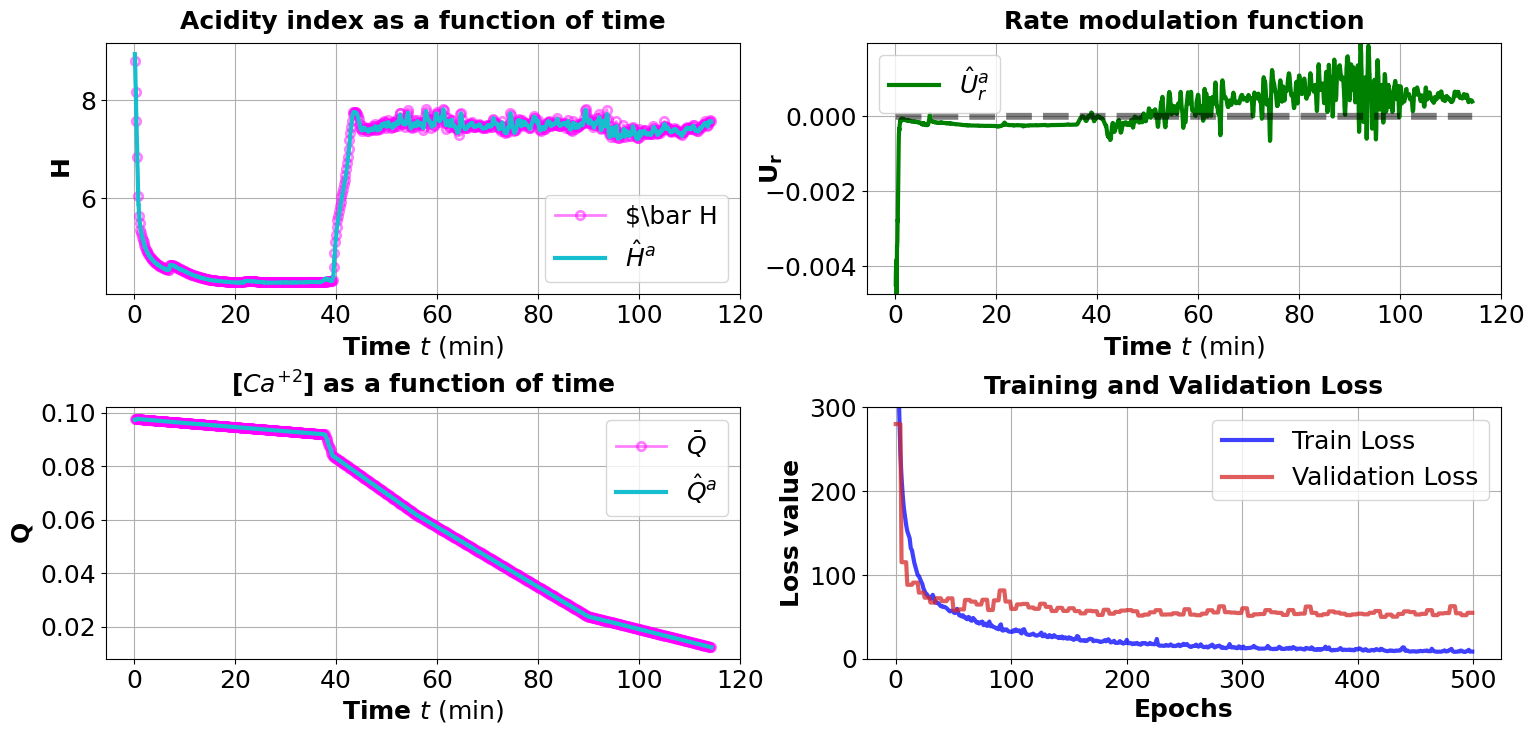}
\end{center}
\caption{\label{fig:ann_fit} Prediction results of the ANN model on the entire dataset.}
\end{figure}

\begin{figure}[!h]
\begin{center}
    \includegraphics[scale=.195]{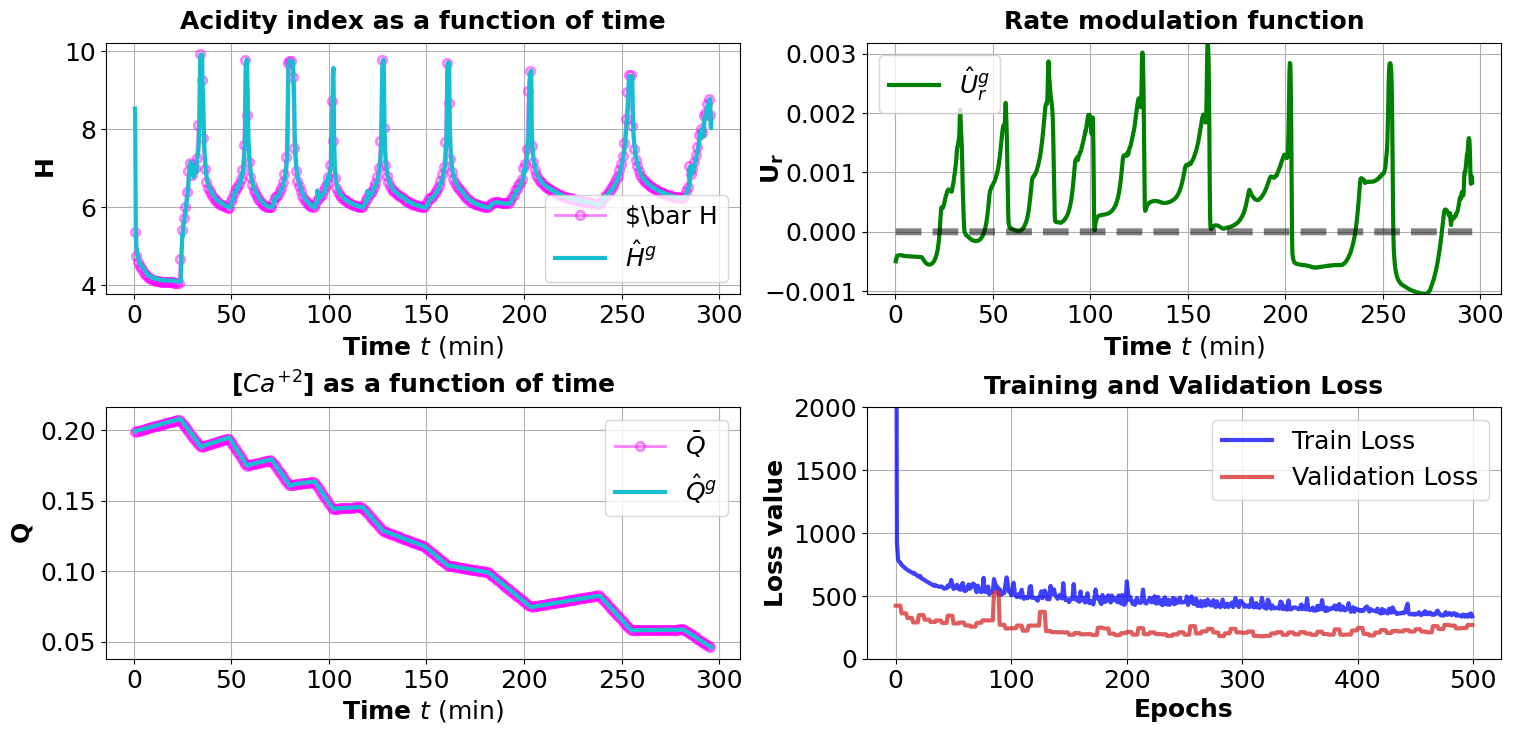} 
    \includegraphics[scale=.195]{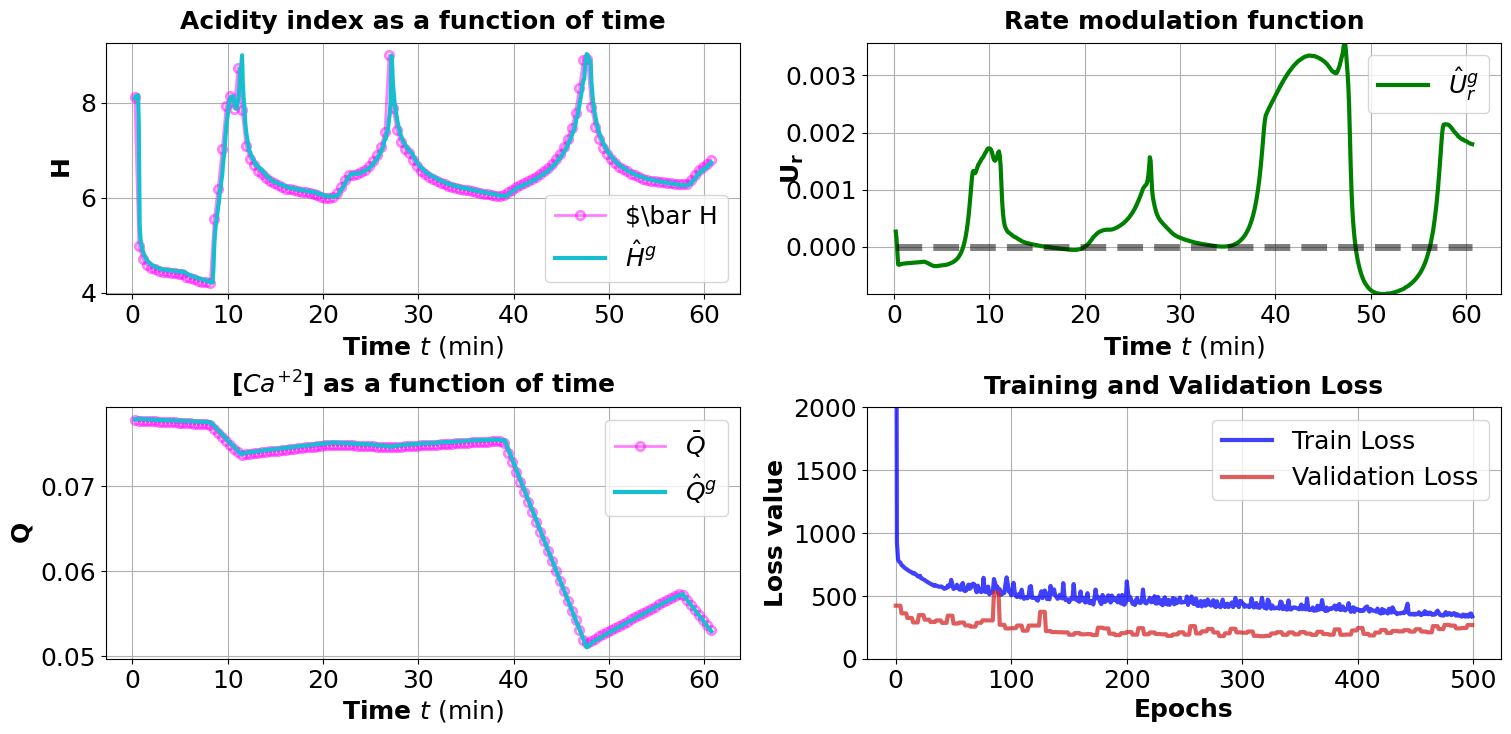} \\
    \includegraphics[scale=.195]{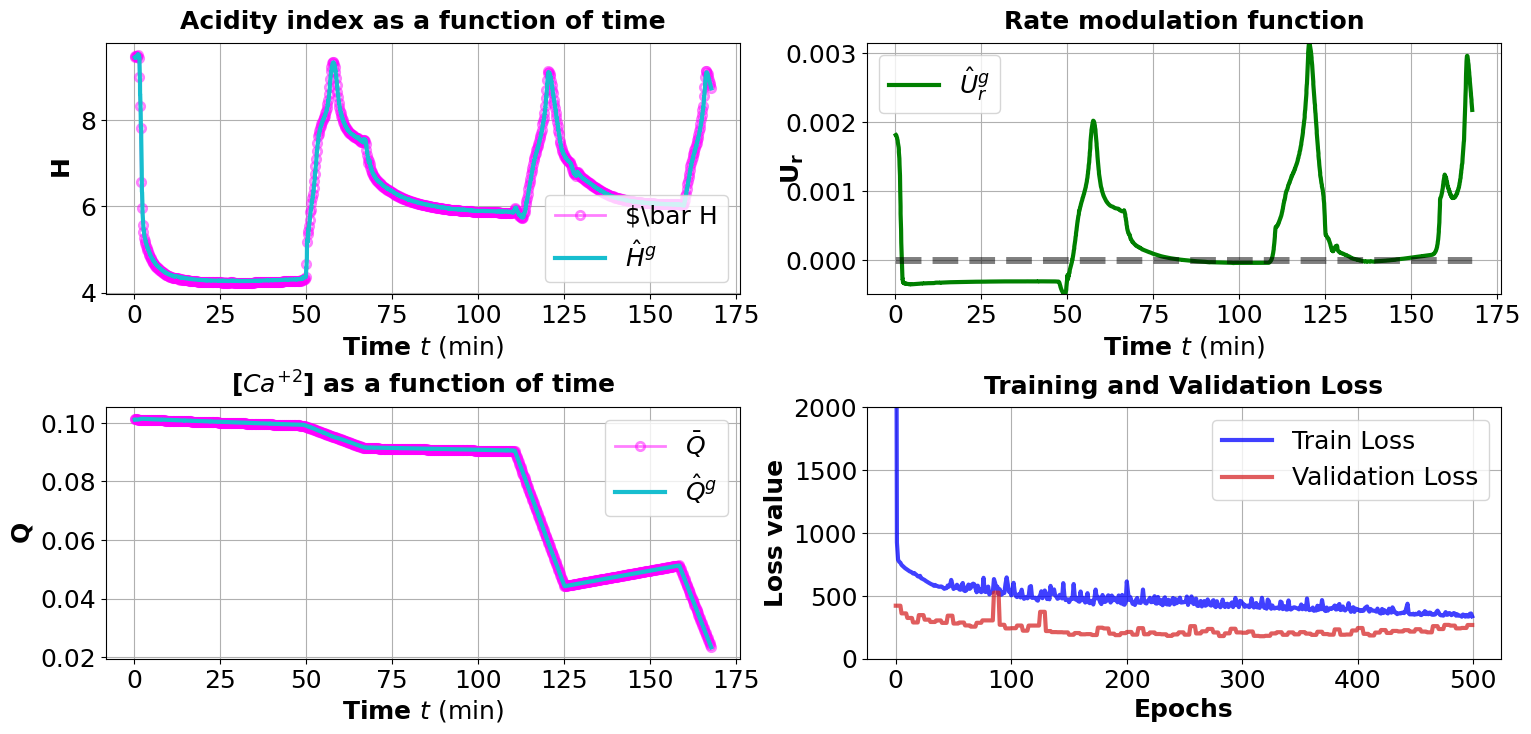} 
    \includegraphics[scale=.195]{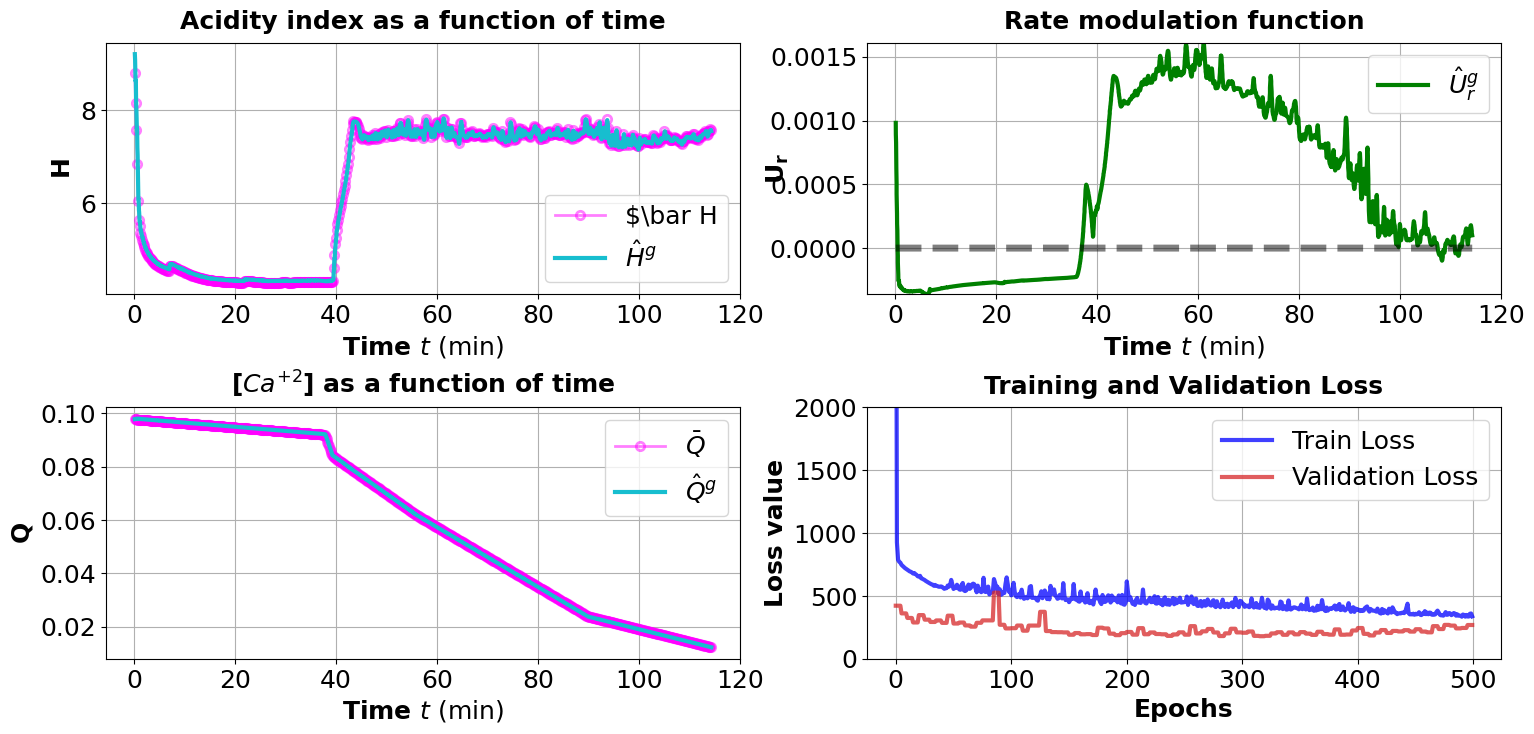}
\end{center}
\caption{\label{fig:gru_fit} Prediction results of the GRUN model on the entire dataset.}
\end{figure}

Both the network architectures were implemented in PyTorch framework. The models were trained using Adamax optimization method with learning rate $0.001$ and other hyperparameters. The loss function used for training is
$$\mathfrak{L} = \ell_1\|\hat U_r - \bar U_r\|^2_{L^2([0,T])} + \ell_2\|\hat U_H - \bar U_H\|^2_{L^2([0,T])} + \ell_3\|\hat H - \bar H\|^2_{L^2([0,T])} + \ell_4\|\hat Q - \bar Q\|^2_{L^2([0,T])},$$ 
where $\bar U_r$ is the manually chosen function (as per \eqref{eq:man_ur}), $\bar H$ and $\bar Q$ are the pH and $\caion$ concentration obtained experimentally, $\bar U_H$ is derived from data via simple finite difference  as $\bar U_H = \bar H_{t_{n+1}}- \bar H_{t_n}$
Based on this and other hyperparameters, listed in \ref{tab:dnn_hparams}, the models were trained for 500 epochs with batch size 128 and the obtained results are as show in the Figures \ref{fig:ann_fit} and \ref{fig:gru_fit}.
\begin{table}[h!]
    \centering
    \begin{tabular}{lll}
        \toprule
        Parameter & Description & Default Value \\
        \midrule
        \multicolumn{3}{l}{\textbf{Training Parameters}} \\
        \midrule
        $L$ & Sequence length for training/testing & 64 \\
        $m$ & Batch size for training & 128 \\
        $N_E$ & Number of epochs for training & 500 \\
        $\eta$ & Learning rate for optimizer & 0.001  \\
        \midrule
        \multicolumn{3}{l}{\textbf{Model Configurations}} \\
        \midrule
        $S_I$ & Input features & 2 \\
        $S_O$ & Output features & 2 \\
        $I_N$ & Input size (ANN) & $S_I \times L$ = 128 \\
        $S_{H1}$ & Hidden layer feature size (ANN) & 128 \\
        $S_{H2}$ & Hidden layer feature size (GRUN) & 32 \\
        $O_N$ & Output size (ANN) & $S_O \times L$ = 128 \\
        $N_L$ & Number of layers (GRUN) & 3 \\
        $Bi$ & Bidirectionality (GRUN) & True \\
        \midrule
        \multicolumn{3}{l}{\textbf{Loss Weights}} \\
        \midrule
        $\ell_1$ & Weight for $\hat U_r$ loss term & 1 \\
        $\ell_2$ & Weight for $\hat U_H$ loss term & 1 \\
        $\ell_3$ & Weight for $\hat H$ loss term & 10 \\
        $\ell_4$ & Weight for $\hat Q$ loss term & $10^6$ \\
        \bottomrule
    \end{tabular}
    \caption{Training and Model Configuration Parameters \label{tab:dnn_hparams}}
\end{table}

\section{Discussion \label{sec:discuss}}
In this section we shall compare and discuss the performance of different methods used to generate the design function $U_r$. As a consequence we also infer the validity of the mathematical model \ref{eq:abs_psd_spde} as a surrogate for the pH driven precipitation process. Firstly, as already mentioned above, the design function $U^m_r$ obtained via manual method \eqref{eq:man_ur} provides only a broad alignment of the qualitative behavior of $\hat{Q}$ with respect to $\bar Q$. The main properties that $U^m_r$ is able to reconstruct are: (i) the decreasing trend of $Q$ during increasing pH and increasing trend of $Q$ during decreasing pH agrees for both $\bar Q$ and $\hat Q$. (ii) The time points when $Q$ changes direction from an increasing function to a decreasing function and vice-versa matches. (iii) Thirdly, the magnitude of increase during increasing pH is relatively smaller compared to the magnitude of decrease during increase in pH. Despite these qualitative properties, the lack of quantitative match, i.e. the prevalent large error $\|\hat Q -\bar Q\|$, is not insignificant to be ignored. This led us to pursue a more mathematical approach for which two approaches are used: (i) classical method of forward-backward SDE sweep method and (ii) a statistical DNN based approaches are used to approximate the optimal design function $U^{*}_r$.

The classical FBSDE method is an online optimization method the computes the optimal design function $\hat U^f_r$ by constructing a sequence of approximates that successively reduces the cost value $J$. The sequence construction is performed for each observation data as well as model parameters kept fixed. Thus, the approximated design $\hat U^f_r$ is a function of the data trace $\bar Q$ and model parameters $\theta$. Thus, if one of them (either $\bar Q$ or $\theta$) changes the algorithm has to be rerun in order to construct a new sequence that converges to a new design function that explains the given data $\bar Q$. The proposed FBSSM algorithm (Algorithm \ref{algo:pgd}) has a sublinear order of convergence and requires up to 500 iterations, which amounts to around 4-5 hours of user runtime, to obtain significant reduction in the error. From Figure \ref{fig:FBSSM} we see that after 500 iterations the quantitative fit for the Experiments 1,3 and 4 are 90\% complete while for Experiment 2 the fit is only 70\% complete. It would require another 500 iterations to reach 90\% completion. This aspect makes the FBSSM quite cumbersome and unsuitable for real-time system identification or state estimation and thus also for autonomous control. 

\begin{figure}[!htbp]
\begin{center}
    \includegraphics[scale=.4]{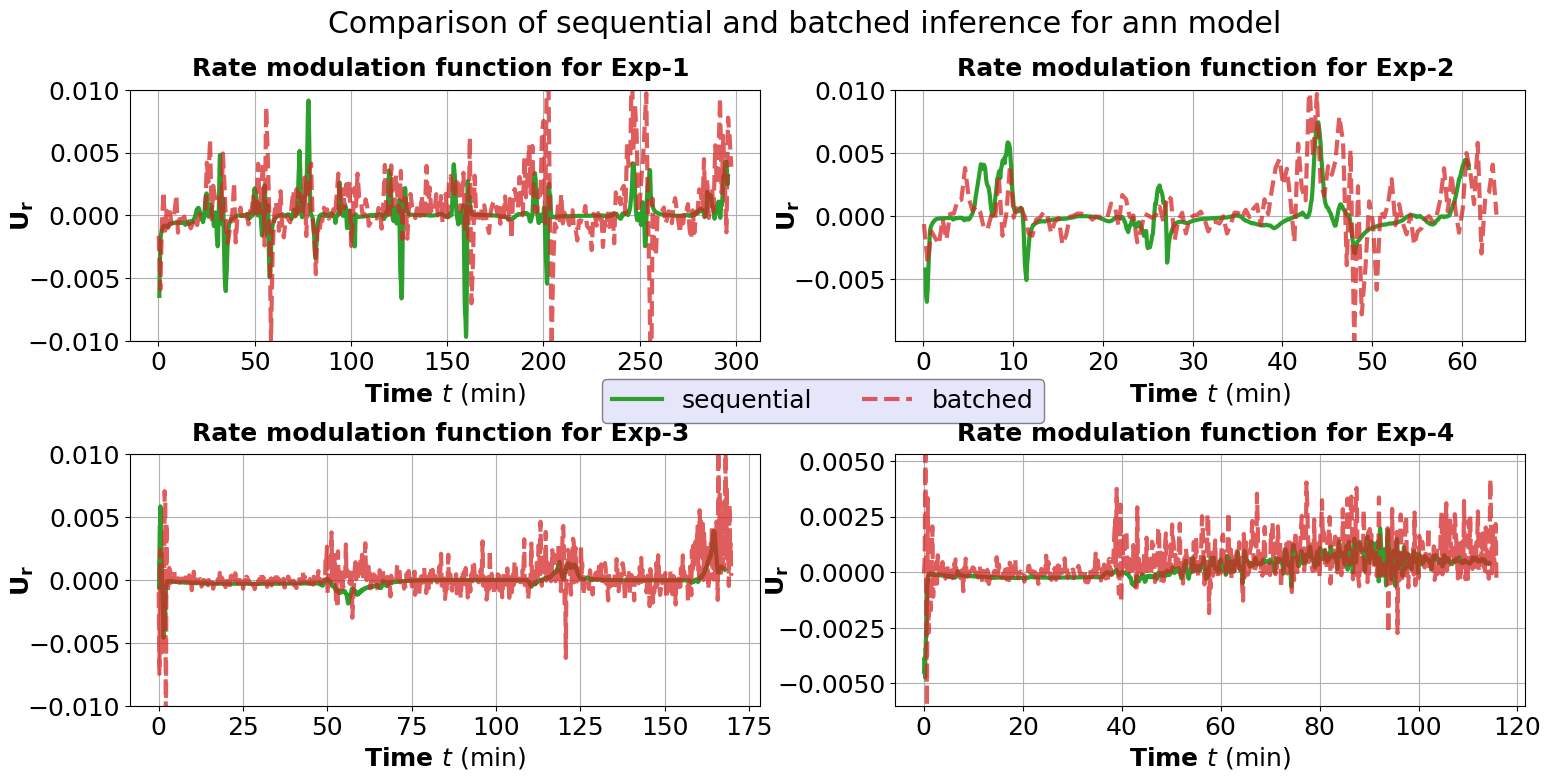} \\
    \includegraphics[scale=.4]{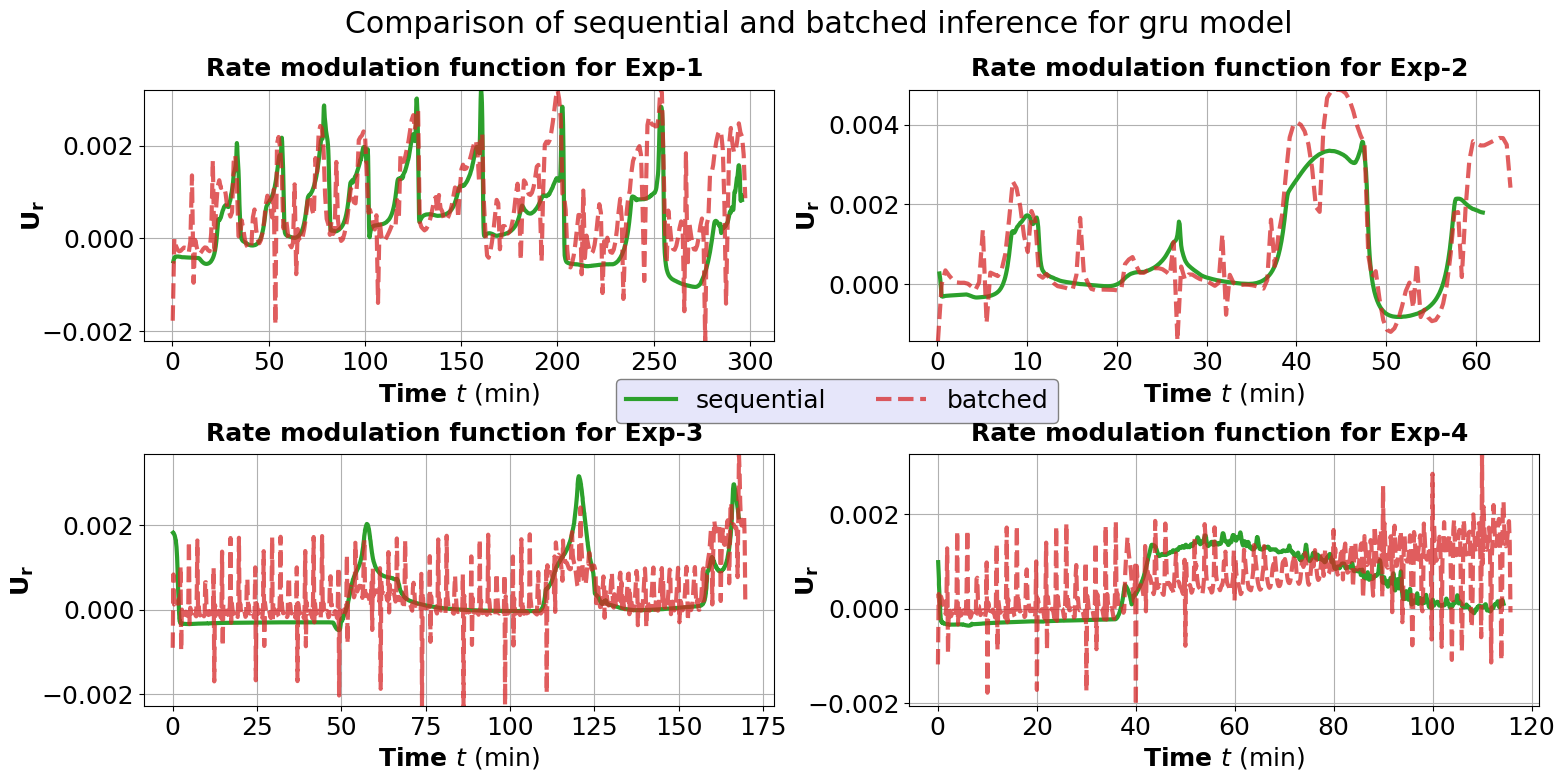}
\caption{\label{fig:ann_vs_gru_fit} Comparison of the design function generated $U_r$ generated ANN and GRUN.}
\end{center}
\end{figure}

To overcome this restriction, we adopted the machine learning based method which considers all the experimental data and aims to approximate the stochastic optimization solver for any experimental observation that is supplied. Accordingly, we trained two different architectures that take the observations as input, while keeping the model parameters $\theta$ fixed,  and predict the input $U_H$ and the design $U_r$ functions simultaneously. Although, one could also generalize the network design to consider all different model parameters, we have kept the scope of the problem limited to a fixed parameter case. The predictions obtained by both the models are as shown in Figures \ref{fig:ann_fit} and \ref{fig:gru_fit} respectively. In both cases we see that the predicted pH and $\caion$ profiles fit the observations perfectly, i.e. they are able to learn from the dataset with training loss reducing significantly and converging to a low value while the validation loss following a similar pattern and also converging near training loss. The convergence of the ANN model is much faster and smoother compared to the GRUN model. This can be attributed to the fact that for GRUN model has a recurrent architecture with continuous activation functions, thus making it sensitive to non-smooth input sequence mainly the input data $U_H$. Consequently, it also results in longer training time in comparison to the ANN model. 
Thus, due to different model architectures, the generated design functions $\hat U_r^a$ and $\hat U_r^g$ from ANN and GRUN based models have slightly different profiles. To compare the two we test them on different types of inference sequences. We consider rolling window sequence and disjoint batched sequence. The former represents a continuously varying data while the latter represents accumulated observation sequence obtained in regular intervals. Both the models have generated smoother functions for the rolling window sequence while for the batched sequences the predicted signal has relatively higher oscillations. In particular $\hat U_r^a$ is more oscillatory in comparison to the $\hat U_r^g$ (as seen in Figure \ref{fig:ann_vs_gru_fit}).

\begin{figure}[!htbp]
\begin{center}
    \includegraphics[scale=.4]{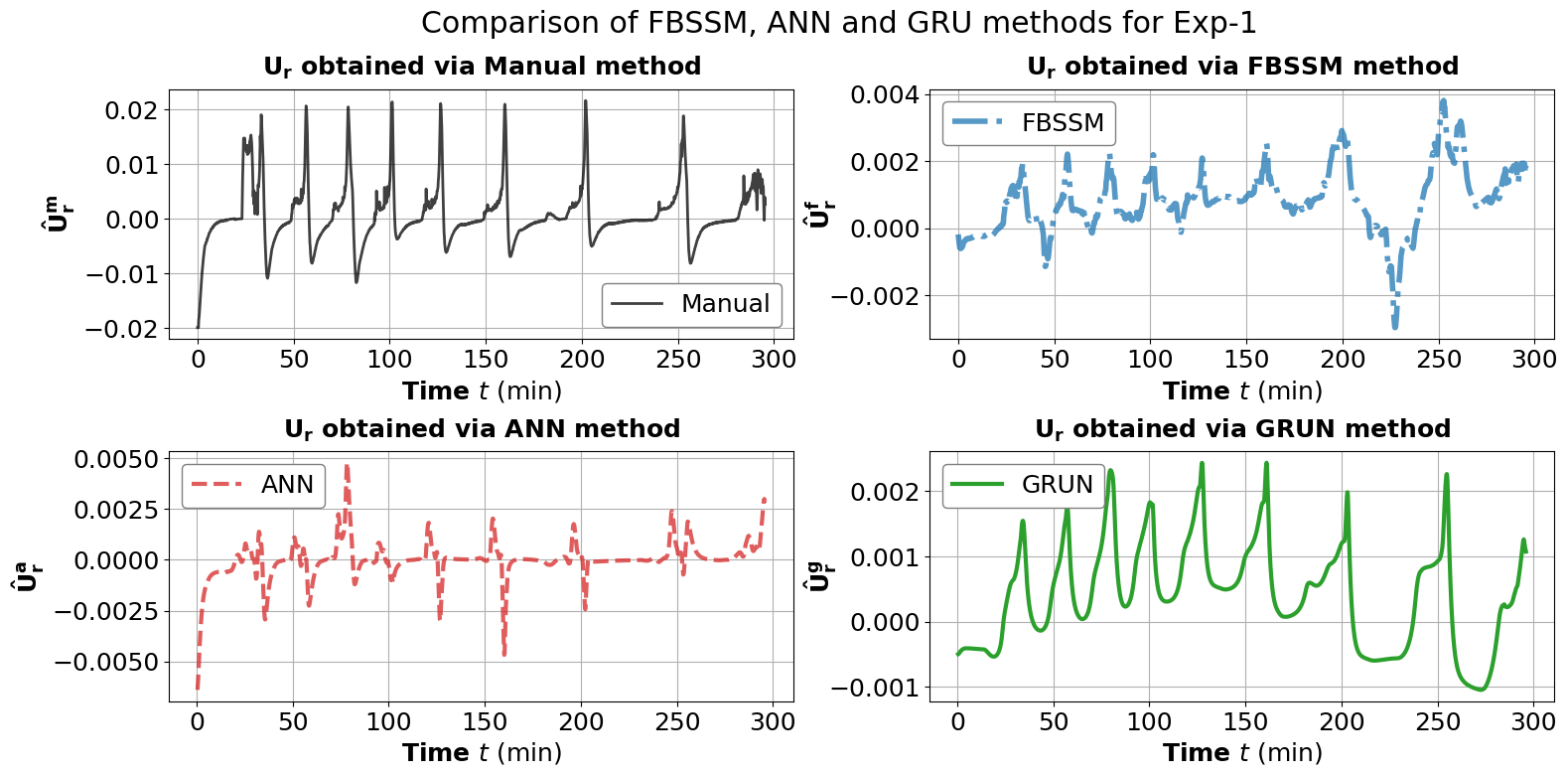}\\
    \includegraphics[scale=.4]{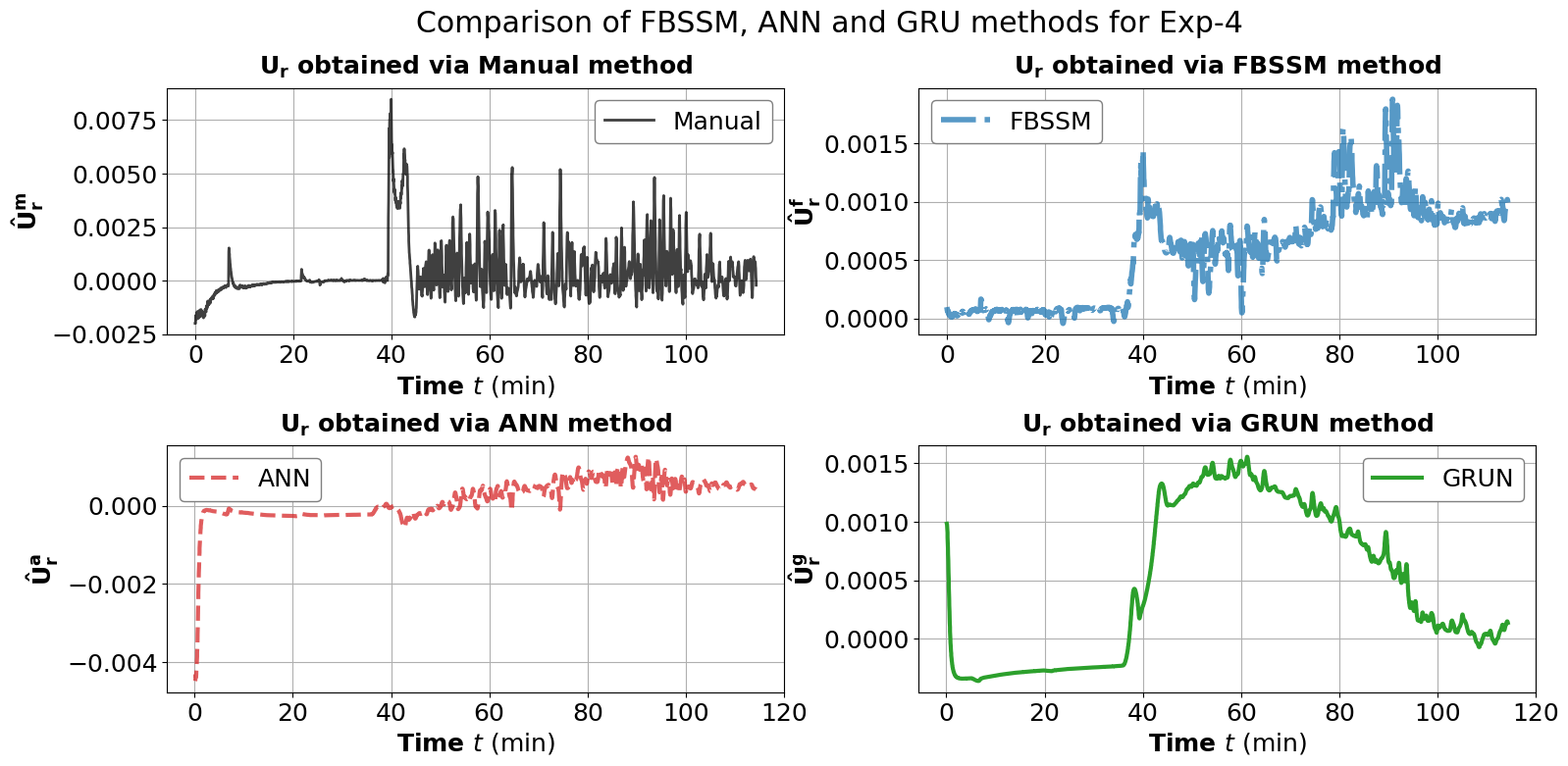}
\caption{\label{fig:ann_vs_gru_vs_fbsde} Comparison of predicted design function $\hat U_r$ by four different methods namely- manual, FBSSM, ANN and GRUN for experiment data 1 (top figure) and 4 (bottom figure).}
\end{center}
\end{figure}

Furthermore, comparing the predicted $U_r$ by all the four methods namely- manual ($\hat U^m_r$), FBSSM ($\hat U^f_r$), ANN ($\hat U^a_r$) and GRUN ($\hat U^g_r$) we observe the following. From Figure \ref{fig:ann_vs_gru_vs_fbsde} observing the plots for Experiments 1 to 3, we see that, $U^g_r$ is a smoother signal in comparison the others and has fewer fluctuations. The GRUN model generates a smoother signal compared to the others. The magnitude of the min and max values for $U^g_r$ is also lesser compared to $\hat U^m_r$ and $\hat U^a_r$. Apart from these difference, the one thing that $\hat U^g_r$, $\hat U^a_r$ and $\hat U^m_r$ have in common is the roughly matching time points where the signal fluctuates in order to change the direction of $Q$ from increasing to decreasing and vice-versa. This also validates the choice of the manual function $\hat U^m_r$ and indeed justifies it as a rough approximation. 
Lastly, the $\hat U^f_r$ signal has a different profile with respect to all others and is an indication that it is still far from the optimal solution and would require much more iterations.

\begin{figure}[!h]
\begin{center}
    \includegraphics[scale=.4]{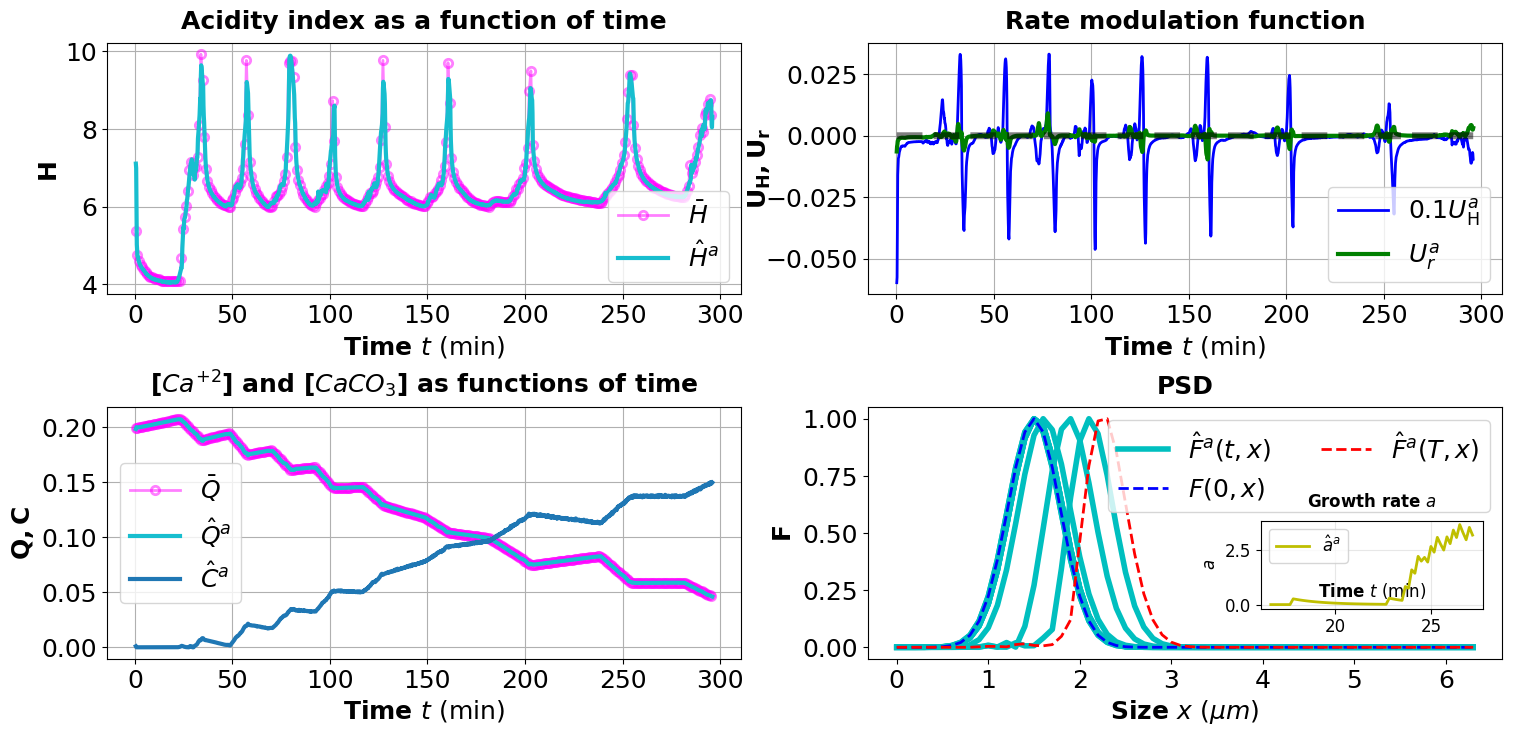}\\ 
    \includegraphics[scale=.4]{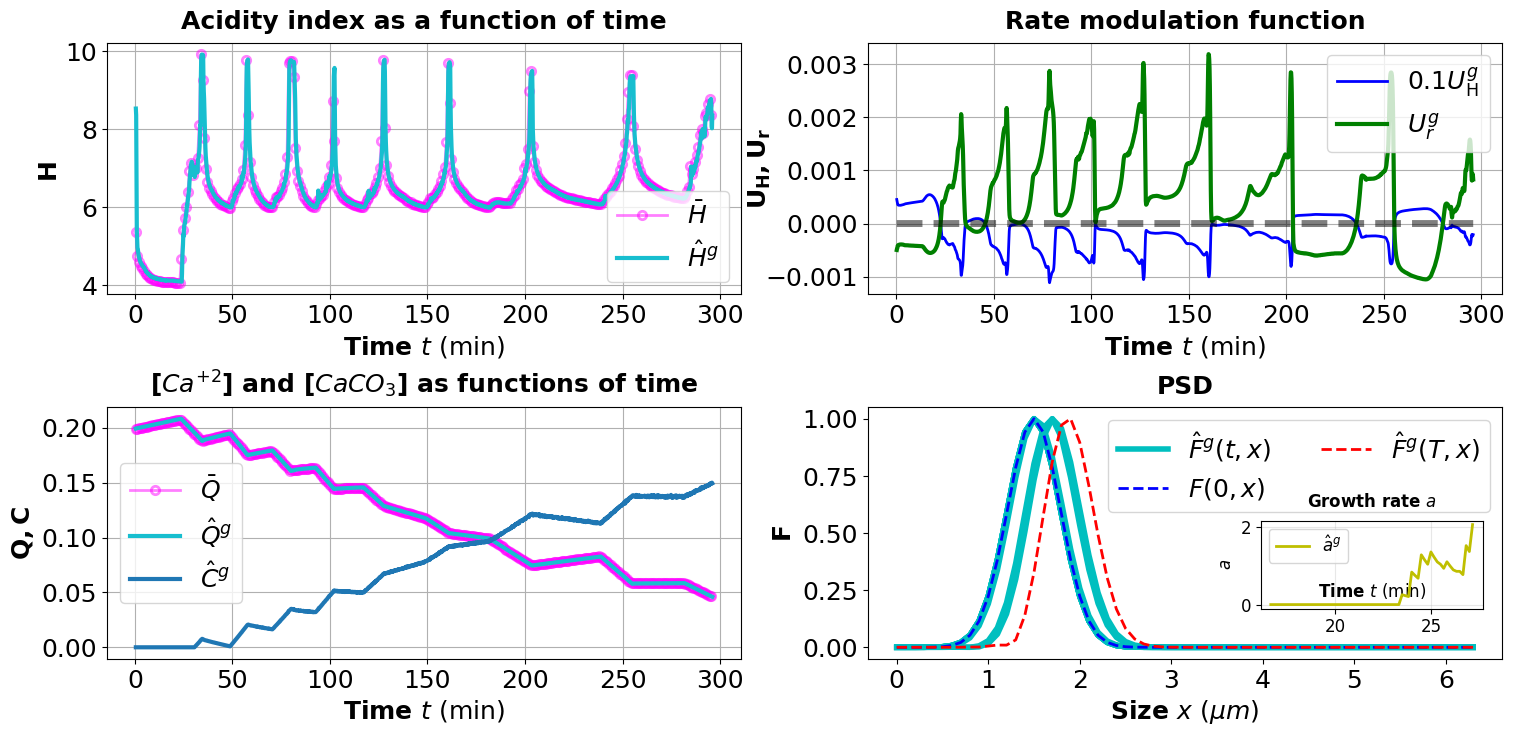}\\[-3ex]
\end{center}
\caption{\label{fig:ann_gru_full_sys_one} Plots of the state variables of the system resulting from $\hat U^g_r$ obtained by ANN (left) and GRUN (right) model, respectively, for Experiment 1.}
\end{figure}

Since, $\hat U^g_r$ and $\hat U^a_r$ have better predictions compared to $\hat U^f_r$ and $\hat U^m_r$ we use the ANN and GRUN model and plot all the state variables to recover the behavior of the entire system. This can be seen in Figure \ref{fig:ann_gru_full_sys_one} where we observe that for each data trace of Experiment 1, the decrease in $\caion$ is accompanied by the corresponding increase in $\caco$ in aqueous phase. As a consequence of which growth of precipitated solids takes place which is depicted by the shift in PSD function $\hat F$  towards the higher particle size. Note that interestingly, the ANN model produces faster growth compared to GRUN model. This can be observed by the large shift in predicted PSD $\hat F^a$ by ANN as compared to $\hat F^g$. This is also made evident by the growth rate inset-subplot, where in the ANN model produces faster increase in growth rate compared to the GRUN model. The underlying reason for this being that, the higher magnitude in $U_H$ produces slight differences in the $pH$ profile which in turn results in the difference in $\csat(H)$ thus resulting in the difference in the growth rates. Although in the current scope of the study, this does not have an impact in the design of the model but would play a significant role for when the objective is to control the process to obtain a PSD of certain profile also subjected to some time constraints. 

\section{Conclusion \label{sec:conlcusion}} 
In this work we have introduced a mathematical model based on coupled degenerate stochastic partial differential equations to elucidate the dynamics of semi-batch calcium carbonate precipitation process due by pH changes induced by NaOH and/or CO$_2$. We provide detailed mathematical analysis for the well-posedness of the system based on which a sufficient condition for the stationarity is obtained in Theorem \ref{thm:exist_uniq_psd} and Theorem \ref{thm:stable_soln} respectively. Subsequently, with the aim of finding an optimal design or control function, we also formulated a stochastic optimal control problem for which necessary and sufficient conditions for the existence a solution was also established. This further enabled us to propose the FBSSM algorithm to construct $U_r$ in an automated manner. Since, the FBSSM is computationally expensive and repetitive, we further proposed machine learning based methods in the form of ANN and GRUN methods which were able to obtain perfect fit on the data. The explicit form of $\hat U_r = \tilde U_r \hat U_h$ in the network implementation confirms the intuition of the manually obtained design function $U^m_r$ which was only able to capture some high level qualitative properties. Altogether, using three different methods we were able to establish the validity of mathematical model as a surrogate for the precipitation process. Of the three methods, manual tuning provided an intuitive initial design but lacked precise control. Algorithmic fitting, however, demonstrated significantly better performance, with the FBSSM achieving reasonable fits after 200-500 iterations. Notably, a DNN-based method surpassed all others, delivering the most accurate fitting and offering substantial computational advantages for potential real-time process monitoring and control. Consequently, this work provides the basis for autonomous control of the precipitation process where the objective is then to design an optimal pH profile that can drive the process to obtain precipitates of desired size profile. 

\section{Appendix \label{sec:appendix}}
\subsection{Ornstein-Uhlenbeck Process}
\begin{lem}[Solution of Ornstein-Uhlenbeck Process]
Let $X_t$ satisfy the SDE
\begin{equation}
\label{eq:ou_sde}
dX_t = \theta (\mu - X_t) \dt + \sigma \dWt{}, \qquad X(0) = X_0,
\end{equation}
where $\theta > 0$, $\mu \in \R$, $\sigma > 0$. Then the solution is
\begin{equation}
\label{eq:ou_solution}
X_t = X_0 e^{-\theta t} + \mu (1 - e^{-\theta t}) + \sigma \sqrt{\frac{1 - e^{-2\theta t}}{2\theta}} Z,
\end{equation}
where $Z \sim \mathcal{N}(0,1)$.
\end{lem}
\begin{proof}
Multiply by the integrating factor $e^{\theta t}$ and integrate:
\[
e^{\theta t} X_t = X_0 + \int_0^t \theta \mu e^{\theta s} ds + \int_0^t \sigma e^{\theta s} dW_s.
\]
Solve for $X_t$ and compute the variance of the stochastic integral:
\begin{align*}
X_t &= X_0 e^{-\theta t} + \mu (1 - e^{-\theta t}) + \sigma e^{-\theta t} \int_0^t e^{\theta s} dW_s, \\
&= X_0 e^{-\theta t} + \mu (1 - e^{-\theta t}) + \sigma \sqrt{\frac{1 - e^{-2\theta t}}{2\theta}} Z, \quad Z \sim \mathcal{N}(0,1)
\end{align*}
since $\int_0^t e^{\theta s} dW_s \sim \mathcal{N}\left(0, \frac{e^{2\theta t} - 1}{2\theta}\right)$. Substitute and simplify.
\end{proof}

\subsection{Exponential OU Process}
\begin{lem}[Solution and Moments of Exponential OU Process]
Let $X_t$ satisfy the SDE
\begin{equation}
\label{eq:Exp_OU_proc}
    d X_t = \left(\mu(t) - \theta X_t\right) dt + \sigma X_t dW_t, \qquad X(0) = X_0,
\end{equation}
where $\theta > 0$, $\sigma > 0$, and $\mu(t)$ is a given function. Then:
\begin{enumerate}
    \item The solution is
    \begin{equation}
    \label{eq:Exp_OU_proc_solved}
        X_t = e^{-t \sigma_X + \sigma W_t} \left( X_0 + \int_0^t e^{s \sigma_X - \sigma W_s} \mu(s) ds \right), \qquad \sigma_X := \theta + \frac{\sigma^2}{2}.
    \end{equation}
    \item The mean is
    \begin{equation}
    \label{eq:Exp_OU_proc_mean}
        \mathbb{E}[X_t] = X_0 e^{-\theta t} + \int_0^t e^{-\theta (t-s)} \mu(s) ds.
    \end{equation}
    \item For $\mu(t) \equiv 0$, the covariance is
    \begin{equation}
    \label{eq:Exp_OU_proc_cov}
        \mathrm{Cov}(X_t, X_s) = X_0^2 e^{-\theta (t+s)} \left( e^{\sigma^2 \min(t,s)} - 1 \right)
    \end{equation}
    and the variance is
    \begin{equation}
    \label{eq:Exp_OU_proc_var}
        \mathrm{Var}(X_t) = X_0^2 e^{-2\theta t} \left( e^{\sigma^2 t} - 1 \right).
    \end{equation}
\end{enumerate}
\end{lem}
\begin{proof}
Let $Y_t = \exp\left(-\sigma_X t + \sigma W_t\right)$, where $\sigma_X = \theta + \frac{\sigma^2}{2}$. By Itô's lemma,
\[
dY_t = -\sigma_X Y_t dt + \sigma Y_t dW_t.
\]
Let $X_t = Y_t Z_t$. Applying Itô's product rule,
\[
dX_t = Y_t dZ_t + Z_t dY_t.
\]
Substitute $dY_t$ and set $dZ_t = \frac{\mu(t)}{Y_t} dt$ so that
\[
dX_t = Y_t \frac{\mu(t)}{Y_t} dt + Z_t \left( -\sigma_X Y_t dt + \sigma Y_t dW_t \right) = \mu(t) dt - \sigma_X X_t dt + \sigma X_t dW_t.
\]
This matches the original SDE. Integrating $dZ_t$ gives
\[
Z_t = X_0 + \int_0^t \frac{\mu(s)}{Y_s} ds.
\]
Therefore,
\[
X_t = Y_t Z_t = e^{-t \sigma_X + \sigma W_t} \left( X_0 + \int_0^t e^{s \sigma_X - \sigma W_s} \mu(s) ds \right).
\]

\noindent For the mean, note that $\mathbb{E}[e^{\sigma W_t - \frac{1}{2}\sigma^2 t}] = 1$ and $\mathbb{E}[e^{\sigma W_t}] = e^{\frac{1}{2}\sigma^2 t}$, so
\[
\mathbb{E}[X_t] = X_0 e^{-\theta t} + \int_0^t e^{-\theta (t-s)} \mu(s) ds.
\]

\noindent For the variance (when $\mu(t) \equiv 0$), $X_t = X_0 e^{-\theta t} e^{\sigma W_t - \frac{1}{2}\sigma^2 t}$, so
\[
\mathrm{Var}(X_t) = X_0^2 e^{-2\theta t} \left( \mathbb{E}[e^{2\sigma W_t - \sigma^2 t}] - 1 \right) = X_0^2 e^{-2\theta t} \left( e^{2\sigma^2 t - \sigma^2 t} - 1 \right) = X_0^2 e^{-2\theta t} (e^{\sigma^2 t} - 1).
\]
\end{proof}

\bibliographystyle{IEEEtran}
\bibliography{literature}

\begin{thebibliography}{10}
\providecommand{\url}[1]{#1}
\csname url@samestyle\endcsname
\providecommand{\newblock}{\relax}
\providecommand{\bibinfo}[2]{#2}
\providecommand{\BIBentrySTDinterwordspacing}{\spaceskip=0pt\relax}
\providecommand{\BIBentryALTinterwordstretchfactor}{4}
\providecommand{\BIBentryALTinterwordspacing}{\spaceskip=\fontdimen2\font plus
\BIBentryALTinterwordstretchfactor\fontdimen3\font minus
  \fontdimen4\font\relax}
\providecommand{\BIBforeignlanguage}[2]{{%
\expandafter\ifx\csname l@#1\endcsname\relax
\typeout{** WARNING: IEEEtran.bst: No hyphenation pattern has been}%
\typeout{** loaded for the language `#1'. Using the pattern for}%
\typeout{** the default language instead.}%
\else
\language=\csname l@#1\endcsname
\fi
#2}}
\providecommand{\BIBdecl}{\relax}
\BIBdecl

\bibitem{Feng2007Effect}
B.~Feng, A.~K. Yong, and H.~An, ``Effect of various factors on the particle
  size of calcium carbonate formed in a precipitation process,''
  \emph{Materials Science and Engineering: A}, vol. 445-446, pp. 170--179,
  2007.

\bibitem{Ma2002Optimal}
D.~L. Ma, D.~K. Tafti, and R.~D. Braatz, ``Optimal control and simulation of
  multidimensional crystallization processes,'' \emph{Computers \& Chemical
  Engineering}, vol.~26, no. 7-8, pp. 1103--1116, 2002.

\bibitem{Woodall2024}
C.~M. Woodall, K.~Vaz~Gomes, A.~Voigt, K.~Sundmacher, and J.~Wilcox,
  ``\BIBforeignlanguage{en}{Tuning acid extraction of magnesium and calcium
  from platinum group metal tailings for {CO$_{2}$} conversion and storage},''
  \emph{\BIBforeignlanguage{en}{RSC Sustain.}}, vol.~2, no.~11, pp. 3320--3333,
  2024.

\bibitem{Seifritz1990CO2}
W.~Seifritz, ``{CO2} disposal by means of silicates,'' \emph{Nature}, vol. 345,
  no. 6275, pp. 486--486, 1990.

\bibitem{Goff1998Environmental}
F.~Goff, W.~Seifritz, T.~McCormick, K.~Lackner, B.~Sass, J.~Smith, B.~Smith,
  L.~Smith, A.~Walder, S.~White, H.~Ziock, G.~Cherepakhin, R.~Rieke, R.~Smith,
  S.~Stegmaier, E.~Goff, K.~Goff, and M.~Goff, ``Environmental geosciences,''
  \emph{Environmental Geosciences}, vol.~5, no.~3, pp. 95--114, 1998.

\bibitem{Bobicki2012Carbon}
E.~R. Bobicki, Q.~Liu, Z.~Xu, and H.~Zeng, ``Carbon capture and storage using
  alkaline industrial wastes,'' \emph{Progress in Energy \& Combustion
  Science}, vol.~38, no.~2, pp. 302--320, 2012.

\bibitem{Olajire2013Carbon}
A.~A. Olajire, ``Carbon capture and storage options: a review,'' \emph{J PETROL
  SCI ENG}, vol. 109, pp. 302--324, 2013.

\bibitem{Sanna2014Carbonation}
A.~Sanna, Y.~Kuva, G.~Perotti, M.~J. Castaldi, and M.~M. Maroto-Valer,
  ``Carbonation of {Ca}- and {Mg}-silicates: Process analysis and
  optimisation,'' \emph{Chem. Soc. Rev.}, vol.~43, no.~23, pp. 7986--8008,
  2014.

\bibitem{Hegde2024Towards}
C.~Hegde, A.~Voigt, and K.~Sundmacher, ``Towards {pH} swing-based {CO2}
  mineralization by {Calcium Carbonate} precipitation: Modeling and
  experimental analysis,'' in \emph{Proceedings of the 34th European Symposium
  on Computer Aided Process Engineering / 15th International Symposium on
  Process Systems Engineering (ESCAPE34/PSE24), June 2-6, 2024, Florence,
  Italy}, F.~Manenti and G.~V. Reklaitis, Eds.\hskip 1em plus 0.5em minus
  0.4em\relax Elsevier B.V., 2024, pp. 1--6.

\bibitem{Koutsoukos1984Precipitation}
P.~G. Koutsoukos and C.~G. Kontoyannis, ``Precipitation of calcium carbonate in
  aqueous solutions,'' \emph{J. Chem. Soc. - Faraday Transactions I}, vol.~80,
  pp. 1181--1192, 1984.

\bibitem{Eisenschmidt2015Face}
H.~Eisenschmidt, A.~Voigt, and K.~Sundmacher, ``Face-specific growth and
  dissolution kinetics of potassium dihydrogen phosphate crystals from batch
  crystallization experiments,'' \emph{Cryst. Growth Des.}, vol.~15, no.~1, pp.
  219--227, 2015.

\bibitem{Borchert2012Efficient}
C.~Borchert and K.~Sundmacher, ``Efficient formulation of crystal shape
  evolution equations,'' \emph{Chemical Engineering Science}, vol.~84, pp.
  85--99, 2012.

\bibitem{Simone2017Systematic}
B.~Simone, G.~Di~Profio, E.~Curcio, E.~Drioli, F.~De~Luca, and P.~Parapari,
  ``Systematic study of the impact of agitation rate on {CaCO3} crystallization
  in a double-jet precipitation,'' \emph{Chemical Engineering Science}, vol.
  173, pp. 334--344, 2017.

\bibitem{Rauscher2005Influence}
F.~Rauscher, Schmaucher, and et~al., ``Influence of process parameters on the
  morphology of precipitated calcium carbonate,'' \emph{Colloids and Surfaces
  A: Physicochemical and Engineering Aspects}, vol. 254, no. 1-3, pp. 149--158,
  2005.

\bibitem{Falini2009Calcium}
G.~Falini, S.~Fermani, G.~Tosi, and E.~Dinelli, ``Calcium carbonate morphology
  and structure in the presence of seawater ions and humic acids,''
  \emph{Cryst. Growth Des.}, vol.~9, no.~5, pp. 2065--2072, 2009.

\bibitem{Park2008Effects}
W.~K. Park, S.~J. Ko, S.~W. Lee, K.~H. Cho, J.~W. Ahn, and C.~Han, ``Effects of
  magnesium chloride and organic additives on the synthesis of aragonite
  precipitated calcium carbonate,'' \emph{J. Cryst. Growth}, vol. 310, no.~10,
  pp. 2593--2601, 2008.

\bibitem{Davis2000Calcium}
K.~J. Davis, A.~K. Lee, and et~al., ``Calcium carbonate precipitation and the
  role of organic matrix,'' \emph{Science}, vol. 290, no. 5494, pp. 1134--1137,
  2000.

\bibitem{Bobicki2012}
E.~R. Bobicki, Q.~Liu, Z.~Xu, and H.~Zeng, ``Carbon capture and storage using
  alkaline industrial wastes,'' pp. 302--320, Apr. 2012.

\bibitem{Khosa2019}
A.~A. Khosa, T.~Xu, B.~Xia, J.~Yan, and C.~Zhao, ``Technological challenges and
  industrial applications of caco3/cao based thermal energy storage system –
  a review,'' \emph{Solar Energy}, vol. 193, pp. 618--636, Nov. 2019.

\bibitem{Gomes2025}
K.~V. Gomes, C.~M. Woodall, H.~Pilorgé, P.~Psarras, and J.~Wilcox,
  ``Techno-economic analysis of indirect carbonation processes for carbon
  sequestration using mining waste,'' \emph{Energy Advances}, 2025.

\bibitem{Ramkrishna2014Population}
D.~Ramkrishna and M.~R. Singh, ``Population balance modeling: Current status
  and future prospects,'' \emph{Annual Review of Chemical and Biomolecular
  Engineering}, vol.~5, no.~1, pp. 123--146, 2014.

\bibitem{Liendo2022Nucleation}
F.~Liendo, M.~Arduino, F.~A. Deorsola, and S.~Bensaid, ``Nucleation and growth
  kinetics of {CaCO3} crystals in the presence of foreign monovalent ions,''
  \emph{J. Cryst. Growth}, vol. 578, p. 126406, 2022.

\bibitem{John2009Numerical}
V.~John, B.~Niemann, and K.~Sundmacher, ``Numerical solution of a population
  balance system for crystallization,'' \emph{Chemical Engineering Science},
  vol.~64, no.~4, pp. 740--749, 2009.

\bibitem{Durr2020Approximate}
R.~Dürr and A.~Bück, ``Approximate moment methods for population balance
  equations in particulate and bioengineering processes,'' \emph{Processes},
  vol.~8, no.~4, p. 414, 2020.

\bibitem{Bartsch2019Multiscale}
C.~Bartsch, T.~Stolze, B.~Niemann, V.~John, and K.~Sundmacher, ``Multiscale
  modeling of crystallization processes: Coupling population balance equations
  with fluid dynamics,'' \emph{Chemical Engineering Science}, vol. 208, p.
  115160, 2019.

\bibitem{Chen1986Generalized}
C.-C. Chen, H.~I. Britt, J.~F. Boston, and L.~B. Evans, ``Generalized
  thermodynamic model for predicting solubility of salt in mixed-solvent
  systems,'' \emph{AIChE Journal}, vol.~32, no.~3, pp. 444--454, 1986.

\bibitem{Qian2012Calculation}
H.~Qian, X.~Zhang, and P.~Li, ``Calculation of caco3 solubility
  (precipitability) in natural waters,'' \emph{Asian Journal of Chemistry},
  vol.~24, pp. 668--672, 02 2012.

\bibitem{Teir2007Dissolution}
S.~Teir, S.~Eloneva, J.~Fagerlund, L.~Tervahauta, J.~Salminen, T.~Uusitalo,
  M.~Ilola, B.~Nyman, S.~Rasi, S.~Hilska, and M.~Vesterinen, ``Dissolution of
  natural serpentine for {CO2} sequestration,'' \emph{INT J MINER PROCESS},
  vol.~83, no. 1-2, pp. 36--44, 2007.

\bibitem{Teir2009Thermodynamic}
S.~Teir, Y.~Kuva, J.~Fagerlund, B.~Nyman, J.~Savolainen, J.~Salminen,
  L.~Tervahauta, S.~Eloneva, M.~Ilola, M.~Vesterinen, S.~Rasi, S.~Hilska, and
  T.~Uusitalo, ``Thermodynamic aspects of using magnesium silicate minerals as
  raw materials for carbon dioxide sequestration,'' \emph{Applied Energy},
  vol.~86, no.~2, pp. 262--271, 2009.

\bibitem{Rudy2017DataDriven}
S.~H. Rudy, S.~L. Brunton, J.~L. Proctor, and J.~N. Kutz, ``Data-driven
  discovery of partial differential equations,'' \emph{Science Advances},
  vol.~3, no.~4, 2017.

\bibitem{Raissi2018Deep}
M.~Raissi, P.~Perdikaris, and G.~E. Karniadakis, ``Deep hidden physics models:
  Deep learning of nonlinear partial differential equations,'' \emph{The
  Journal of Machine Learning Research}, vol.~19, no.~1, pp. 932--955, 2018.

\bibitem{E2021Deep}
W.~E, Z.~Ma, and C.~Zheng, ``Deep learning for high-dimensional partial
  differential equations,'' \emph{Nonlinearity}, vol.~35, no.~1, pp. R1--R48,
  2021.

\bibitem{Han2018Solving}
J.~Han, A.~Jentzen, and W.~E, ``Solving high-dimensional partial differential
  equations using deep learning,'' \emph{PNAS}, vol. 115, no.~34, pp.
  8505--8510, 2018.

\bibitem{HiremathECC2024}
M.~Kakanov, S.~A. Hiremath, A.~Voigt, K.~Sundmacher, and N.~Bajcinca,
  ``Learnable adaptive and robust controller for a two particle carbonate
  precipitation process,'' in \emph{2024 European Control Conference
  ({ECC})}.\hskip 1em plus 0.5em minus 0.4em\relax IEEE, jun 2024.

\bibitem{HiremathEscape2024}
S.~Hiremath, M.~Kakanov, A.~Voigt, K.~Sundmacher, and N.~Bajcinca, ``Learning
  based adaptive robust control of a precipitation process,'' in \emph{34th
  European Symposium on Computer Aided Process Engineering / 15th International
  Symposium on Process Systems Engineering}, ser. Computer Aided Chemical
  Engineering, F.~Manenti and G.~V. Reklaitis, Eds.\hskip 1em plus 0.5em minus
  0.4em\relax Elsevier, 2024, vol.~53, pp. 1801--1806.

\bibitem{Larson1973}
T.~E. Larson, F.~W. Sollo, F.~F. Mcgurk, and T.~E. Larson, ``Complexes
  affecting the solubility of calcium carbonate in water wrc research report
  no. 68 complexes affecting the solubility of calcium carbonate in water,''
  Tech. Rep., 1973.

\bibitem{Liendo2022}
F.~Liendo, M.~Arduino, F.~A. Deorsola, and S.~Bensaid, ``Nucleation and growth
  kinetics of caco3 crystals in the presence of foreign monovalent ions,''
  \emph{Journal of Crystal Growth}, vol. 578, Jan. 2022.

\bibitem{Lax1964}
P.~D. Lax and B.~Wendroff, ``\BIBforeignlanguage{en}{Difference schemes for
  hyperbolic equations with high order of accuracy},''
  \emph{\BIBforeignlanguage{en}{Commun. Pure Appl. Math.}}, vol.~17, no.~3, pp.
  381--398, aug 1964.

\bibitem{Gikhman2007}
I.~I. Gikhman and A.~V. Skorokhod, \emph{\BIBforeignlanguage{en}{The theory of
  stochastic processes {III}}}, ser. Classics in mathematics.\hskip 1em plus
  0.5em minus 0.4em\relax Berlin, Germany: Springer, mar 2007.

\bibitem{Zeller1956}
E.~J. Zeller and J.~L. Wray, ``Factors influencing precipitation of calcium
  carbonate,'' \emph{AAPG Bulletin}, vol.~40, pp. 140--152, Jan. 1956.

\bibitem{Chassagneux2017}
J.~F. Chassagneux, H.~Chotai, and M.~Muuls, \emph{\BIBforeignlanguage{en}{A
  forward-backward {SDEs} approach to pricing in carbon markets}}, 1st~ed.,
  ser. SpringerBriefs in Mathematics of Planet Earth.\hskip 1em plus 0.5em
  minus 0.4em\relax Cham, Switzerland: Springer International Publishing, oct
  2017.

\end{thebibliography}

\end{document}